\documentclass[11pt,twoside, a4paper]{article}
	
\usepackage{latexsym}
\usepackage[english]{babel}
\usepackage{t1enc}
\usepackage{mathrsfs}
\usepackage{ifthen}
\usepackage{url}
\usepackage{epsfig}
\usepackage{amsbsy}
\usepackage{extarrows}
\usepackage{amsfonts}
\usepackage{amsmath}
\usepackage{amssymb}
\usepackage{amsthm}
\usepackage{amsxtra}
\usepackage{mathabx}
\usepackage{bbm}
\usepackage{color}

\usepackage[hmargin=2.28cm,vmargin=2.26cm]{geometry}  

\usepackage{verbatim}

\usepackage{hyperref} 

\numberwithin{equation}{section}
\numberwithin{figure}{section}

\newtheoremstyle{thm-style-oskari}
{7pt}      
{7pt}      
{\itshape} 
{}         
{\scshape} 
{.}        
{.5em}     
{}         

\theoremstyle{thm-style-oskari}
    \newtheorem{theorem}{Theorem}[section]
    \newtheorem{proposition}[theorem]{Proposition}
    \newtheorem{corollary}[theorem]{Corollary}
    \newtheorem{lemma}[theorem]{Lemma}
    \newtheorem{definition}[theorem]{Definition}
    
    \newtheorem{convention}[theorem]{Convention}
    \newtheorem{remark}[theorem]{Remark}

\newenvironment{Proof}[1][Proof]{\begin{proof}[\sc{#1}]}{\end{proof}}

\newcommand{\NCorollary}[2] {
        \begin{corollary}[#1] \label{crl:#1}
                #2
        \end{corollary}
        }

\newcommand{\bels}[2] {
        \begin{equation} \label{#1} \begin{split} 
                #2 
        \end{split} \end{equation}
        }

\newcommand{\bea}[1]{
	\begin{align*}
		#1
	\end{align*}
	}

\definecolor{olivegreen}{rgb}{0,0.6,0.1}

\newcommand{\bs}[1]{\boldsymbol{\mathrm{#1}}} 
\newcommand{\bb}{\mathbb} 
\newcommand{\bbm}{\mathbbm} 
\renewcommand{\rm}{\mathrm} 
\renewcommand{\cal}{\mathcal} 
\newcommand{\scr}{\mathscr} 
\renewcommand{\frak}{\mathfrak} 
\newcommand{\ul}[1]{\underline{#1} \!\,} 
\newcommand{\ol}[1]{\overline{#1} \!\,} 
\newcommand{\wh}{\widehat}
\newcommand{\wt}{\widetilde}

\newcommand{\eps}{\varepsilon}
\newcommand{\ord} {\mathcal{O}}

\newcommand{\ceil}[1]  {\lceil  {#1} \rceil}

\renewcommand{\P}{\mathbb{P}}
\newcommand{\E}{\mathbb{E}}
\newcommand{\R}{\mathbb{R}}
\newcommand{\C}{\mathbb{C}}
\newcommand{\N}{\mathbb{N}}

\newcommand{\Cp}{\mathbb{H}}
\newcommand{\D}{\mathbb{D}}

\newcommand{\ee}{\mathrm{e}} 
\newcommand{\ii}{\mathrm{i}} 
\newcommand{\dd}{\mathrm{d}}

\newcommand{\p}[1]{({#1})}
\newcommand{\pb}[1]{\bigl({#1}\bigr)}
\newcommand{\pB}[1]{\Bigl({#1}\Bigr)}
\newcommand{\pbb}[1]{\biggl({#1}\biggr)}

\newcommand{\s}[1]{[{#1}]}
\renewcommand{\sb}[1]{\bigl[{#1}\bigr]}
\newcommand{\sB}[1]{\Bigl[{#1}\Bigr]}
\newcommand{\sbb}[1]{\biggl[{#1}\biggr]}

\renewcommand{\c}[1]{\{{#1}\}}
\newcommand{\cb}[1]{\bigl\{{#1}\bigr\}}
\newcommand{\cB}[1]{\Bigl\{{#1}\Bigr\}}

\newcommand{\sett}[1] { \{ {#1} \} }

\newcommand{\setb}[1] { \bigl\{ {#1} \bigl\} }
\newcommand{\setB}[1] { \Bigl\{ {#1} \Bigr\} }

\newcommand{\abs}[1]{\lvert #1 \rvert}
\newcommand{\absb}[1]{\big\lvert #1 \big\rvert}
\newcommand{\absB}[1]{\Big\lvert #1 \Big\rvert}
\newcommand{\absbb}[1]{\bigg\lvert #1 \bigg\rvert}

\newcommand{\norm}[1]{\lVert #1 \rVert}
\newcommand{\normb}[1]{\big\lVert #1 \big\rVert}
\newcommand{\normB}[1]{\Big\lVert #1 \Big\rVert}

\newcommand{\avg}[1]{\langle #1 \rangle}
\newcommand{\avgb}[1]{\big\langle #1 \big\rangle}

\newcommand{\scalar}[2]{\langle{#1} \mspace{2mu}, {#2}\rangle}

\DeclareMathOperator{\tr}{Tr}

\DeclareMathOperator{\supp}{supp}

\DeclareMathOperator{\re}{Re}
\DeclareMathOperator{\im}{Im}

\DeclareMathOperator{\dist} {dist}                
\DeclareMathOperator*{\spec}{Spec}						

\newcommand{\1} {\mspace{1 mu}}
\newcommand{\2} {\mspace{2 mu}}
\newcommand{\msp}[1] {\mspace{#1 mu}}

\newcommand{\tsfrac}[2]{{\textstyle \frac{\1 #1 \1}{\1 #2 \1}}} 

\newcommand{\titem}[1] {\item[\emph{(#1)}]} 

\begin{document}
\title{\bf Stability of the Matrix Dyson Equation and Random Matrices with Correlations
\vspace{0.3cm}}
\author{
\begin{minipage}{0.3\textwidth}
 \begin{center}
Oskari H. Ajanki\footnotemark[1] \\
\footnotesize 
{IST Austria}\\
{\url{oskari.ajanki@iki.fi}}
\end{center}
\end{minipage}
\begin{minipage}{0.3\textwidth}
\begin{center}
L\'aszl\'o Erd{\H o}s\footnotemark[1]  \\
\footnotesize {IST Austria}\\
{\url{lerdos@ist.ac.at}}
\end{center}
\end{minipage}
\begin{minipage}{0.3\textwidth}
 \begin{center}
Torben Kr\"uger\footnotemark[1] \\
\footnotesize 
{IST Austria}\\
{\url{torben.krueger@ist.ac.at}}
\end{center}
\end{minipage}
}

\date{\vspace{0.3cm}\today}
\maketitle
\thispagestyle{empty} 

\footnotetext[1]{Partially supported by ERC Advanced Grant RANMAT No.\ 338804}
\vspace{-0.5cm}
\begin{abstract} 
We consider real symmetric or complex hermitian random matrices 
with correlated entries. We prove
local laws for the resolvent  and  universality of the local eigenvalue
statistics in the bulk of the spectrum.
The correlations have  fast decay but are otherwise  of general form.
The key novelty is the detailed stability analysis of the corresponding  matrix valued Dyson equation
whose solution is  the deterministic  limit of the resolvent. 
\end{abstract}
\vspace{0.3cm}
{\bf AMS Subject Classification (2010):} \texttt{60B20}, \texttt{15B52}, \texttt{46T99}.

\section{Introduction}

E. Wigner's vision on the ubiquity of random matrix spectral statistics
in quantum systems posed a main challenge to mathematics. The basic conjecture
is that the distribution of the eigenvalue gaps of a large self-adjoint matrix  with sufficient disorder is universal in the sense that
it is independent of the details of the system and it depends only on the symmetry
type of the model. This universal statistics has been computed by Dyson, Gaudin and Mehta
for the Gaussian Unitary and Orthogonal Ensembles (GUE/GOE) in the limit as the dimension of the matrix
goes to infinity. GUE and GOE are the simplest
mean field random matrix models in their respective symmetry classes. They
have centered Gaussian entries that are identically distributed and independent (modulo the hermitian symmetry). 
The celebrated Wigner-Dyson-Mehta (WDM) universality conjecture, 
as formulated in the classical book of Mehta \cite{Mehta-book},
 asserts that the  same gap statistics holds if the matrix elements  are independent and  have 
 arbitrary   identical  distribution  (they are called {\it Wigner ensembles}).  
 The WDM conjecture has recently been proved in increasing generality
 in a series of papers \cite{EKYY2,EPRSY, ESY, ESYY} for both the real symmetric and complex hermitian 
 symmetry classes via the \emph{Dyson Brownian motion}.  An alternative approach introducing the \emph{four-moment comparison theorem} was presented in \cite{TaoVuActa2011,TaoVuWDM2011,TV}.
  In this paper we only discuss universality in the bulk of the spectrum, but we remark
  that a similar development took place for the edge universality.

The next step towards Wigner's vision is to drop the assumption of identical distribution 
in the WDM conjecture but still maintain the mean field character of the model
by requiring a uniform lower and upper bound on the variances of the
matrix elements.   This generalization has been achieved in two steps. If the matrix of variances 
is stochastic,  then universality was proved in  \cite{EKYY2, EYY, EYYrigi}, in parallel with the proof of the
original WDM conjecture for Wigner ensembles.  Without the stochasticity condition on the variances
the limiting eigenvalue density is not  the Wigner semicircle any more;
the correct density was analyzed in \cite{AEK1,AEK1cpam} and the universality was proved \cite{AEK2}. We remark that one may also depart from the semicircle law by adding a large diagonal component to Wigner matrices;
universality for such deformed Wigner matrices was obtained in \cite{LSSY2014}.
 Finally we mention  a separate direction 
to generalize the original WDM conjecture that aims at departing from the mean field condition:
bulk universality for general band matrices with a band width comparable to the matrix size 
was proved in \cite{BEYYband}, see also \cite{ShBand2} for Gaussian block-band matrices.

In this paper we drop the third key condition in the original WDM conjecture, 
the independence of the matrix
elements, i.e. we consider matrices with correlated entries.
Correlations come in many different forms and if they are extremely strong and long range, the
universality may even be violated. We therefore consider random matrix models with a suitable decay of correlations.
These models still carry sufficiently many random degrees of freedom for Wigner's vision to hold and, indeed,
our main result yields spectral universality for such matrices.

We now describe the key points of the current work. Our main result 
is the \emph{local law} for the resolvent
\vspace{-0.5cm}
\bels{intro: def of resolvent}{
\bs{G}(\zeta) := (\1\bs{H}-\zeta\1\bs{1})^{-1}
,
} 
of the random matrix $\bs{H}=\bs{H}^* \in \C^{N \times N}$ with the spectral parameter  $ \zeta$  in the complex upper half plane  $\Cp :=\{ \zeta\in \C\; : \; \im \zeta>0\}$    that lies 
very close to the real axis. We show that, as the size $N$ of the random matrix tends to infinity, $\bs{G}=\bs{G}(\zeta)$ is well approximated by a deterministic matrix $\bs{M}=\bs{M}(\zeta)$ that satisfies a nonlinear {\it matrix equation} of the form 
\bels{MDE in intro}{
\bs{1}+(\zeta \2\bs{1}-\bs{A}+\cal{S}[\bs{M}])\bs{M}\,=\,\bs{0}\,.
}
Here the self-adjoint matrix  $\bs{A}$ and the operator  $\cal{S}: \C^{N \times N} \to  \C^{N \times N}$  on the space of matrices are determined by the first two moments of the random matrix 
\bels{Intro definition of S}{
\bs{A}\,:= \, \E\2\bs{H}\,,\qquad \mathcal{S}[\1\bs{R}\1] \2:=\2  \E\,(\bs{H}-\bs{A}) \1\bs{R}(\bs{H}-\bs{A})\,.
}
The central role of \eqref{MDE in intro} 
 in the context of random matrices has been recognized   by several authors  \cite{Girko-book,Helton2007-OSE,PasturShcerbinaAMSbook,WegnerNorb}. We will call \eqref{MDE in intro} {\it Matrix Dyson Equation} (MDE) since the analogous equation for 
the resolvent is sometimes called Dyson equation in perturbation theory.

Local laws have become a cornerstone in the analysis of spectral properties of large random matrices \cite{AEK2,FixedDegree2016,EKYY,ESY3,EYYrigi,GNTT,YAR2016,TaoVuConcentration}. 
In its simplest form, a local law  considers the normalized trace $  \frac{1}{N} \tr \bs{G}(\zeta)$ of the resolvent. Viewed as a  Stieltjes transform, it describes the empirical density of eigenvalues on the scale determined by $\eta =\im\, \zeta$. Assuming a normalization such that the spectrum of $\bs{H}$ remains bounded as $N\to\infty$, the typical eigenvalue spacing in the bulk is of order $1/N$. 
The local law  asserts that this normalized trace  approaches
a deterministic function $m(\zeta)$ as the size $N$ of the matrix tends to infinity
and this convergence holds uniformly even if $\eta =\eta_N$ depends on $N$ as long as $\eta\gg 1/N$.
Equivalently, the empirical density of the  eigenvalues converges on any scales slightly above $1/N$ to a  deterministic limit measure on $\R$ with Stieltjes transform $m(\zeta)$.

Since $\bs{G}$ is asymptotically close to $\bs{M}$, the deterministic limit of the Stieltjes transform of the empirical spectral measure is given by $m(\zeta)= \frac{1}{N}\tr \bs{M}(\zeta)$. Already in the case of random matrices with centered independent entries (\emph{Wigner-type} matrices) the limiting measure $\rho(\rm{d} \omega)$ and its Stieltjes transform $m(\zeta)$ typically depend
 on the entire matrix of variances $s_{xy}:= \E |h_{xy}|^2$ and the only known way to determine $\rho$ is to solve \eqref{MDE in intro}. However, in this setting the problem simplifies considerably because the off-diagonal elements of $\bs{G}$ tend to zero, $\bs{M}$ is a diagonal matrix and \eqref{MDE in intro} reduces to  a vector equation for  its diagonal elements. 
In case the variance matrix is doubly stochastic, $\sum_y s_{xy}=1$ (\emph{generalized Wigner matrix}), the problem simplifies yet again, leading to $\bs{M} = m_{\rm{sc}}\bs{1}$, where $m_{\rm{sc}}=m_{\rm{sc}}(\zeta)$ is the Stieltjes transform of the celebrated semicircle law. 

The main novelty of this work is  to handle  general correlations that do not allow to simplify \eqref{MDE in intro}.
The off-diagonal matrix elements $G_{x y}$,  $x\neq y$, do not vanish in general, even in the $N\to\infty$ limit.  
The proof of the local law consists of two major  parts.   
First,  we derive an approximate equation
\vspace{-0.5cm}
\bels{perturbed MDE in intro}{
\bs{1}+(\zeta \2\bs{1}-\bs{A}+\cal{S}[\bs{G}])\bs{G}\,\approx\,\bs{0}\,,
}
for the resolvent of $\bs{H}$.  
 To avoid confusion we stress that the expectation over the random matrix $\bs{H}$ in \eqref{Intro definition of S} is only used to define the deterministic operator $\cal{S}$. If the argument $\bs{G}$ of $\cal{S}$ itself is random as in \eqref{perturbed MDE in intro}, then $\cal{S}[\bs{G}]$ is still random and 
we have $\cal{S}[\bs{G}]=\wt{\E}\2(\wt{\bs{H}}-\bs{A})\bs{G}(\wt{\bs{H}}-\bs{A})$, where the expectation $\wt{\E} $ acts only on an independent copy $\wt{\bs{H}}$ of $\bs{H}$.

Second, we show that the Matrix Dyson Equation  \eqref{MDE in intro} is stable under small perturbations, concluding that $\bs{G} \approx \bs{M}$.  The nontrivial correlations and the non commutativity of the matrix structure in the Dyson equation pose major difficulties compared to  the uncorrelated case.  
 
Local laws are the first step of a general 
{\it three step strategy} developed in \cite{ESY, ESYY,EYY, EYYrigi}
for proving universality.   The second step is to  add
a tiny independent Gaussian component and prove universality for this 
slightly deformed model via analyzing the fast convergence of the
Dyson Brownian motion (DBM) to local equilibrium. Finally, the third step is
a perturbation argument showing that the tiny Gaussian component
does not alter the local statistics.

In fact, the second and the third steps are very robust arguments and they easily extend to the correlated case. They do not use any properties of the original ensemble other  than the a priori bounds encoded in  the local laws, provided that the variances of the matrix elements have a positive lower bound (see \cite{BY13EigenMoment, ES, LSY,LY}). 
Therefore our work  focuses on the first step,  establishing the stability of \eqref{MDE in intro} and thus obtaining a  local law. 

Prior to the current paper, bulk universality 
has already been established for several random matrix models which carry some 
specific correlation from their construction. These include  sample covariance
matrices \cite{ESYY}, adjacency matrices of large regular graphs  \cite{BHKY}
and  invariant  $\beta$-ensembles at various levels of generality \cite{BFG,BEY,DeiftRiemannHilbert,DeiftInvariant,FGmulti,GVrep,MSchBetaUni}.
However, neither of these papers aimed at understanding the effect of a general correlation
nor were their methods  suitable to deal with it. 
Universality for Gaussian matrices with  a {\it translation invariant}
covariance structure was established in \cite{AEK3}.
For general distributions  of the matrix entries,  but with a specific two-scale finite range
correlation structure  that is smooth on the large scale and translation invariant on the short scale,  universality was proved in \cite{Che2016}, independently of the current work.

Finally,  we mention that  there exists an extensive literature on the limiting eigenvalue distribution for random matrices with correlated entries on the global scale
(see e.g. \cite{AZdep,BMP2015,BKV1996,Girko-book,HLN2005,Schenker2005} 
and references therein), however these works either dealt with Gaussian random matrices or more  specific correlation
 structures that allow   one  to effectively reduce \eqref{MDE in intro} to a vector or scalar equation.
While the Matrix Dyson Equation in full generality was introduced for the analysis on the global scale before us, we are not aware of
a  proof establishing that  the empirical density of states converges to the deterministic density given by the solution of the MDE for a similarly broad class of models that we consider in this paper. 
This convergence is expressed by the fact that $\frac{1}{N} \tr \bs{G}(\zeta) \approx \frac{1}{N} \tr \bs{M}(\zeta)$ holds for any fixed $\zeta\in \Cp$. 
We thus believe that our proof identifying the limiting eigenvalue distribution is a new result even on the global scale for ensembles with general short range correlations and non-Gaussian distribution.  

  We present the stability of the MDE and its application to random matrices with correlated entries separately. Our findings on the MDE are given in Section~\ref{subsec:The Matrix Dyson equation}, while Section~\ref{subsec:Correlated random matrices} contains the results about random matrices with correlated entries. These sections can be read independently of each other,  with the latter relying on the former only through some basic definitions. 
 In Section~\ref{sec:Local law for correlated random matrices} we prove the local law for  random matrices with correlations. The proof relies on the results  stated  in Section~\ref{subsec:The Matrix Dyson equation}. These results concerning the MDE are established in Section~\ref{sec:The Matrix Dyson Equation}, which can be read independently of any other section.
 Besides the results from Section~\ref{subsec:The Matrix Dyson equation} on the MDE  the  main ingredients of the proof  in Section~\ref{sec:Local law for correlated random matrices}
 are (i) estimates on the random error term appearing in the approximate MDE \eqref{perturbed MDE in intro} and 
 (ii) the  \emph{fluctuation averaging mechanism} 
 for this error term. These two inputs  (Lemma~\ref{lmm:Smallness of error matrix} and Proposition~\ref{prp:fluctuation averaging})  are established in 
 Sections~\ref{sec:Estimating the error term} and~\ref{sec:Fluctuation averaging}, respectively.  However,  Section~\ref{sec:Local law for correlated random matrices} can be understood without reading the ensuing sections, taking these inputs for granted. 
  Finally, we apply the local law to establish  the rigidity of eigenvalues and bulk universality in Section~\ref{sec:Bulk universality and rigidity}.  The appendix collects the proofs for auxiliary results of generic nature that are not directly concerned with either the MDE or the random matrices.

\newpage

\section{Main results}
\label{sec:Main results}

\subsection{The Matrix Dyson Equation}
\label{subsec:The Matrix Dyson equation}
In this section we present our main results on the \emph{Matrix Dyson Equation} and its stability. The corresponding proofs are carried out in Section~\ref{sec:The Matrix Dyson Equation}.
 We consider the linear space $ \C^{N \times N}$ of $N\times N$ complex matrices $\bs{R}=(r_{xy})_{x,y =1}^N$, 
and make it a Hilbert space by  equipping it with the standard normalized scalar product
\bels{matrix scalar product}{
\scalar{\1\bs{R}}{\bs{T}\1}:= \frac{1}{N}\msp{-2}\tr \bs{R}^{\msp{-2}*}\bs{T}\,.
}
We denote the cone of strictly positive definite matrices by
\[
\mathscr{C}_+\,:=\,\{\2\bs{R}\in \C^{N\times N}:\,  \bs{R}\1>\1 \bs{0}\2\}\,,
\]
and by $\overline{\mathscr{C}}_+$ its closure, the cone of positive semidefinite matrices.

Let $\bs{A}=\bs{A}^{\msp{-2}*} \in \C^{N \times N}$ be a self-adjoint matrix. We will refer to $\bs{A}$ as the \emph{bare matrix}. Furthermore, let $\cal{S}:\C^{N \times N} \to \C^{N \times N}$ be a linear operator 
 that is
\begin{itemize}
\item 
\emph{self-adjoint} w.r.t. the scalar product \eqref{matrix scalar product}: $\tr\bs{R}^{\msp{-2}*}\cal{S}[\bs{T}]=\tr\cal{S}[\bs{R}]^*\bs{T}$ for any $\bs{R},\bs{T} \in \C^{N \times N}$; 
\item 
\emph{positivity preserving}: $\cal{S}[\bs{R}]\ge \bs{0}$ for any $\bs{R}\ge \bs{0}\2$.
\end{itemize}
Note that in particular $\cal{S}$ commutes with taking the adjoint,  $\cal{S}[\bs{R}]^*=\cal{S}[\1\bs{R}^{\msp{-1}*}]$, and hence it is real symmetric, $\tr\bs{R}\2\cal{S}[\bs{T}]=\tr\cal{S}[\bs{R}]\bs{T}$, for all $\bs{R},\bs{T}\in \C^{N \times N}$.  We will refer to $\cal{S}$ as the \emph{self-energy operator}. 

We call a pair $(\bs{A}, \cal{S})$ consisting of a bare matrix and a self-energy operator with the properties above a \emph{data pair}.
For a given {data pair} $(\bs{A}, \cal{S})$ and a \emph{spectral parameter} $\zeta \in \Cp$
 in the upper half plane  we consider the associated \emph{Matrix Dyson Equation (MDE)},
\bels{MDE}{
-\2\bs{M}(\1\zeta\1)^{-1} \,=\,  \1\zeta\1\bs{1}-\bs{A}+ \cal{S}\1[\1\bs{M}(\1\zeta\1)] \,, 
}
for a solution matrix $\bs{M}=\bs{M}(\1\zeta\1) \in \C^{N\times N}$ with positive definite imaginary part, 
\bels{im M positive definite}{
\im \bs{M} :=
\frac{1}{2\ii}(\bs{M}-\bs{M}^*)
\in \scr{C}_+
\,.
}

The question of existence and uniqueness of solutions to \eqref{MDE} with the constraint \eqref{im M positive definite} has been answered in \cite{Helton2007-OSE}. 
The MDE has a unique solution matrix $\bs{M}(\1\zeta\1)$ for any spectral parameter $\zeta \in\Cp$ and these matrices constitute a holomorphic function $\bs{M}: \Cp \to \C^{N \times N}$. 

On the space of matrices $\C^{N \times N}$  we 
consider three norms.  For $\bs{R}\in \C^{N \times N}$ we denote by $\norm{\bs{R}}$ the operator norm induced by the standard Euclidean norm $\norm{\2\cdot\2}$ on $\C^N$, by $\norm{\bs{R}}_{\rm{hs}}:=\sqrt{\scalar{\bs{R}}{\bs{R}}}$ the norm associated with the scalar product \eqref{matrix scalar product} and by
\bels{definition entrywise max norm}{
\norm{\bs{R}}_{\rm{max}} :=\max_{\,x,\1y=1}^N\2\abs{\1r_{xy}}
\,,
}
the entrywise maximum norm on $\C^{N \times N}$. We also denote the normalized trace of $\bs{R}$ by $\avg{\bs{R}}:=\scalar{\bs{1}}{\bs{R}}$. 

For linear operators $\cal{T}:\C^{N \times N} \to \C^{N \times N}$ we denote by $\norm{\cal{T}}$ the operator norm induced by the norm $\norm{\2\cdot\2}$ on $\C^{N \times N}$ and by $\norm{\cal{T}}_{\rm{sp}}$ the operator norm induced by $\norm{\2\cdot\2}_{\rm{hs}}$.

The following proposition provides a representation of the solution $\bs{M}$ as the Stieltjes-transform of a measure with values in $\overline{\scr{C}}_+$. This is a standard result for matrix-valued Nevanlinna functions (see e.g. \cite{GT2000}). For the  convenience  of the reader we provide a proof which also 
 gives an effective  control on the boundedness of the  support of this matrix-valued measure. 

\begin{proposition}[Stieltjes transform representation]
\label{prp:Stieltjes transform representation}
 Let $\bs{M}: \Cp \to \C^{N \times N}$ be the unique solution of \eqref{MDE} with $\im \bs{M}\in \scr{C}_+$. Then $\bs{M}$ admits a  Stieltjes transform representation,
\bels{M as Stieltjes transform}{
m_{xy}(\zeta)\,=\,\int_\R\frac{v_{xy}(\dd \tau)}{\tau-\zeta}\,, \qquad \zeta \in \Cp
\,,\;x, y=1,\dots,N
\,.
}
The measure $\bs{V}(\dd \tau)=(v_{x y}(\dd \tau))_{x,y=1}^N$ on the real line with values in positive semidefinite matrices is unique. It satisfies the normalization 
$\bs{V}(\R)=\bs{1}$ and has support in the interval $[-\kappa,\kappa]$, where
\bels{kappa definition}{
\kappa\,:=\, \norm{\bs{A}}+2\1\norm{\cal{S}}^{1/2}\,.
}
\end{proposition}

We will now make additional quantitative assumptions on the data pair $(\bs{A},\cal{S})$ that ensure a certain regularity of the measure $\bs{V}(\dd \tau)$. Our assumptions, labeled {\bf A1} and {\bf A2}, always come together with a set of \emph{model parameters} $\scr{P}_1$ and $\scr{P}_2$, respectively,
 that control them effectively.  Estimates will typically be uniform in all data pairs that satisfy these assumptions with the given set of model parameters. In particular, they are uniform in the size $N$ of the matrix, which is of great importance in the application to random matrix theory.

\begin{itemize}
\item[{\bf A1}] \emph{Flatness}:  Let $\scr{P}_1=(p_1,P_1)$ with $p_1,P_1>0$. The self-energy operator $\cal{S}$ is called \emph{flat} (with model parameters $\scr{P}_1$) if  it satisfies the lower and upper bound 
\bels{Flatness}{
p_1\2\avg{\bs{R}} \,\bs{1}\,\le\, \cal{S}[\bs{R}]\,\le\, P_1\2\avg{\bs{R}} \,\bs{1}\,,\qquad \bs{R} \in \overline{\scr{C}}_+\,.
}
\end{itemize} 

\begin{proposition}[Regularity of  self-consistent  density of states]
\label{prp:Density of states} 
Assume that $ \cal{S}$ is flat, i.e. it satisfies {\bf A1} with some model parameters $\scr{P}_1$ and that the bare matrix has a bounded spectral norm,
\bels{spectral norm bound on A}{
\norm{\bs{A}}\,\leq\, P_0
\,,
}
for some constant $P_0>0$. 
Then the holomorphic function $\avg{\bs{M}}: \Cp \to \Cp$ is the Stieltjes transform of a H\"older-continuous probability density $\rho$ with respect to the Lebesgue-measure, 
\bels{definition of rho through V}{
\avg{\bs{V}(\dd \tau)}\,=\, \rho(\tau)\2\dd \tau\,.
}
More precisely, 
\[
\abs{\1\rho(\tau_1)-\rho(\tau_2)}\,\leq\, C\2\abs{\tau_1-\tau_2}^{c}\,,\qquad \tau_1,\tau_2 \in \R\,,
\]
where $c>0$ is a universal constant and the constant $C>0$ depends only on the model parameters $\scr{P}_1$ and $P_0$.  Furthermore, $\rho$ is real analytic on the open set $\sett{\tau \in \R : \rho(\tau)>0}$. 
\end{proposition}

\begin{definition}[Self-consistent density of states]
\label{def:Density of states}
Assuming a flat self-energy operator, the probability density $\rho: \R \to [0,\infty)$, defined through \eqref{definition of rho through V}, is called the \emph{self-consistent density of states} (of the MDE with data pair $(\bs{A},\cal{S})$). We denote by $\supp \rho \subseteq \R$ its support on the real line and call it the \emph{self-consistent spectrum}  . With a slight abuse of notation we also denote by
\bels{definition extended DOS}{
\rho(\1\zeta\1)\,:=\, \frac{1}{\pi} \im\2 \avg{\bs{M}(\1\zeta\1)}\,,\qquad \zeta \in \Cp\,,
} 
the harmonic extension of $\rho$ to the complex upper half plane.
\end{definition}

The second set of assumptions describe the decay properties of the data pair $(\bs{A},\cal{S})$. To formulate them, we need to equip the index set $\{1,\dots,N\}$ with a concept of distance.
Recall that a \emph{pseudometric} $ d $ on a set $A$ is a symmetric function $ d : A \times A \to [0,\infty] $ such that $ d(x,y) \leq d(x,z) + d(z,y) $ for all $ x,y,z \in A$.
We say that the pseudometric space $ (A,d) $ with a finite set $A$  has \emph{sub-$ P$-dimensional volume}, for some constant $ P>0 $, if the metric balls $ B_\tau(x):=\{y\,:\, d(x,y)\2\leq\2 \tau\} $, satisfy
\bels{polynomial ball growth}{  
|B_\tau(x)|\,\leq\, \tau^P\,,
\qquad \tau \ge 2\;,\,x\in A
\,.
}

\begin{itemize}
\item[{\bf A2}] \emph{Faster than power law decay:} 
Let $\scr{P}_2=(P,\ul{\pi}_1,\ul{\pi}_2)$, where $P>0$ is a constant and $\ul{\pi}_k=(\pi_k(\nu))_{\nu=0}^\infty$, $k=1,2$ are  sequences of positive constants. The data pair $(\bs{A},\cal{S})$ is said to have \emph{faster than power law decay} (with model parameters $\scr{P}_2$) if there exists a pseudometric $d$ on the index space $ \{1,\dots,N\} $ such that the pseudometric space $\bb{X}=(\{1,\dots,N\},d)$ has {sub-$ P$-dimensional volume}  (cf. \eqref{polynomial ball growth}) and 
\begin{align}
\label{Expectation decay}
\abs{a_{xy}} 
\,&\leq\, \frac{\pi_1(\nu)}{(1+d(x,y))^\nu}+\frac{\pi_1(0)}{N}
\\
\label{S operator decay}
\qquad\abs{\1\cal{S}[\bs{R}]_{xy}} 
\,&\leq\, 
\biggl(\frac{\pi_2(\nu)}{(1+d(x,y))^\nu}+\frac{\pi_2(0)}{N}\,\biggr)\,\norm{\bs{R}}_{\rm{max}}
\,,\qquad
\bs{R} \in \C^{N\times N},
\end{align}
holds for any $\nu\in \N$ and $x,y \in \bb{X}$. 
\end{itemize}

In order to state bounds of the form \eqref{Expectation decay} and \eqref{S operator decay} more conveniently we introduce the following matrix norms.

\begin{definition}[Faster than power law decay] 
Given a pseudometric $ d $ on $\{1,\dots,N\}$ and a sequence $\ul{\pi} = (\pi(\nu))_{\nu=0}^\infty $  of positive constants, we define:
\bels{def of decay norm}{
\norm{\bs{R}}_{\ul{\pi}}
\,:=\, 
\sup_{\nu \in \N}\,\max_{x,y =1}^N \,\biggl(\frac{\pi(\nu)}{(1+d(x,y))^\nu}+\frac{\pi(0)}{N}\biggr)^{\!-1}\abs{\1r_{xy}}
\,,\qquad 
\bs{R} \in \C^{N \times N}
\,.
}
If $\norm{\bs{R}}_{\ul{\pi}} \leq 1 $, for some sequence $ \ul{\pi} $, we say that $ \bs{R} $ has \emph{faster than power law decay}  (up to level $\frac{1}{N}$)  in the pseudometric space $
\bb{X} := (\sett{1,\dots,N},d\2)  $.
\end{definition}

This norm expresses the typical behavior of many matrices in this paper that they have an off-diagonal decay faster than any power, up to a possible mean-field term of order $1/N$.
Using this norm the bounds \eqref{Expectation decay} and \eqref{S operator decay} take the simple forms:
\[ 
\qquad\norm{\bs{A}}_{\ul{\pi}_1} \leq 1\,,
\qquad
\norm{\1\cal{S}[\bs{R}]}_{\ul{\pi}_2} \leq \norm{\bs{R}}_{\rm{max}} 
\,.
\]

Our main result, the stability of the MDE, holds uniformly for all spectral parameters that are either away from the  self-consistent spectrum, $\supp \rho$,  or where the  self-consistent  density of states takes positive values. 
Therefore, for any $\delta>0$ we set 
\[
\D_\delta\,:=\, \setb{\2\zeta \in \Cp: \rho( \zeta)+\dist(\zeta, \supp \rho)> \delta\,}
\,.
\]

\begin{theorem}[Faster than power law decay of solution]
\label{thr:Arbitrarily high polynomial decay of solution}
Assume {\bf A1} and {\bf A2} and let $\delta>0$. Then there exists a positive sequence $\ul{\gamma}$ such that
\bels{Arbitrarily high polynomial decay of solution}{
\norm{\1\bs{M}(\1\zeta\1)}_{\ul{\gamma}}\,\le\, 1\,, \qquad \zeta \in \D_\delta\,.
}
The  sequence $\ul{\gamma}$ depends only on $\delta$ and the model parameters $\scr{P}_1$ and  $\scr{P}_2$.
\end{theorem}

Our main result on the MDE is its stability 
with respect to the entrywise maximum norm on $\C^{N \times N}$, see \eqref{definition entrywise max norm}.
The choice of this norm is especially useful for applications in random matrix theory, since the matrix valued error terms are typically controlled in this norm. We denote by 
\[
B^{\rm{max}}_\tau(\bs{R})\,:=\,\setb{\bs{Q} \in \C^{N \times N}: \norm{\bs{Q}-\bs{R}}_{\rm{max}} \le \tau\,}
\,,
\]
the ball of radius $ \tau > 0 $ around $\bs{R} \in \C^{N \times N}$ w.r.t. the entrywise maximum norm.

\begin{theorem}[Stability]
\label{thr:Stability}
Assume {\bf A1} and {\bf A2}, 
let $\delta>0$ and $\zeta \in \D_\delta$.  Then there exist constants $c_1,c_2>0$ and a unique function $\bs{\frak{G}}=\bs{\frak{G}}_\zeta: B^{\rm{max}}_{c_1}(\bs{0}) \to B^{\rm{max}}_{c_2}(\bs{M})$ such that 
\bels{perturbed MDE2}{
-\bs{1}\,=\, (\zeta\2\bs{1}-\bs{A}+\cal{S}[\bs{\frak{G}}(\bs{D})])\bs{\frak{G}}(\bs{D})+\bs{D}\,,
}
for all $\bs{D} \in B^{\rm{max}}_{c_1}(\bs{0})$, where $\bs{M}=\bs{M}(\1\zeta\1)$. The function $\bs{\frak{G}}$ is analytic. 
In particular, there exists a constant $C>0$ such that 
\bels{MDE stability2}{
\norm{\bs{\frak{G}}(\bs{D}_1)-\bs{\frak{G}}(\bs{D}_2)}_\rm{max}\le\, C\2\norm{\bs{D}_1-\bs{D}_2}_\rm{max}
\,.
}
for all $\bs{D}_1,\bs{D}_2 \in B^{\rm{max}}_{c_1/2}(\bs{0})$. 

 Furthermore, there is a sequence $\ul{\gamma}$ of positive constants, and a linear operator $\cal{Z}=\cal{Z}_\zeta: \C^{N \times N} \to \C^{N \times N}$ such that 
the derivative of $\bs{\frak{G}}$, evaluated at $\bs{D}=\bs{0}$, has the form
\bels{derivative representation with Z}{
\nabla \bs{\frak{G}}(\bs{0})\,=\, \cal{Z}+\bs{M}\,\rm{Id}\,,
}
and $\cal{Z}$, as well as its adjoint $\cal{Z}^*$ with respect to the scalar product \eqref{matrix scalar product}, satisfy
\bels{bounds on Z}{
\qquad \norm{\cal{Z}[\bs{R}]}_{\ul{\gamma}}+\norm{\cal{Z}^*[\bs{R}]}_{\ul{\gamma}}\,\le\, \norm{\bs{R}}_{\rm{max}}
\,,
\qquad\forall\2\bs{R} \in \C^{N \times N}
\,.
}
Here $c_1, c_2, C $ and $\ul{\gamma}$ 
depend only on $\delta$ and the model parameters $\scr{P}_1$, $\scr{P}_2$ from assumptions {\bf A1} and {\bf A2}.
\end{theorem}
Theorem~\ref{thr:Stability} states quantitative regularity properties of the analytic map $\bs{\frak{G}}$. These estimates yield strong stability properties of the MDE. 
For a concrete application in the proof of the local law, see Corollary~\ref{crl:Stability for local laws} below.

\subsection{Random matrices with correlations}
\label{subsec:Correlated random matrices}

In this section we present our results on local eigenvalue statistics of random matrices with correlations. 
Let $\bs{H}=(h_{x,y})_{x,y=1}^N \in \C^{N \times N}$ be a self-adjoint random matrix.
For a spectral parameter $\zeta \in \Cp$ we consider the associated \emph{Matrix Dyson Equation (MDE)},
\bels{RM MDE}{
-\2\bs{M}(\1\zeta\1)^{-1}\,&=\,\zeta\2\bs{1}-\bs{A} +\cal{S}[\bs{M}(\1\zeta\1)]
\,, 
\\
\\
\qquad\bs{A}\,:=\, \E \2\bs{H} \,,\qquad& \cal{S}[\bs{R}]\,:=\, \E\2 \bs{H} \bs{R} \bs{H}-\bs{A} \bs{R} \bs{A}
\,,
}
for a solution matrix $\bs{M}(\1\zeta\1)$ with positive definite imaginary part (cf. \eqref{im M positive definite}). 
The linear \emph{self-energy operator}  $\cal{S}: \C^{N \times N} \to  \C^{N \times N}$ from \eqref{RM MDE} preserves the cone $\ol{\scr{C}}_+$ of positive semidefinite matrices and the MDE therefore has a unique solution \cite{Helton2007-OSE} whose properties  have been presented  in Subsection~\ref{subsec:The Matrix Dyson equation}. 

Our main result states that under natural assumptions on the correlations of the entries within the random matrix $\bs{H}$, the resolvent 
\bels{definition resolvent}{
\bs{G}(\1\zeta\1)\,:=\, (\1\bs{H}- \zeta \2\bs{1})^{-1},
}
is close to the non-random solution $\bs{M}(\1\zeta\1)$ of the MDE \eqref{RM MDE}, provided $N$ is large enough. In order to list these assumptions, we write $\bs{H}$ as a sum of its expectation and fluctuation
\bels{H structure}{
\bs{H} \,=\, \bs{A} + \frac{1}{\sqrt{N}}\bs{W}\,. 
}
Here, the \emph{bare matrix} $\bs{A}$ is a non-random self-adjoint matrix and $\bs{W}$ is a self-adjoint random matrix with centered entries, $\E\2\bs{W}=\bs{0}$. The normalization factor $ N^{-1/2} $  in \eqref{H structure} ensures that the spectrum of the \emph{fluctuation matrix} $\bs{W}$, with entries of a typical size of order one, remains bounded.

In the following we will assume that there exists some pseudometric $ d $ on the index set $\{1, \dots,N\}$, such that the resulting pseudometric space
\[
\bb{X}\,=\, (\{1, \dots,N\},d\2)\,,
\]
has \emph{sub-$ P$-dimensional volume} for some constant $P>0$, i.e. $d$  satisfies \eqref{polynomial ball growth}, and that the bare and fluctuation matrices satisfy the following assumptions: 
\begin{itemize}
\item[{\bf B1}]  \emph{Existence of moments:} Moments of all orders of $\bs{W}$ exist, i.e.,  there is a sequence of positive constants $\ul{\kappa}_1=(\kappa_1(\nu))_{\nu \in \N}$ such that 
\bels{bounded moments}{
\qquad
\E\,\abs{w_{xy}}^\nu\,\le\, \kappa_1(\nu)
\,,
}
for all $x,y\in \bb{X}$ and $\nu \in \N$.
\item[{\bf B2}] \emph{Decay of expectation:}  The entries $a_{xy}$ of the bare matrix $\bs{A}$ decay in the distance of the indices $x$ and $y$, i.e., there is a sequence of positive constants $\ul{\kappa}_2=(\kappa_2(\nu))_{\nu \in \N}$ such that 
\bels{decay of expectation}{
\abs{\1a_{xy}}\,\le\, \frac{\kappa_2(\nu)}{(1+d(x,y))^\nu}\,,
}
for all $x,y \in \bb{X}$ and $\nu \in \N$.
\item[{\bf B3}]\emph{Decay of correlation:}  
The correlations in $ \bs{W} $ are fast decaying, i.e., there is a sequence of positive constants $\ul{\kappa}_3=(\kappa_3(\nu))_{\nu \in \N}$ such that for all symmetric sets $A,B \subseteq \bb{X}^2$ ($A$ is symmetric if $(x,y) \in A$ implies $(y,x) \in A$), and all smooth functions $\phi: \C^{\abs{A}} \to \R$ and $\psi: \C^{\abs{B}} \to \R$, we have
\bels{decay of correlation}{
\abs{\2\mathrm{Cov}(\1\phi(\rm{W}_{\!A})\1,\psi(\rm{W}_{\!B}))\1}
\,\leq\, 
\kappa_3(\nu)\frac{\norm{\nabla \phi}_\infty\norm{\nabla \psi}_\infty}{(1+d_2(A,B))^{\1\nu}} 
\,,
\qquad
\forall\2\nu \in\N.
}
Here, $ \mathrm{Cov}(Z_1,Z_2):=\2\E \2Z_1Z_2-\E\2Z_1\,\E\2Z_2$ is the covariance, $ \rm{W}_{\!A} := (w_{xy})_{(x,y)\in A} $, and
\[
d_2(A,B) \,:= \min\cb{\max\cb{d(x_1,x_2),d(y_1,y_2)}:(x_1,\1y_1) \2\in\1 A\,,\;  (x_2,\1y_2) \2\in\1 B}\,,
\]
is the distance between $A$ and $B$ in the product metric on $\bb{X}$. The supremum norm on vector valued functions $\Phi=(\phi_i)_i$ is  $\norm{\Phi}_\infty:=\sup_{Y}\max_{i}| \phi_i(Y)|$.

\item[{\bf B4}]  \emph{Flatness:}   There is a positive constant $\kappa_4$ such that for any two deterministic vectors $\bs{u},\bs{v} \in \C^N$ we have 
\bels{RM Flatness}{
\E\,\abs{ \bs{u}^*\2\bs{W}\bs{v}}^2\,\geq\, \kappa_4\, \norm{\bs{u}}^2\norm{\bs{v}}^2,
}
where $\norm{\2\cdot\2}$ denotes the standard Euclidean norm on $\C^N$.
\end{itemize}

We consider the constants 
\bels{definition of model parameters K}{
\scr{K}\,:=\,(P,\ul{\kappa}_1,\ul{\kappa}_2, \ul{\kappa}_3,\kappa_4)\,,
}
appearing in the above assumptions \eqref{polynomial ball growth} and \eqref{bounded moments}-\eqref{RM Flatness}, as \emph{model parameters}. These parameters are regarded as fixed and our statements are  uniform in the  ensemble of  all correlated random matrices of all dimensions $N$  satisfying {\bf B1-B4}  with given $\scr{K}$. 

Under the assumptions {\bf B1}-{\bf B4} the function $\rho: \Cp \to [0,\infty)$, given in terms of the solution $\bs{M}$ to \eqref{RM MDE} by
\[
\rho(\1\zeta\1)\,=\, \frac{1}{\pi N}\im \tr \bs{M}(\1\zeta\1)\,,
\] 
is the harmonic extension of a H\"older-continuous probability density $\rho: \R \to [0,\infty)$ (cf. Proposition~\ref{prp:Density of states}), which is called the \emph{self-consistent density of states} (cf. Definition~\ref{def:Density of states}).

\begin{theorem}[Local law for correlated random matrices] 
\label{thr:Local law for correlated random matrices}
Let $\bs{G}$ be the resolvent of a random matrix $\bs{H}$ written in the form \eqref{H structure} that satisfies {\bf B1}-{\bf B4}. 
For all $\delta,\eps>0$ and $\nu \in \N$ there exists a positive constant $C$ such that in the bulk,
\bels{entrywise local law in bulk}{
\P\Biggl[\,
\exists\2 \zeta \in \Cp\,\text{ s.t. }\rho(\zeta)\2\ge\2 \delta\,,\;  \im \zeta \ge N^{-1+\eps},\, \max_{x\1,\2y\1=\11}^N\abs{\1G_{xy}(\1\zeta\1)-m_{xy}(\1\zeta\1)}\2\ge\2 \frac{N^\eps}{\sqrt{N \im \zeta}\2}
\Biggr]
\leq
\frac{C}{N^{\nu}\!}\,.
}
Furthermore, the normalized trace converges with the improved rate 
\bels{trace local law in bulk}{
\P\Biggl[\,
\exists\2 \zeta \in \Cp\,\text{ s.t. }\rho(\zeta)\2\ge\2 \delta\,,\;  \im \zeta \ge N^{-1+\eps},\,  \absbb{\frac{1}{N}\!\tr\bs{G}(\1\zeta\1)-\frac{1}{N}\!\tr\bs{M}(\1\zeta\1)}\ge \frac{N^\eps}{{N \im \zeta}\2}
\Biggr]
\leq
\frac{C}{N^{\nu}\!}
\,.
}
The constant $C$ depends only on the model parameters $\scr{K}$ in addition to $\delta$, $\eps$ and $\nu$. 
\end{theorem}

In Section~\ref{sec:Local law for correlated random matrices} we present the proof of Theorem~\ref{thr:Local law for correlated random matrices} that is based on the results from Section~\ref{subsec:The Matrix Dyson equation} about the Matrix Dyson Equation.
As a  standard consequence of the local law \eqref{entrywise local law in bulk} and the uniform boundedness of $\im m_{xx}$ 
from Theorem~\ref{thr:Arbitrarily high polynomial decay of solution}, the eigenvectors of $\bs{H}$ in the bulk  are completely delocalized. 
This directly follows from the uniform boundedness of $\im G_{xx}(\1\zeta\1)$ and spectral decomposition of the resolvent (see e.g. \cite{EKYY}). 
 
\begin{corollary}[Delocalization of eigenvectors]
Pick any  $\delta,\eps,\nu>0$  and let $\bs{u}$ be  a normalized, $\norm{\bs{u}}=1$, eigenvector of $\bs{H}$, corresponding to an eigenvalue $\lambda \in \R$ in the bulk, i.e., $\rho(\lambda)\ge \delta$. Then 
\[
\P\biggl[\2
\max_{x =1}^N|\1u_x|\2\ge\2 \frac{N^\eps}{\!\sqrt{N\2}}%
\biggr]
\,\le\, \frac{C}{N^\nu\!}
\;,
\] 
for a positive constant $C$, depending only on the model parameters $\scr{K}$ in addition to $\delta$, $\eps$ and $\nu$. 
\end{corollary}
 The averaged local law \eqref{trace local law in bulk} directly implies the rigidity of the eigenvalues in the bulk. 
For any $\tau\in \R$, we define 
\vspace{-0.3cm}
\bels{energy index}{
i(\tau) \,:=\, \bigg\lceil  N \!\int_{-\infty}^\tau\msp{-10}\rho(\omega)\2\dd \omega\bigg\rceil
\,.
}

This is the index of an eigenvalue that is typically close to a spectral parameter $\tau$ in the bulk.  Then the standard argument presented in Section~\ref{subsec:Rigidity} proves the following result. 
\begin{corollary}[Rigidity]
\label{crl:Rigidity}
For any $\delta,\eps,\nu>0$ we have
\bels{bulk rigidity}{
\P\2\sbb{\2
\sup\big\{\,\abs{\1\lambda_{i(\tau)} - \tau\2}\; : \, \tau \in \R\,, \;\rho(\tau)\ge \delta\2\big\}
\ge 
\frac{N^\eps\!}{N}\,
\2}
\,\le \frac{C}{\2N^\nu\!}
\;,
}
for a positive constant $C$, depending only on the model parameters $\scr{K}$ in addition to $\delta$, $\eps$ and $\nu$.
\end{corollary}

Another consequence of Theorem~\ref{thr:Local law for correlated random matrices} is the universality of the local eigenvalue statistics in the bulk of the spectrum of $\bs{H}$
 both in the sense of averaged correlation functions and in the sense of gap universality. 
For the universality statement we make the following additional assumption that is stronger than {\bf B4}:
\begin{itemize}
\item[{\bf B5}] \emph{Fullness:} We say that $\bs{H}$ is $\beta=1$ ($\beta =2$) - \emph{full} if $\bs{H}\in \R^{N \times N}$ is real symmetric ($\bs{H}\in \C^{N \times N}$ is complex hermitian) and there is a positive constant $\kappa_5$ such that
\[
\E\2\abs{\1\tr \bs{R}\bs{W}\1}^2\,\ge\, \kappa_5 \tr \bs{R}^2
\,,
\]
for any real symmetric $\bs{R}\in \R^{N \times N}$ (any complex hermitian $\bs{R}\in \C^{N \times N}$).
\end{itemize}
When {\bf B5} is assumed 
we consider $\kappa_5$ as an additional model parameter. 

The first formulation of the bulk universality states that the $ k $-point correlation functions $\rho_k$ of the eigenvalues   of $ \bs{H}  $, rescaled around an energy parameter $\omega $ in the bulk, converge weakly to those of the GUE/GOE. The latter are given by  the correlation functions of well known determinantal processes. The precise statement is the following:

\begin{corollary}[Correlation function bulk universality]
\label{crl:Bulk universality} 
Let $\bs{H}$ satisfy {\bf B1}-{\bf B3} and {\bf B5} with $\beta=1$ ($\beta=2$). 
Pick any  $ \delta > 0 $  and choose any $\omega \in \R $ with $\rho(\omega) \geq \delta $.
Fix $ k \in \N $ and $ \eps \in (0,1/2)  $.
Then for any smooth, compactly supported test function $ \Phi : \R^k \to \R $ the $ k $-point local correlation functions  $ \rho_k : \R^k \to  [\10,\infty) $  of the eigenvalues of $ \bs{H}$ 
converge to the $ k $-point correlation function $ \Upsilon_k: \R^k \to [\10,\infty)$ of the GOE(GUE)-determinantal point process, 
\[
\Biggl|\;
\int_{\R^k} \Phi(\bs{\tau})\,
\Biggr[\,
\frac{1}{\1\rho(\omega)^k\!}\msp{5}
\rho_k\msp{-2}\bigg(\omega + \frac{\tau_1}{N\rho(\omega)}\2,\dots,\omega + \frac{\tau_k}{N \rho(\omega)}\bigg)-\Upsilon_{\!k}(\bs{\tau})\,\Biggr]\,
\dd \bs{\tau}
\,\Biggr|
\,\leq\, 
\frac{C}{N^c\!}
\,,
\]
where $ \bs{\tau} = (\tau_1,\dots,\tau_k) $, and the positive constants $ C,c $ depend only on $\delta$, $\Phi$ and the model parameters. 
\end{corollary}

The second formulation compares the joint distributions of gaps between consecutive eigenvalues of $ \bs{H}$ in the bulk
with those of the GUE/GOE. The proofs of Corollaries~\ref{crl:Bulk universality} and \ref{crl:Bulk Gap Universality} are presented in Section~\ref{subsec:Bulk universality}.

\begin{corollary}[Gap universality in bulk]
\label{crl:Bulk Gap Universality}
Let $\bs{H}$ satisfy {\bf B1}-{\bf B3} and {\bf B5} with $\beta=1$ ($\beta=2$). 
Pick any $\delta>0$, an energy $\tau$ in the bulk, i.e.  $\rho(\tau) \ge \delta $, and let $i=i(\tau)$ be 
the corresponding index defined in \eqref{energy index}.
Then for all $ n \in \N $ and all smooth compactly supported observables $ \Phi:\R^n \to \R$, there are two positive constants $ C $ and $ c $, 
depending on $n$,  $ \delta $, $ \Phi $ and the model parameters, such that the local eigenvalue distribution is universal, 
\bea{
&\absbb{\,\E\2\Phi\Bigr( \bigl(\2N\msp{-1}\rho(\1\lambda_i)\2(\1\lambda_{i+j}-\lambda_i)\2\bigr)_{j=1}^n \Bigr)
-
\E_{\rm G}\2\Phi\Bigl( \bigl(\2N\msp{-1}\rho_{\rm sc}(0)\2(\1\lambda_{\ceil{N/2}+j}-\lambda_{\ceil{N/2}})\2\bigr)_{j=1}^n \Bigr)
}
\,\leq\,
\frac{C}{N^c\!}\,.
}
Here the second expectation $\E_{\rm G}$ is 
with respect to GUE and GOE in the cases of complex Hermitian and real symmetric $\bs{H}$, respectively, and $\rho_{\rm sc}(0)=1/(2\1\pi)$ is the value of Wigner's semicircle law at the origin. 
\end{corollary}

During the final preparation of this manuscript and after announcing our theorems, we learned that a similar universality result but  with a special correlation structure was proved independently in  \cite{Che2016}. The  covariances   in \cite{Che2016}  have  a specific 
finite range and translation  invariant structure,
\begin{equation}
\label{checov}
\E\, h_{xy} h_{uv} = \psi \Big( \frac{x}{N}, \frac{y}{N}, u-x, v-y\Big),
\end{equation}
where $\psi$ is a piecewise Lipschitz function with finite support in the third and fourth variables.
The short scale translation invariance in \eqref{checov} allows one to use partial Fourier transform after effectively decoupling the slow variables from the fast ones. This 
renders the matrix equation \eqref{MDE in intro} into a vector equation  for $N^2$ variables and  the necessary stability result directly
 follows from \cite{AEK1}. The main difference between the current work and \cite{Che2016}
is that here we analyze  \eqref{MDE in intro} as a genuine matrix equation without relying on translation invariance and thus arbitrary short range correlations are allowed. 

\section{Local law for random matrices with correlations}
\label{sec:Local law for correlated random matrices}

In this section we show how the stability of the MDE, Theorem~\ref{thr:Stability},
can be combined with probabilistic estimates for random matrices with correlated entries to obtain a conceptually simple proof of the local law, Theorem~\ref{thr:Local law for correlated random matrices}.
We state these probabilistic estimates in Lemma~\ref{lmm:Smallness of error matrix} and Proposition~\ref{prp:fluctuation averaging} below before applying them to establish the local law.
Their proofs  are postponed to Sections~\ref{sec:Estimating the error term} and \ref{sec:Fluctuation averaging}, respectively. 

Consider any self-adjoint random matrix $ \bs{H} $, and let $(\bs{A},\cal{S})$ be the data pair for the MDE generated by the first two moments of $ \bs{H} $ through \eqref{RM MDE}.
Clearly, the self-energy operator $ \cal{S}: \C^{N \times N} \to \C^{N \times N} $ generated by \eqref{RM MDE} is self-adjoint with respect to the scalar product \eqref{matrix scalar product}, and preserves the cone of positive semidefinite matrices. 
The next lemma, whose proof is postponed to end of this section, shows that also the other assumptions with regards to our MDE results in Section~\ref{subsec:The Matrix Dyson equation} are satisfied for random matrices considered in Section~\ref{subsec:Correlated random matrices}.

\begin{lemma}[MDE data generated by random matrices] 
\label{lmm:MDE data from RM} 
If $\bs{H}$ satisfies {\bf B1}-{\bf B4}, then the data pair $(\bs{A},\cal{S})$ generated through \eqref{RM MDE} satisfies {\bf A1} and {\bf A2}. The corresponding model parameters $\scr{P}_1$ and $\scr{P}_2$ depend only on $\scr{K}$. 
\end{lemma}

In order to apply the stability of the MDE, we first write the defining equation for the resolvent \eqref{definition resolvent}, namely $-\bs{1}=(\zeta\1\bs{1}-\bs{H}\1)\1\bs{G}(\zeta)$,  of $ \bs{H} $ into the form
\begin{subequations}
\label{D-perturbed MDE and D}
\bels{G MDE with error D}{
-\bs{1}\,=\, (\1\zeta\2\bs{1}-\bs{A}+\cal{S}[\bs{G}(\1\zeta\1)]\1)\bs{G}(\1\zeta\1) + \bs{D}(\1\zeta\1)
\,,
}
a perturbed version of the MDE \eqref{MDE}.
Here the \emph{error matrix} $\bs{D} : \Cp \to \C^{N\times N} $  is given by
\bels{definition of error matrix D}{
\bs{D}(\1\zeta\1)\,:=\, -\1(\2\cal{S}[\bs{G}(\1\zeta\1)]+\bs{H}-\bs{A})\1\bs{G}(\1\zeta\1)
\,.
}
\end{subequations}
We view the resolvent $ \bs{G}(\zeta) $ as a perturbation of the deterministic matrix $ \bs{M}(\zeta) $  induced by the random perturbation $ \bs{D}(\zeta)$. 
Using the notation $\bs{\frak{G}}= \bs{\frak{G}}_\zeta $ from Theorem~\ref{thr:Stability}, we identify from \eqref{MDE}
and \eqref{G MDE with error D},  $ \bs{M}(\zeta) = \bs{\frak{G}}_\zeta(\bs{0}) $ and $ \bs{G}(\zeta) = \bs{\frak{G}}_\zeta(\bs{D}(\zeta)) $. Thus Theorem~\ref{thr:Stability} yields the following:

\NCorollary{Stability for local laws}{
Assume $\bs{H} $ satisfies {\bf B1}-{\bf B4}, fix $\delta>0$ and $ \zeta \in \D_\delta$. There exist constants $ c,C > 0 $, depending only on the model parameters and $\delta$, such that on the event where the a-priori bound
\bels{a-priori bound for stability}{
\norm{\bs{G}(\zeta)-\bs{M}(\zeta)}_{\rm{max}} \leq c
\,,
} 
holds, the difference $ \bs{G}(\zeta)-\bs{M}(\zeta) $ is bounded in terms of the perturbation $ \bs{D}(\zeta) $ by the two estimates:
\begin{align}
\label{LL max-stability}
\norm{\bs{G}(\zeta)-\bs{M}(\zeta)}_{\rm{max}}
&\leq\, C\2\norm{\bs{D}(\zeta)}_{\rm{max}}
\\
\label{LL avg-stability}
{\textstyle \frac{1}{N}}\abs{\1\tr(\bs{G}(\zeta)-\bs{M}(\zeta))}
\,&\leq\, 
\abs{\avg{\1\bs{J}(\zeta),\bs{D}(\zeta)}}+C\2\norm{\bs{D}(\zeta)}_{\rm{max}}^2
\,,
\end{align}
for some non-random $ \bs{J}(\zeta)\in \C^{N\times N} $ with fast decay, $\norm{\bs{J}(\zeta)}_{\ul{\gamma}} \leq C $, where the sequence $ \ul{\gamma} $ is from Theorem~\ref{thr:Stability}.
} 
\begin{Proof}
By Lemma~\ref{lmm:MDE data from RM} assumptions {\bf A1} and {\bf A2}  are satisfied for the data pair $(\bs{A}, \cal{S})$.
Hence, the first bound \eqref{LL max-stability} follows directly from \eqref{MDE stability2} with $\bs{D}_1=\bs{0}$ and $\bs{D}_2=\bs{D}(\zeta)$.
For the second bound \eqref{LL avg-stability} we first write $\frac{1}{N}\tr[\1\bs{G}-\bs{M}]= \avg{\bs{1}, \bs{\frak{G}}(\bs{D})-\bs{\frak{G}}(\bs{0})} $, then use the analyticity of $ \bs{\frak{G}} $ and the representation \eqref{derivative representation with Z} of its derivative  to obtain
\[
\avg{\2\bs{1}, \bs{\frak{G}}(\bs{D})-\bs{\frak{G}}(\bs{0})\1}
\,=\, \scalar{\2\bs{1}}{\cal{Z}[\bs{D}]+\bs{M}\bs{D}\1}+ \ord\bigl(\2\norm{\bs{D}}_{\rm{max}}^2\bigr)
\,=\, \scalar{\cal{Z}^*[\bs{1}]+\bs{M}^*\!}{\bs{D}\1}+ \ord\bigl(\2\norm{\bs{D}}_{\rm{max}}^2\bigr)
\,.
\]
Identifying $ \bs{J}:= \cal{Z}^\ast[\bs{1}] + \bs{M}^\ast $ yields \eqref{LL avg-stability}. The fast off-diagonal decay of the entries of $\bs{J}$ follows from \eqref{Arbitrarily high polynomial decay of solution} and \eqref{bounds on Z}. 
\end{Proof} 
 
Corollary~\ref{crl:Stability for local laws} shows, that on the event where the rough a-priori bound  \eqref{a-priori bound for stability} holds, the proof of the local law \eqref{entrywise local law in bulk} and \eqref{trace local law in bulk} is reduced to bounding the error $ \bs{D} $ on the right hand sides of \eqref{LL max-stability} and \eqref{LL avg-stability} by $ (N\im \zeta)^{-1/2} $ and $ (N\im \zeta)^{-1} $, respectively.   
In order, to state such estimates for the error matrix  we use the notion of \emph{stochastic domination}, first introduced in \cite{EKYY}, that is designed to compare random variables up to $N^\eps$-factors on very high probability sets. 

\begin{definition}[Stochastic domination] 
\label{def:Stochastic domination}
Let $X=X^{(N)}$, $Y=Y^{(N)}$ be sequences of non-negative random variables. We say $X$ is \emph{stochastically dominated by} $Y$ if 
\[
\P\bigl[X \2>\2 N^\eps Y\bigr]\,\le\, C(\eps, \nu)N^{-\nu}, \qquad N \in \N\,,
\]
for any $\eps>0$, $\nu \in \N$ and some ($N$-independent) family of positive constants $C$. In this case we write $X\prec Y$. 
\end{definition}

In this paper the family $C$ of constants in Definition~\ref{def:Stochastic domination} will always be an explicit function of the model parameters \eqref{definition of model parameters K} and possibly some additional parameters that are considered fixed and apparent from the context. However, the constants are always uniform in the spectral parameter $\zeta$ on the domain under consideration and indices $x,y$ in case $X=r_{xy}$ is the element of a matrix $\bs{R}=(r_{xy})_{x,y}$.  To use the notion of stochastic domination, we will think of $\bs{H}=\bs{H}^{(N)}$ as embedded into a sequence of random matrices with the same model parameters. 

The following lemma asserts that the error matrix $\bs{D}$ from \eqref{definition of error matrix D} converges to zero as the size $N$ of the random matrix grows to infinity.

\begin{lemma}[Smallness of perturbation in max-norm]
\label{lmm:Smallness of error matrix} Let $C>0$  and $\delta, \eps>0$ be fixed. 
Away from the real axis the error matrix $\bs{D}$ is small without regardless of an  a-priori bound on $ \bs{G}-\bs{M} $: 
\bels{bound on error D away from real line}{
\norm{\bs{D}(\tau+\ii \eta)}_{\rm{max}}\,\prec\,\frac{1}{\sqrt{N}}\,,\qquad \tau \in  [-C,C]\,,\; \eta \in [1, C]\,.
} 
Near the real axis and  in the regime where the harmonic extension of the self-consistent density of states is bounded away from zero, we have 
\bels{bound on error D in the bulk}{
\norm{\1\bs{D}(\1\zeta\1)}_{\rm{max}}\,\bbm{1}\bigl(\2\norm{\bs{G}(\zeta)-\bs{M}(\zeta)}_{\rm{max}}
\,\leq\, 
N^{-\eps}\bigr)\,\prec\,\frac{1}{\sqrt{N \im \zeta}}\,,
}
for all $\zeta \in \Cp$ with $\rho(\zeta)\ge \delta$ and $\im \zeta \ge N^{-1+\eps}$.
\end{lemma}

The proof of this key technical result is postponed to Section~\ref{sec:Estimating the error term}.
In order to bound the first term on the right hand side of \eqref{LL avg-stability} we use the following \emph{fluctuation averaging mechanism} (introduced in \cite{EYYber} for Wigner matrices) to improve the bound \eqref{bound on error D in the bulk} to a better bound for the inner product $ \avg{\bs{J},\bs{D}} $, given a version of the entry-wise local law.

\begin{proposition}[Fluctuation averaging]
\label{prp:fluctuation averaging}
Assume {\bf B1-B4}, and let $ [\kappa_-,\kappa_+] $ be the convex hull of $ \supp \rho $. 
Let $\delta,C>0$ and $\zeta\in \Cp$ with $\delta \le \dist(\zeta,[\kappa_-,\kappa_+])+\rho(\1\zeta\1)\le \delta^{-1}$  and  $\dist(\zeta, \spec(\bs{H}^B))^{-1}\prec  N^{C}$  for all $B \subsetneq\bb{X} $.
Let $\eps\in(0,1/2)$ be a constant and $\Psi$ a non-random control parameter with  $N^{-1/2}\le\Psi\le N^{-\eps}$. Suppose that the entrywise local law holds in the form
\bels{Assumption of FA on G-M}{
\norm{\bs{G}(\1\zeta\1)-\bs{M}(\1\zeta\1)}_{\rm{max}}\,\prec\, \Psi\,.
}
Then the error matrix $\bs{D}$, defined in \eqref{definition of error matrix D}, satisfies
\bels{general FA}{
\abs{\avg{\1\bs{R}\1,\bs{D}(\1\zeta\1)}}\,\prec\, \Psi^2\,,
}
for every non-random $\bs{R} \in \C^{N \times N}$ with faster than power law decay.
\end{proposition}
Note that Proposition~\ref{prp:fluctuation averaging} is stated on a slightly larger domain of spectral parameters than Theorem~\ref{thr:Local law for correlated random matrices} as it allows $\zeta$ to be away from the convex hull of  $\supp\rho$ even if $\rho(\zeta)$ is not bounded away from zero. This slight extension will be needed in Section~\ref{sec:Bulk universality and rigidity}.
The proof of Proposition~\ref{prp:fluctuation averaging} is carried out in Section~\ref{sec:Fluctuation averaging}. 
We have now stated all the results needed to prove the local law. 
In order to keep formulas short,  will use the notation:
\bels{definition of Lambda}{
\Lambda(\1\zeta\1)\,:=\, \norm{\bs{G}(\1\zeta\1)-\bs{M}(\1\zeta\1)}_{\rm{max}}
\,.
}
\begin{Proof}[Proof of Theorem~\ref{thr:Local law for correlated random matrices}]
We will start with the proof of \eqref{entrywise local law in bulk}.
By the Stieltjes transform representation \eqref{M as Stieltjes transform} of $\bs{M}$ and the trivial bound $\norm{\bs{G}(\zeta)} \le \frac{1}{\im \zeta}$ the norm of the difference $ \Lambda(\zeta) $ converges to zero as $ \zeta $ moves further away from the real axis.
In particular, the a-priori bound \eqref{a-priori bound for stability} needed for Corollary~\ref{crl:Stability for local laws} automatically holds for sufficiently large $ \im \zeta $.
Thus combining the corollary with the unconditional error bound \eqref{bound on error D away from real line} the estimate \eqref{LL max-stability}  takes the form 
\bels{rough bound for large im z}{
\Lambda(\tau + \ii\1\eta_\ast)\,\prec\, \frac{1}{\sqrt{N}}\,,\qquad \tau \in [-C_1,C_1]\,,
}
for any fixed constant $ C_1>0 $ and sufficiently large $\eta_\ast$.

Now let $\tau \in \R$, $\eta_0 \in [N^{-1+\eps},\eta_\ast]$ and $\zeta_0=\tau +\ii \1\eta_0 \in \Cp$ such that $\rho(\zeta_0)\ge \delta$ for some $\delta \in (0,1]$. Note that $\rho(\zeta_0)\ge \delta$ and $\eta_0 \le \eta_\ast$ imply $\tau \in [-C_1,C_1]$ for some positive constant $C_1$ because $\rho$ is the harmonic extension of the  self-consistent density of states with compact support  in $[-\kappa,\kappa]$ (Proposition~\ref{prp:Stieltjes transform representation}).
Since in addition the  self-consistent density of states is uniformly H\"older continuous (cf. Proposition~\ref{prp:Density of states}), there is a constant $c_1$, depending on $\delta$ and $\scr{P}$, such that
$
\inf_{\eta \in [\eta_0,\eta_\ast]}\rho(\tau+\ii\1\eta)\ge c_1
$. 
Therefore, by \eqref{bound on error D in the bulk} and \eqref{MDE stability2} we infer that
\bels{Gap in the values of G - M max-norm}{
\Lambda(\1\zeta\1)\,\bbm{1}(\Lambda(\1\zeta\1)\le N^{-\eps/4})\,\prec\, \frac{1}{\sqrt{N\im \zeta}}
\,,
\qquad 
\zeta \in \tau + \ii\1[\eta_0,\eta_\ast]\,.
}
Since $N^{-\eps/2} \ge (N\im \zeta)^{-1/2}$,  
the inequality \eqref{Gap in the values of G - M max-norm} establishes on a high probability event a gap in the set possible values that $\Lambda(\zeta)$ can take. 
The indicator function in \eqref{Gap in the values of G - M max-norm} is absent for $\zeta=\tau +\ii\1\eta_\ast$ because of \eqref{rough bound for large im z}, i.e. at that point the value lies below the gap. From the Lipshitz-continuity of $\zeta \mapsto \Lambda(\zeta)$ with Lipshitz-constant bounded by $2N^2$ for $\im \zeta \ge \frac{1}{N}$ and a standard continuity argument together with a union bound (e.g. Lemma~A.1 in \cite{AEK2}), we conclude that $ \Lambda(\zeta) $ lies below the gap for any $\zeta$ with $\im\zeta \in [\eta_0, \eta_\ast]$ with very high probability. 
 Thus, using the definition of stochastic domination, we see that $ \max_{\zeta \in \D_\delta} \Lambda(\zeta) \prec (N\im\zeta)^{-1/2} $, i.e., the entrywise  local law  \eqref{entrywise local law in bulk} holds.

Now we prove \eqref{trace local law in bulk}. 
Let $\zeta \in \Cp$ with $\im \zeta \ge N^{-1+\eps}$ and $\rho(\zeta)\ge \delta $, so that the entrywise local law \eqref{Assumption of FA on G-M} holds at $ \zeta $ with $ \Psi := (N\im\zeta)^{-1/2} $.
Applying the fluctuation averaging (Proposition~\ref{prp:fluctuation averaging}) with $ \Psi := (N\im\zeta)^{-1/2} $ and $ \bs{R} := \bs{J}(\zeta) $ yields $ \abs{\avg{\bs{J},\bs{D}}} \prec (N\im\zeta)^{-1} $. 
Plugging this estimate for the first term into the right hand side of \eqref{LL avg-stability}, and recalling the definition of stochastic domination,  yields \eqref{trace local law in bulk}. This finishes the proof of Theorem~\ref{thr:Local law for correlated random matrices}.  
\end{Proof}

\begin{Proof}[Proof of Lemma~\ref{lmm:MDE data from RM}]
The condition \eqref{Expectation decay} on the bare matrix $\bs{A}$ is clearly satisfied by \eqref{decay of expectation}. The lower bound on $\cal{S}$ in \eqref{Flatness} follows from \eqref{RM Flatness}. To show this, let $\bs{R}=\sum_{\1i} \varrho_i \2\bs{r}_i\bs{r}_i^* \in \ol{\scr{C}}_+$,  where the sum is over the orthonormal basis $ (\bs{r}_i)_{i=1}^N $. 
Then $ \bs{v}^*\cal{S}[\bs{R}]\1\bs{v} =  \sum_{\1i}\varrho_i\2 \E\2 \abs{\1\bs{r}_i^*\bs{W}\bs{v}}^2 \ge \kappa_4\sum_{\1j} \varrho_i $, 
for any normalized vector $\bs{v}\in \C^N$. 

We will now verify the upper bounds on $\cal{S}$ in \eqref{Flatness} and \eqref{S operator decay}. Both bounds follow from the decay of covariances 
\bels{decay of covariances}{
\abs{\1\E\2 w_{xu}w_{vy}}\,\le\,
\kappa_3(2\nu)\,\pb{\2
(\1q_{xy}\1q_{uv})^{\nu} + (\1q_{xv}\1q_{uy})^{\nu} 
}
\,,\qquad q_{\1xy}\,:=\, 
{\textstyle \frac{1}{\11\2+\2d(x,\2y)\1}} 
\,, 
\qquad \nu \in \N
\,,
}
which is an immediate consequence of \eqref{decay of correlation} with the choices $W_A=(w_{xu},w_{ux})$, $W_B=(w_{vy},w_{yv})$, $\phi(\xi_1,\xi_2)=\ol{\xi}_1$ and $\psi(\xi_1,\xi_2)={\xi}_1$. 

Indeed, to see the upper bound in  \eqref{Flatness} it suffices to show 
\bels{bound on S on rank one}{
\norm{\1\cal{S}[\1\bs{r}\1\bs{r}^*]\1}\,\le\, C\2N^{-1}\,,
}
for a constant $C>0$, depending on $\scr{K}$, and 
 any normalized vector $\bs{r} \in \C^N$, because we can use for any $\bs{R} \in \ol{\scr{C}}_+$  the spectral decomposition $\bs{R}=\sum_i \varrho_i \2\bs{r}_i\bs{r}_i^* $ as above. The estimate \eqref{decay of covariances} yields
\bels{Bound S by Qnu}{
\norm{\1\cal{S}[\1\bs{r}\1\bs{r}^*]\1}
\,\le\,
\frac{\kappa_3(2 \nu)}{N}
\pB{\,
\norm{\1 \bs{Q}^{(\nu)}}\msp{5} \abs{\bs{r}}^*\bs{Q}^{(\nu)}\msp{-1}\abs{\bs{r}} \,+\,\normb{(\bs{Q}^{(\nu)}\abs{\bs{r}})(\bs{Q}^{(\nu)}\abs{\bs{r}})^*}\,}
\,,
}
where we defined the matrix $\bs{Q}^{(\nu)}$ with entries $q_{xy}^{(\nu)}:=q_{xy}^{\2\nu}$ and $\abs{\bs{r}}:=(\abs{r_x})_{x\in \bb{X}}$. Since
\bels{Estimate on Qnu}{
\norm{\bs{Q}^{(\nu)}}\,\le\, 
\max_x \sum_y\2 q_{xy}^{\2\nu}
\,,
}
the inequality \eqref{bound on S on rank one} follows from the sub-$ P$-dimensional volume \eqref{polynomial ball growth} by choosing $\nu$ sufficiently large. 

To show  \eqref{S operator decay} we fix any $\bs{R} \in \C^{N \times N}$ and estimate the entries $\cal{S}_{xy}[\bs{R}]$ of $\cal{S}[\bs{R}]$ by  using \eqref{decay of covariances},
\[
\abs{\cal{S}_{xy}[\bs{R}]}\,\le\, \kappa_3(2 \nu)\2
\pbb{
\pB{\frac{1}{N}\sum_{u,v}q_{uv}^{\2\nu}}\2q_{xy}^{\2\nu}+\frac{1}{N}\pB{\sum_vq_{xv}^{\2\nu}}\pB{ \sum_u q_{uy}^{\2\nu}}
}
\,\norm{\bs{R}}_{\rm{max}}
\]
The bound \eqref{S operator decay} follows because the right hand side of \eqref{Estimate on Qnu} is finite. 
\end{Proof}

\section{The Matrix Dyson Equation}
\label{sec:The Matrix Dyson Equation}
This section is dedicated to the analysis of the MDE \eqref{MDE}. In particular, it is thus independent of the probabilistic results established in Sections~\ref{sec:Estimating the error term} - \ref{sec:Bulk universality and rigidity}. In Section~\ref{subsec:The solution of the Matrix Dyson Equation} we establish a variety of properties of the solution $\bs{M}$ to the MDE. The section starts with the proof of Proposition~\ref{prp:Stieltjes transform representation} and ends with the proof of Theorem~\ref{thr:Arbitrarily high polynomial decay of solution}. In Section~\ref{subsec:Stability of the Matrix Dyson Equation} we prove Proposition~\ref{prp:Density of states} and the stability of the MDE,  Theorem~\ref{thr:Stability}.

\subsection{The solution of the Matrix Dyson Equation}
\label{subsec:The solution of the Matrix Dyson Equation}

Most of the inequalities in this and the following section are uniform in the data pair $(\bs{A},\cal{S})$ that determines the MDE and its solution, given a fixed set of model parameters $\scr{P}_k$ corresponding to the assumptions {\bf Ak}. We therefore  introduce a convention for inequalities up to constants, depending only on the model parameters.

\begin{convention}[Comparison relation and constants]
\label{con:Comparison relation}
Suppose a set of model parameters $\scr{P}$ is given. Within the proofs we will write $C$ and $c$ for  generic positive constants, depending on $\scr{P}$.  In particular, $C$ and $c$ may change their values from inequality to inequality. If $C,c$ depend on additional parameters $\scr{L}$, we will indicate this by writing $C(\scr{L}),c(\scr{L})$.  We also use the comparison relation $\alpha \lesssim \beta$ or $\beta \gtrsim \alpha$  for any positive $\alpha$ and $\beta$ if there exists a constant $C>0$ that depends only on $\scr{P}$, but is otherwise uniform in the data pair $(\bs{A},\cal{S})$, such that $\alpha \le C \beta$. In particular, $C$ does not depend on the dimension $N$ or the spectral parameter $\zeta$. In case $\alpha\lesssim \beta\lesssim \alpha$ we write $\alpha \sim \beta$. For two matrices $\bs{R},\bs{T}\in \overline{\scr{C}}_+$ we similarly write $\bs{R}\lesssim\bs{T}$ if the inequality $\bs{R}\le C \bs{T}$ in the sense of quadratic forms holds with a constant $C>0$ depending only on the model parameters.  
\end{convention}

In the upcoming analysis many quantities depend on the spectral parameter $\zeta$. We will often suppress this dependence in our notation and write e.g. $\bs{M}=\bs{M}(\1\zeta\1)$, $\rho=\rho(\1\zeta\1)$, etc. 

\begin{Proof}[Proof of Proposition \ref{prp:Stieltjes transform representation}] 
In this proof we will  generalize 
the proof of Proposition~2.1 from \cite{AEK1cpam} to our  matrix  setup.
By taking the imaginary part of both sides of the MDE and using $\im \bs{M}\ge \bs{0}$ and $\bs{A}=\bs{A}^{\!*}$ we see that
\[
-\im \sb{\2\bs{M}(\1\zeta\1)^{-1}}\,=\, \bs{M}^*(\1\zeta\1)^{-1}\im \bs{M}(\1\zeta\1)\2\bs{M}(\1\zeta\1)^{-1} \,\ge\, \im \zeta\,\bs{1}\,.
\]
In particular, this implies the trivial bound on the solution to the MDE,
\bels{trivial bound on M}{
\norm{\bs{M}(\1\zeta\1)}\,\leq\, \frac{1}{\im\2\zeta}\,, \qquad \zeta \in \Cp\,.
} 

Let $\bs{w}\in \C^{N}$ be normalized, $\bs{w}^*\bs{w}=1$. Since $\bs{M}(\zeta)$ has positive imaginary part, the analytic function $\zeta \mapsto  \bs{w}^*\bs{M}(\zeta)\bs{w}$ takes values in $\Cp$. From the trivial upper bound \eqref{trivial bound on M} and the MDE itself, we infer the asymptotics 
$
\ii\1\eta\1  \bs{w}^*\bs{M}(\ii\1\eta)\bs{w}\to -1
$ as $\eta\to  \infty$. 
By the characterization of Stieltjes transforms of probability measures on the complex upper half plane (cf. Theorem 3.5 in \cite{Garnett-BA2007}), we infer
\[
 \bs{w}^*\bs{M}(\zeta)\bs{w}\,=\, \int\frac{{v}_{\bs{w}}(\dd \tau)}{\tau-\zeta}\,,
\]
where $v_{\bs{w}}$ is a probability measure on the real line.
By polarization, we find the general representation \eqref{M as Stieltjes transform}. 

We now show that $\supp\bs{V} \subseteq [-\kappa,\kappa]$, where $\kappa=\norm{\bs{A}}+2\norm{\cal{S}}^{1/2}$ (cf. \eqref{kappa definition}). 
Note that {\bf A1} implies $\norm{\cal{S}}\lesssim 1$. 
Indeed, letting $(\2\cdot\2)_\pm$ denote the positive and negative parts, we find 
\bels{bounded operator norm of S}{
\norm{\1\cal{S}[\bs{R}]}\,\le\,  P_1\2\bigl(\,\avg{\1(\re \bs{R})_+}+\avg{\1(\re \bs{R})_-}+\avg{\1(\im \bs{R})_+}+\avg{\1(\im \bs{R})_-}
\,\bigr)
\,\le\,
2\1P_1\norm{\bs{R}}_{\rm{hs}}\,,
}
for any $ \bs{R} \in \C^{N\times N}$.
Since $\norm{\bs{R}}_{\rm{hs}} \le \norm{\bs{R}}$ the bound $\norm{\cal{S}}\lesssim 1$ follows. 
The following argument will prove that $\norm{\im \bs{M}(\1\zeta\1)} \to 0$ as $\im \zeta \downarrow 0$ locally uniformly for all $\zeta \in \Cp$ with $ \abs{\zeta}>\kappa$. This implies $\supp\bs{V} \subseteq [-\kappa,\kappa]$.

Let us fix $\zeta \in \Cp$ with $ \abs{\zeta}>\kappa$ and suppose that $\norm{\bs{M}}$ satisfies the upper bound
\bels{assume bound on M}{
\norm{\bs{M}}\,<\, \frac{\abs{\zeta}-\norm{\bs{A}}}{2 \norm{\cal{S}}}\,.
}
Then by taking the inverse and then the norm on both sides of \eqref{MDE} we conclude that
\bels{get better bound on M}{
\norm{\bs{M}}\,\le\, \frac{1}{\abs{\zeta}-\norm{\bs{A}}-\norm{\cal{S}}\norm{\bs{M}}}\,\le\, \frac{2}{\2\abs{\zeta}-\norm{\bs{A}}}\,.
}
Therefore, \eqref{assume bound on M} implies \eqref{get better bound on M} and we see that there is a gap in the possible values of $\norm{\bs{M}}$, namely
\[
\qquad
\textstyle
\norm{\bs{M}(\1\zeta\1)} \not \in \pB{\frac{2}{\2\abs{\1\zeta}\1-\1\norm{\bs{A}}\2}\,,\2\frac{\2\abs{\zeta}\2-\1\norm{\bs{A}}}{2\1 \norm{\cal{S}}}}\qquad 
\text{for} \quad \abs{\zeta}\2>\2\kappa\,.
\]
Since $\zeta \mapsto\norm{\bs{M}(\1\zeta\1)} $ is a continuous function and for large $\im \zeta$ the values of this function lie below the gap by the trivial bound \eqref{trivial bound on M}, we infer
\bels{norm M below the gap}{
\norm{\bs{M}(\1\zeta\1)}\,\le\, \frac{2}{\abs{\zeta}-\norm{\bs{A}}}\qquad 
\text{for} \quad \abs{\zeta}\2>\2\kappa\,.
}
Let us now take the imaginary part of the MDE and multiply  it  with $\bs{M}^*$ from the left and with $\bs{M}$ from the right,
\bels{imaginary part of MDE}{
\im \bs{M}\,=\, (\1\im \zeta\1)\, \bs{M}^*\bs{M} +\bs{M}^*\cal{S}[\im \bs{M}]\bs{M}\,.
}
By taking the norm on both sides of \eqref{imaginary part of MDE}, using a trivial estimate on the right hand side and rearranging the resulting terms, we get 
\bels{norm im M estimate}{
\norm{\im \bs{M}}\,\le\,\frac{{ \im \zeta}\2\norm{\bs{M}}^2}{1-\norm{\bs{M}}^2\norm{\cal{S}}}\,. 
}
Here we used $\norm{\bs{M}}^2\norm{\cal{S}}<1$, which is satisfied by \eqref{norm M below the gap} for $ \abs{\zeta}>\kappa$. 
We may estimate the right hand side of \eqref{norm im M estimate} further by applying \eqref{norm M below the gap}. Thus we find
\bels{final bound on norm im M}{
\norm{\im \bs{M}}\,\le\,\frac{4\2{ \im \zeta}}{(\abs{\zeta}-\norm{\bs{A}})^2-4\1\norm{\cal{S}}}
\,=\, 
\frac{4\2{ \im \zeta}}{(\abs{\zeta}-\kappa+ 2 \norm{\cal{S}}^{1/2})^2-4\1\norm{\cal{S}}}
\,.
}
The right hand side of \eqref{final bound on norm im M} converges to zero locally uniformly for all $\zeta \in \Cp$ with $ \abs{\zeta}>\kappa$ as $\im \zeta \downarrow 0$.
This finishes the proof of Proposition~\ref{prp:Stieltjes transform representation}.
\end{Proof}

The following proposition lists bounds on $ \bs{M} $ that, besides the ones stated in Section~\ref{subsec:The Matrix Dyson equation}, constitute the only properties of $ \bs{M}$ that we need outside this section.

\begin{proposition}[Properties of the solution]
\label{prp:Properties of the solution}
Assume {\bf A1} and that $\norm{\bs{A}}\le P_0$ for some constant $P_0>0$.
 Then uniformly for all spectral parameters $\zeta \in \Cp$ the following bounds hold:
\begin{itemize}
\item[\emph{(i)}] The solution is bounded in the spectral norm,
\bels{M upper bound}{
\norm{\1\bs{M}(\1\zeta\1)}\,\lesssim\, \frac{1}{\rho(\1\zeta\1)+\dist (\zeta, \supp \rho)}\,.
}
\item[\emph{(ii)}] The inverse of the solution is bounded in the spectral norm,
\bels{M lower bound}{
\norm{\1\bs{M}(\1\zeta\1)^{-1}}\,\lesssim\,1+\abs{\zeta}\,.
}
\item[\emph{(iii)}] 
The imaginary part of $\bs{M}$ is comparable to the harmonic extension of the self-consistent density of states, 
\bels{im M lower and upper bound}{
\rho(\1\zeta\1)\2\bs{1}\,\lesssim\, \im \bs{M}(\1\zeta\1)\,\lesssim\, (1+\abs{\zeta}^2)\norm{\1\bs{M}(\1\zeta\1)}^2\rho(\1\zeta\1)\2\bs{1}\,.
}
\end{itemize}
\end{proposition}

\begin{Proof}
The inequalities \eqref{M upper bound} and \eqref{M lower bound} provide upper and lower bounds on the singular values of the solution, respectively.
Before proving these bounds we show that $\bs{M}$ has a bounded normalized Hilbert-Schmidt norm,
\bels{Bounded HS norm}{
\norm{\1\bs{M}(\1\zeta\1)}_{\rm{hs}}\,\lesssim\, 1\,,
\qquad \forall\1\zeta \in \Cp
\,.
}
For this purpose we take the imaginary part of \eqref{MDE} (cf. \eqref{imaginary part of MDE}) and find
$
\im \bs{M} 
\ge
\bs{M}^*\cal{S}[\im \bs{M}]\bs{M}
$,  where $ \bs{M} = \bs{M}(\zeta) $.  
The lower bound on $\cal{S}$ from \eqref{Flatness} implies
\bels{imaginary part of M lower bounded by rho}{
\im \bs{M} 
\,\gtrsim\, 
\rho\,\bs{M}^*\bs{M}\,,
}
where we used the definition of $\rho$ in \eqref{definition extended DOS}. Taking the normalized trace on both sides of \eqref{imaginary part of M lower bounded by rho} shows \eqref{Bounded HS norm}. 
\medskip
\\
\textit{Proof of }(ii):
Taking the norm on both sides of \eqref{MDE} yields 
\bels{Bound in proof of M2}{
\norm{\1\bs{M}^{-1}}\,\le\, \abs{\zeta}+\norm{\bs{A}}+\norm{\cal{S}}_{\rm{hs} \to \norm{\1\cdot\1}}\norm{\bs{M}}_{\rm{hs}}
\,\lesssim\, 1+\abs{\zeta}\,,
}
where $\norm{\cal{S}}_{\rm{hs} \to \norm{\1\cdot\1}}$ denotes the norm of $\cal{S}$ from $\C^{N\times N}$ equipped with the norm $\norm{\2\cdot\2}_{\rm{hs}}$ to $\C^{N\times N}$ equipped with $\norm{\2\cdot\2}$.
For the last inequality in \eqref{Bound in proof of M2} we used \eqref{Bounded HS norm} and that by {\bf A1}  we have $\norm{\cal{S}}_{\rm{hs} \to \norm{\1\cdot\1}}\lesssim 1$ (cf. \eqref{bounded operator norm of S}).
\medskip
\\
\textit{Proof of }(iii): First we treat the simple case of large spectral parameters, $\abs{\zeta}\ge 1+\kappa$, where $\kappa$ was defined in \eqref{kappa definition}. Recall that the matrix valued measure $\bs{V}(\dd \tau)$ (cf. \eqref{M as Stieltjes transform}) is supported in $[-\kappa,\kappa]$ by Proposition~\ref{prp:Stieltjes transform representation}. The normalization, $\bs{V}(\R)=\bs{1}$ implies that for any vector $\bs{u}\in \C^N$ with $\norm{\bs{u}}=1$ the function $\zeta \mapsto \frac{1}{\pi} \im [\bs{u}^*\bs{M}(\1\zeta\1)\bs{u}]$ is the harmonic extension of a probability measure with support in $[-\kappa,\kappa]$, hence it behaves as $-\zeta^{-1}$ for large $|\zeta|$. We conclude that 
$
\im \bs{M}(\1\zeta\1)\sim\rho(\1\zeta\1)\sim\abs{\zeta}^{-2}\im \zeta
$,
for $\abs{\zeta}\ge 1+\kappa$. Since for these $\zeta$ we also have $\norm{\bs{M}(\1\zeta\1)}\sim \abs{\zeta}^{-1}$ by the Stieltjes transform representation 
\eqref{M as Stieltjes transform} we conclude that \eqref{im M lower and upper bound} holds in this regime. 

Now we consider $\zeta \in \Cp$ with $\abs{\zeta}\le 1+\kappa$. We start with the lower bound on $\im \bs{M}$. 
From \eqref{imaginary part of M lower bounded by rho} we see that
$
\im \bs{M} 
\,\gtrsim\, 
\rho\,\norm{\bs{M}^{-1}}^{-2}\2\bs{1}
$,
and since $\norm{\1\bs{M}^{-1}}\lesssim 1$ by {(ii)}, the lower bound in \eqref{im M lower and upper bound} is proven.

 For the upper bound,   taking the imaginary part of the MDE (cf. \eqref{imaginary part of MDE}) and using {\bf A1} and that $\im \bs{M} \gtrsim \im \zeta \2\bs{1}$ by the Stieltjes transform representation \eqref{M as Stieltjes transform}, we get
\[
\im \bs{M} \,=\, \im \zeta\, \bs{M}^*\bs{M} + \bs{M}^*S[\1\im \bs{M}]\1 \bs{M}
\,\lesssim\, (\1\im \zeta +\rho\2)\2\bs{M}^*\bs{M}
\,\lesssim\, 
\rho\,\norm{\bs{M}}^2\2\bs{1}\,.
\]
\textit{Proof of }(i): In the regime $\abs{\zeta}\ge 1+\kappa$ the bound \eqref{M upper bound} follows from the Stieltjes transform representation \eqref{M as Stieltjes transform}. Thus we consider   $\abs{\zeta}\le 1+\kappa$. We take the imaginary part on both sides of \eqref{MDE} and use the lower bound in \eqref{im M lower and upper bound} and $\cal{S}[\bs{1}]\gtrsim \bs{1}$ to get
\[
-\im \bs{M}(\1\zeta\1)^{-1}\,\ge\, \cal{S}[\1\im \bs{M}(\1\zeta\1)]\,\gtrsim\,\rho(\1\zeta\1)\2\bs{1}\,.
\]
Since in general $\im \bs{R}^{-1}\ge \bs{1}$ implies $\norm{\bs{R}}\le 1$ for any $\bs{R} \in \C^{N \times N}$, 
we infer that $\norm{ \bs{M}(\1\zeta\1)}\lesssim \rho(\1\zeta\1)^{-1}$. On the other hand, 
$ \norm{ \bs{M}(\1\zeta\1)} \lesssim \dist(\1\zeta,\1 \supp \rho\1)^{-1} $ 
follows from \eqref{M as Stieltjes transform} again.
\end{Proof}
In order to show the fast decay of off-diagonal entries of $ \bs{M} $, Theorem~\ref{thr:Arbitrarily high polynomial decay of solution}, we rely on the following general result on matrices with decaying off-diagonal entries.

\begin{lemma}[Perturbed Combes-Thomas estimate]
\label{lmm:Perturbed Combes-Thomas} Let $\bs{R}\in \C^{N \times N}$ be such that  
\[
\abs{\1r_{xy}}\,\le\, \frac{\beta(\nu)}{(1+d(x,y))^\nu}+\frac{\beta(0)}{N}
\,,\qquad \forall\2x,y \in \bb{X}\,,\;\forall\1\nu \in \N
\,,
\]
with some positive sequence $\ul{\beta}=(\beta(\nu))_{\nu=0}^\infty$, and $\norm{\1\bs{R}^{-1}}\le1$.  

Then there exists a sequence  $\ul{\alpha}=(\alpha(\nu))_{\nu=0}^\infty$, depending only on $\ul{\beta}$ and $P$  (cf. \eqref{polynomial ball growth}), such that 
\bels{R inverse CT inequality}{
\abs{\1(\bs{R}^{-1})_{xy}}\,\le\, \frac{\alpha(\nu)}{(1+d(x,y))^\nu}+\frac{\alpha(0)}{N}
\,,\qquad \forall\2x,y \in \bb{X}\,,\;\forall\1\nu \in \N\,.
}
\end{lemma}

This lemma is reminiscent of a standard Combes-Thomas estimate: An off-diagonal decay of the entries
of a matrix $\bs{R}$ implies a similar decay for its inverse, $\bs{R}^{-1}$, provided the smallest singular value
is bounded away from zero. Indeed, in the case of  $\alpha(0)=\beta(0)=0$ the proof of this lemma 
directly follows from the standard strategy for establishing Combes-Thomas estimates, see e.g. Proposition~13.3.1. in \cite{PasturShcerbinaAMSbook}; we omit the details. We now explain how to extend this standard result to our case, where Lemma~\ref{lmm:Perturbed Combes-Thomas} allows for a
nondecaying component. The detailed proof will be given in the appendix, here we only present the basic idea.

Write $\bs{R}=\bs{S}+\bs{T}$, where $\bs{S}$ has a fast off-diagonal decay and $\bs{T}$ has entries of size $|t_{xy}|\lesssim N^{-1}$.
Note that $\bs{T}$ cannot simply be considered as a small
perturbation since its norm can be of order one, i.e. comparable with that of $\bs{S}$. 
Instead, the proof  relies on  showing that $\bs{S}$ inherits the lower bound on its singular values from $\bs{R}$
and  then applying the standard $\alpha(0)=\beta(0)=0$ version of the Combes-Thomas estimate to $\bs{S}$ to
generate the decaying component of $\bs{R}^{-1}$. The point is that $\bs{T}$ can potentially change only finitely many 
singular values by a significant amount since $\norm{\bs{T}}_{\rm{max}} \lesssim N^{-1}$.
If these few singular values were close to zero, then they would necessarily be isolated, hence
the corresponding singular vectors would be strongly localized. However, because of its small entries, $\bs{T}$ acts trivially on localized vectors which implies that isolated singular values are essentially stable under adding or subtracting $\bs{T}$. This argument excludes the creation of singular values close to zero by subtracting $\bs{T}$ from $\bs{R}$. The details are found in the appendix.
Putting all these ingredients together, we can now complete the proof of Theorem~\ref{thr:Arbitrarily high polynomial decay of solution}.

\begin{Proof}[Proof of Theorem~\ref{thr:Arbitrarily high polynomial decay of solution}]
Recall the model parameters $\ul{\pi}_1,\ul{\pi}_2$ from {\bf A2}. We consider  the MDE \eqref{MDE} entrywise and see that
\bea{
\abs{\1(\1\bs{M}^{-1})_{xy}} 
\;&\le\; 
\abs{\zeta}\2\delta_{xy}\,+
\frac{\pi_1(\nu)+\pi_2(\nu)\norm{\bs{M}}}{(1+d(x,y))^\nu} \,+\, \frac{\pi_1(0)+\pi_2(0)\norm{\bs{M}}}{N}
\,,
}
where we used  the assumptions \eqref{Expectation decay} and \eqref{S operator decay}, as well as  $\norm{\bs{M}}_{\rm{max}} \le \norm{\bs{M}}$.
By  \eqref{M upper bound} and $\zeta \in \D_\delta$, we have  $ \norm{\bs{M}} \lesssim \delta^{-1}$.  Furthermore, for large $|\zeta|$ we also  have $\norm{\bs{M}(\zeta)} \lesssim \abs{\zeta}^{-1}$.
We  can now  apply Lemma~\ref{lmm:Perturbed Combes-Thomas} with the choice $\bs{R}:=\norm{\bs{M}}\2\bs{M}^{-1}$ 
to see the existence of a positive sequence $\ul{\gamma}$ such that \eqref{Arbitrarily high polynomial decay of solution} holds. This finishes the proof of Theorem~\ref{thr:Arbitrarily high polynomial decay of solution}. 
\end{Proof}

\subsection{Stability of the Matrix Dyson Equation}
\label{subsec:Stability of the Matrix Dyson Equation}

The goal of this section is to prove Proposition~\ref{prp:Density of states} and  Theorem~\ref{thr:Stability}. The main technical result, which is needed for these proofs, is the linear stability of the MDE. 
For its statement we introduce for any $\bs{R} \in \C^{N \times N}$ the \emph{sandwiching operator} $\cal{C}_{\bs{R}}: \C^{N \times N}\to \C^{N \times N}$ by
\bels{sandwiching operator}{
\cal{C}_{\1\bs{R}}[\1\bs{T}\1]\,:=\, \bs{R}\bs{T}\bs{R}\,.
}
Note that $\cal{C}_{\bs{R}}^{-1}=\cal{C}_{\bs{R}^{-1}}$ and $\cal{C}_{\bs{R}}^*=\cal{C}_{\bs{R}^*}$ for any $\bs{R} \in \C^{N \times N}$, where $\cal{C}_{\bs{R}}^*$ denotes the adjoint with respect to the scalar product \eqref{matrix scalar product}.

\begin{proposition}[Linear stability]
\label{prp:Linear stability} Assume {\bf A1} and $\norm{\bs{A}}\le P_0$ for some constant $P_0>0$ (cf. \eqref{spectral norm bound on A}). 
There exists a universal numerical constant $C>0$ such that  uniformly for all $\zeta \in \Cp$: 
\bels{Linear stability}{
\norm{(\1\rm{Id}-\cal{C}_{\bs{M}(\zeta)}\cal{S}\1)^{-1}}_{\rm{sp}}
\,\lesssim\; 
1\,+\frac{1}{(\1\rho(\1\zeta\1)+\dist(\zeta,\supp\rho))^C}
\,.
}
\end{proposition}

Before we show a few technical results that prepare the proof of Proposition~\ref{prp:Linear stability}, we give a heuristic argument that explains how the operator $\rm{Id}-\cal{C}_{\bs{M}}\cal{S}$ on the left hand side of \eqref{Linear stability} is connected to the stability of the MDE \eqref{MDE}, written in the form 
$-\bs{1} = ( \zeta\1\bs{1}-\bs{A}+ \cal{S}[\bs{M}])\bs{M}$,
with respect to perturbations.  Suppose that the perturbed MDE
\bels{heuristic perturbed MDE}{ 
-\bs{1} \,=\, (\2\zeta\1\bs{1}-\bs{A}+ \cal{S}[\1\bs{\frak{G}}(\bs{D})]\2)\bs{\frak{G}}(\bs{D}) + \bs{D}\,,
}
with perturbation matrix $\bs{D}$ has a unique solution $\bs{\frak{G}}(\bs{D})$, depending differentiably on $\bs{D}$. Then by differentiating on both sides of \eqref{heuristic perturbed MDE} with respect to $\bs{D}$, setting $\bs{D}=\bs{0}$ and using the MDE for $\bs{M}(\zeta)=\bs{\frak{G}}(\bs{0})$, we find
\bels{Derivative of MDE at D eq 0}{
\bs{0}\,=\,-\bs{M}(\1\zeta\1)^{-1} \nabla_{\bs{R}}\bs{\frak{G}}(\bs{0})
+\cal{S}[\nabla_{\bs{R}}\bs{\frak{G}}(\bs{0})]\1\bs{M}(\1\zeta\1)+ \bs{R}\,,
}
where $ \nabla_{\bs{R}}$ denotes the directional derivative with respect to $\bs{D}$ in the direction $\bs{R} \in \C^{N \times N}$. Rearranging the terms in \eqref{Derivative of MDE at D eq 0} and multiplying with $\bs{M}=\bs{M}(\zeta)$ from the left yields
\bels{Derivative equation at D eq 0}{
(\1\rm{Id}-\cal{C}_{\bs{M}}\cal{S}\1)\nabla_{\bs{R}}\bs{\frak{G}}(\bs{0})\,=\,\bs{M}\bs{R}\,.
}
Thus $\bs{\frak{G}}(\bs{D})$ has a bounded derivative at $\bs{D}=\bs{0}$, i.e., the MDE is stable with respect to the perturbation $\bs{D}$ to linear order, whenever the operator $\rm{Id}-\cal{C}_{\bs{M}}\cal{S}$ is invertible and its inverse is bounded. 
In order to extend the linear stability to the full stability of the MDE for non-infinitesimal perturbations, the linear stability bound \eqref{Linear stability} is fed as an input into a quantitative implicit function theorem (cf. \hyperref[b:IFT 1-DD bound]{(b)} of Lemma~\ref{lmm:Implicit function theorem} and \eqref{stability proof: inverse bound} below). The implicit function theorem then yields the existence of the analytic map $ \bs{D} \mapsto \bs{\frak{G}}(\bs{D}) $ appearing in Theorem~\ref{thr:Stability}.

The following definition will play a crucial role in the upcoming analysis. 
\begin{definition}[Saturated self-energy operator]
Let $\bs{M}=\bs{M}(\1\zeta\1)$ be the solution of the MDE at some spectral parameter $\zeta \in \Cp$. 
We define the linear operator $\cal{F}=\cal{F}(\1\zeta\1): \C^{N \times N} \to \C^{N \times N}$ by
\begin{subequations}
\label{definition of cal F}
\bels{def of F}{
\cal{F}\,:=\, \cal{C}_{\bs{W}}\2\cal{C}_{\!\sqrt{\im \bs{M}}}\,\cal{S}\,\cal{C}_{\!\sqrt{\im \bs{M}}}\,\cal{C}_{\bs{W}}
\,,
}
where we have introduced an auxiliary matrix
\bels{def of W}{
\bs{W}
\,:=\, \pb{\,\bs{1}+\p{\2\cal{C}_{\!\sqrt{\im \bs{M}}}^{\2-1}\1[\1\re \bs{M}\1]\2}^{\12}\2}^{\msp{-2}1/4}
\,.
}
\end{subequations}
We call $\cal{F}$ the \emph{saturated self-energy operator} or the  \emph{saturation} of $\cal{S}$ for short. 
\end{definition}

The operator $\cal{F}$ inherits the self-adjointness with respect to \eqref{matrix scalar product} and the property of mapping $\overline{\mathscr{C}}_+$ to itself from  the self-energy operator  $\cal S$. We will now briefly discuss the reason for introducing $\cal{F}$. In order to invert $\rm{Id}-\cal{C}_{\bs{M}}\cal{S}$ in \eqref{Derivative equation at D eq 0} we have to  show that $\cal{C}_{\bs{M}}\cal{S}$ is dominated by $\rm{Id}$ in some sense. 
 Neither $\cal{S}$ nor $\bs{M}$ can be directly related to the identity operator, but their specific combination $\cal{C}_{\bs{M}}\cal{S}$ can.
We extract this delicate information from the MDE via a Perron-Frobenius argument.  The key observation is that as $\im \zeta \downarrow 0$ the imaginary part of the MDE \eqref{imaginary part of MDE} becomes an eigenvalue equation for the operator $\bs{R} \mapsto \bs{M}^*\cal{S}[\bs{R}]\bs{M}$ with eigenvalue $1$ and corresponding eigenmatrix $\im \bs{M}$. Since this operator is positivity preserving and $\im \bs{M}\in {\mathscr{C}}_+$, its spectral radius is $1$. Naively speaking, through the replacement of $\bs{M}^*$ by $\bs{M}$, the operator $\cal{C}_{\bs{M}}\cal{S}$ gains an additional phase which reduces the spectral radius further and thus guarantees the invertibility of $\rm{Id}-\cal{C}_{\bs{M}}\cal{S}$. However, the non-selfadjointness of the aforementioned operators makes it hard to turn control on their spectral radii into norm-estimates. 
It is therefore essential  to find an appropriate symmetrization of these operators before Perron-Frobenius is applied.
 A similar problem appeared in a simpler commutative setting in \cite{AEK1cpam}. There, $\bs{M}=\rm{diag}(\bs{m})$ was a diagonal matrix and the MDE became a vector equation. In this case the 
problem of inverting $\rm{Id}-\cal{C}_{\bs{M}}\cal{S}$ reduces to inverting a matrix $\bs{1}-\rm{diag}(\bs{m})^2\bs{S}$, where $\bs{S} \in \R^{N \times N}$ is a matrix with non-negative entries that plays the role of the  self-energy  operator $\cal{S}$ in the current setup. The idea in \cite{AEK1cpam} was to write 
\bels{Rewriting derivative in CPAM}{
\bs{1}-\rm{diag}(\bs{m})^2\bs{S}\,=\,\bs{R}\1(\bs{U}-\bs{F})\bs{T}\,,
}
with invertible diagonal matrices $\bs{R}$ and $\bs{T}$, a diagonal unitary matrix $\bs{U}$ and a self-adjoint  matrix  $\bs{F}$, playing the role of the operator $ \cal{F} $,  with positive entries that satisfies the bound $\norm{\bs{F}}\le 1$. It is then possible to see that $\bs{U}-\bs{F}$ is invertible as long as $\bs{U}$ does not leave the Perron-Frobenius eigenvector of $\bs{F}$ invariant. In this commutative setting it is possible to choose $\bs{F}=\rm{diag}(\abs{\bs{m}})\1\bs{S}\2\rm{diag}(\abs{\bs{m}})$, where the absolute value is taken in each component. In our current setting we will achieve a decomposition similar to \eqref{Rewriting derivative in CPAM} on the level of operators acting on $\C^{N \times N}$ (cf. \eqref{expressing stability derivative in terms of cal F} below). The definition \eqref{definition of cal F} ensures that  the saturation  $\cal{F}$ is self-adjoint, positivity-preserving and satisfies $\norm{\cal{F}}\le 1$, as we will establish later. 

\begin{lemma}[Bounds on $\bs{W}$] 
\label{lmm:Bounds on W} Assume  {\bf A1} and $\norm{\bs{A}}\le P_0$ for some constant $P_0>0$. Then uniformly for all spectral parameters $\zeta \in \Cp$ with $\abs{\zeta}\leq 3\1(1+\kappa)$ the matrix $\bs{W}=\bs{W}(\1\zeta\1) \in {\mathscr{C}}_+$, defined in \eqref{def of W}, fulfills the bounds
\bels{bounds on W}{
\norm{\bs{M}}^{-1}\2\bs{1}\,\lesssim\, \rho^{1/2}\2\bs{W}\,\lesssim\, \norm{\bs{M}}^{1/2}\2\bs{1}.
}
\end{lemma}
\begin{Proof}
We write $\bs{W}^{\14}$ in a form that follows immediately from its definition \eqref{def of W},
\[
\bs{W}^4\,=\,
\cal{C}_{\msp{-5}\sqrt{\1\im \bs{M}}}^{\2-1}\2(\2\cal{C}_{\2\im \bs{M}}+\2\cal{C}_{\1\re \bs{M}})[\2(\im\bs{M})^{-1}]\,.
\]
We estimate $(\im\bs{M})^{-1}$ from above and below  by  employing \eqref{im M lower and upper bound} in the regime $\abs{\zeta}\lesssim 1$,
\[
\frac{1}{\rho\2\norm{\bs{M}}^2}\,\cal{C}_{\!\sqrt{\im \bs{M}}}^{\,-1}[\1\bs{M}^*\bs{M}+\bs{M}\bs{M}^*]
\,\lesssim\, 
\bs{W}^{\14}\,\lesssim\, 
\frac{1}{\rho}\,\cal{C}_{\!\sqrt{\im \bs{M}}}^{\,-1}[\1\bs{M}^*\bs{M}+\bs{M}\bs{M}^*]
\,.
\]
Using 
the trivial bounds $2\norm{\bs{M}^{-1}}^{-2}\bs{1}\le\bs{M}^*\bs{M}+\bs{M}\bs{M}^* \le 2\norm{\bs{M}}^2\1\bs{1}$ and 
 $\norm{\bs{M}^{-1}}\lesssim1$ from 
\eqref{M lower bound} 
 as well as \eqref{im M lower and upper bound} again, we find
$
\norm{\bs{M}}^{-4}\1\rho^{-2}\lesssim \bs{W}^4\lesssim\norm{\bs{M}}^2\1\rho^{-2}
$.
This is equivalent to \eqref{bounds on W}.
\end{Proof}

\begin{lemma}[Spectrum of $\cal{F}$] 
\label{lmm:Spectrum of cal F}
Assume {\bf A1} and $\norm{\bs{A}}\le P_0$ for some constant $P_0>0$. 
Then the  saturated self-energy  operator $\cal{F}=\cal{F}(\1\zeta\1)$, defined in \eqref{definition of cal F}, has a unique normalized, $\norm{\bs{F}}_{\rm{hs}}=1$, eigenmatrix $\bs{F}=\bs{F}(\1\zeta\1) \in \mathscr{C}_+$, corresponding to its largest eigenvalue, 
 $ \cal{F}[\1\bs{F}] = \norm{\cal{F}}_{\rm{sp}}\1\bs{F}$. 
Furthermore, the following properties hold uniformly for all spectral parameters $\zeta \in \Cp$  such that 
 $\abs{\zeta}\leq 3\1(1+\kappa)$ and $\norm{\cal{F}(\1\zeta\1)}_{\rm{sp}}\ge 1/2$.
\begin{itemize}
\item[\emph{(i)}] The spectral radius of $\cal{F}$ is given by
\bels{spectral radius of F identity}{
\norm{\cal{F}}_{\rm{sp}}
\,=\, 
1\,- \frac{\scalar{\2\bs{F}}{\1\cal{C}_{\1\bs{W}}[\1\im\bs{M}]\1}}{\avg{\2\bs{F},\!\bs{W}^{-2}}}
\2\im \zeta
\,.
}
\item[\emph{(ii)}] The eigenmatrix $\bs{F}$ is  controlled by the solution of the MDE:  
\bels{bounds on eigenmatrix F}{
\norm{\bs{M}}^{-7}\bs{1}\,\lesssim\,\bs{F}\,\lesssim\,\norm{\bs{M}}^{6}\2\bs{1}\,.
}
\item[\emph{(iii)}] The  operator  $\cal{F}$ has the uniform spectral gap  $ \vartheta \gtrsim \norm{\bs{M}}^{-42}$,  i.e., 
\bels{bound on spectral gap}{
\spec\bigl(\2\cal{F}/\norm{\cal{F}}_{\rm{sp}}\bigr)\,\subseteq\, [-1+\vartheta, 1-\vartheta\2]\cup\{\11\1\}
\,. 
}
\end{itemize}
\vspace{-0.2cm}
\end{lemma}
\begin{Proof} Since $\cal{F}$ preserves the cone $\overline{\mathscr{C}}_+$ of  positive semidefinite  matrices, a version of the Perron-Frobenius theorem for cone preserving operators  implies that there exists a normalized $\bs{F} \in \overline{\mathscr{C}}_+$ such that  $ \cal{F}[\1\bs{F}] = \norm{\cal{F}}_{\rm{sp}}\1\bs{F}$. 
We will show uniqueness of this eigenmatrix later in the proof. First we will prove that \eqref{spectral radius of F identity} holds for any such $\bs{F}$. 
\medskip
\\
\textit{Proof of }(i): 
We define for any matrix $\bs{R} \in \C^{N \times N}$ the {operator} $\cal{K}_{\bs{R}}:\C^{N \times N} \to \C^{N \times N}$ via
\bels{conjugation operator}{
\cal{K}_{\bs{R}}[\bs{T}]\,:=\, \bs{R}^{\msp{-2}*}\bs{T}\1\bs{R}
\,.
}
Note that for self-adjoint $\bs{R}\in \C^{N\times N}$ we have $\cal{K}_{\bs{R}}=\cal{C}_{\bs{R}}$ (cf. \eqref{sandwiching operator}).
Using definition \eqref{conjugation operator}, the imaginary part of the MDE \eqref{imaginary part of MDE} can be written in the form
\bels{imaginary part of MDE using conjugation}{
\im \bs{M}\,=\, (\1\im \zeta\1) \2\cal{K}_{\bs{M}}[\bs{1}]+\cal{K}_{\bs{M}}\cal{S}[\1\im \bs{M}]
\,.
}
We will now write up the equation \eqref{imaginary part of MDE using conjugation} in terms of $\im \bs{M}$, $\cal{F}$ and $\bs{W}$. In order to express $\bs{M}$ in terms of $\bs{W}$, we introduce the unitary matrix 
\bels{definition of U}{
\bs{U}
\,:=\, 
\frac{
\cal{C}_{\msp{-7}\sqrt{\im \bs{M}}}^{\2-1}[\re \bs{M}]-\ii\1\bs{1}
}{
\absb{\2\cal{C}_{\msp{-7}\sqrt{\im \bs{M}}}^{\2-1}[\re \bs{M}]-\ii\1\bs{1}}
}
\,,
}
via the spectral calculus of the self-adjoint matrix $\cal{C}_{\!\sqrt{\im \bs{M}}}^{-1}[\re \bs{M}]$. With \eqref{definition of U} and the definition of $\bs{W}$ from \eqref{def of W} we may write $\bs{M}$ as
\bels{M in terms of U}{
\bs{M}\,=\,\cal{C}_{\!\sqrt{\im \bs{M}}}\2\cal{C}_{\bs{W}}[\bs{U}^*]\,.
}
Here, the matrices  $\bs{W}$ and $\bs{U}$ commute. 
 The identity \eqref{M in terms of U} should be viewed as a \emph{balanced polar decomposition}. Instead of having unitary matrices $\bs{U}_1$ or $\bs{U}_2$ on the left or right of the decompositions $\bs{M}=\bs{U}_1\bs{Q}_1$ or $\bs{M}=\bs{Q}_2\bs{U}_2$, respectively, the unitary matrix $\bs{U}^*$ appears in the middle of $\bs{M}=\bs{Q}^*\bs{U}^*\bs{Q}$ with $\bs{Q}=\bs{W}\sqrt{\im \bs{M}}$. 
Using \eqref{M in terms of U} we also find an expression for $\cal{K}_{\bs{M}}$, namely
\bels{KM in terms of KU}{
\cal{K}_{\bs{M}}\,=\,\cal{C}_{\!\sqrt{\im \bs{M}}}\,\cal{C}_{\1\bs{W}}\2\cal{K}_{\bs{U}^*}\2\cal{C}_{\1\bs{W}}\2\cal{C}_{\!\sqrt{\im \bs{M}}}
\,.
}
Plugging \eqref{KM in terms of KU} into \eqref{imaginary part of MDE using conjugation} and applying the inverse of $\cal{C}_{\!\sqrt{\im \bs{M}}}\2\cal{C}_{\1\bs{W}}\2\cal{K}_{\bs{U}^*}$ on both sides, yields
\bels{imaginary part of MDE with W and F}{
\bs{W}^{-2}\,=\,   \cal{C}_{\2\bs{W}}[\1\im\bs{M}] \2\im \zeta+\cal{F}[\bs{W}^{-2}]
\,,
}
where we used the definition of $\cal{F}$ from \eqref{definition of cal F} and $\cal{K}_{\bs{U}^*}^{-1}[\bs{W}^{-2}]=\bs{W}^{-2}$, which holds because $\bs{U}$ and $\bs{W}$ commute.  We project both sides of \eqref{imaginary part of MDE with W and F} onto the eigenmatrix $\bf{F}$ of $\cal{F}$. Since $\cal{F}$ is self-adjoint with respect to the scalar product \eqref{matrix scalar product} and by  $ \cal{F}[\1\bs{F}] = \norm{\cal{F}}_{\rm{sp}}\1\bs{F}$  
we get
\[
\scalar{\2\bs{F}}{\!\bs{W}^{-2}}
\,=\,
\scalar{\2\bs{F}}{\cal{C}_{\bs{W}}[\im\bs{M}]\1} \2\im \zeta \2+ \norm{\cal{F}\1}_{\rm{sp}}\2\scalar{\2\bs{F}}{\!\bs{W}^{-2}}
\,.
\]
Solving this identity for $\norm{\cal{F}}_{\rm{sp}}$ yields \eqref{spectral radius of F identity}.
\medskip
\\
\textit{Proof of }(ii)\textit{ and }(iii): Let $\zeta \in \Cp$ with $\abs{\zeta}\leq 3\1(1+\kappa)$ and $\norm{\cal{F}(\zeta)}_{\rm{sp}}\ge 1/2$.
The bounds on the eigenmatrix \eqref{bounds on eigenmatrix F} and on the spectral gap \eqref{bound on spectral gap} are a consequence of the  estimate
\bels{connectedness of F}{
\norm{\bs{M}}^{-4} \avg{\1\bs{R}}\1\bs{1}
\,\lesssim\,\cal{F}[\bs{R}]
\,\lesssim\, 
\norm{\bs{M}}^{6} \avg{\bs{R}}\2\bs{1}
\,,\qquad 
\forall\2\bs{R} \in \overline{\scr{C}}_+ 
\,.
}
We verify \eqref{connectedness of F} below.
Given \eqref{connectedness of F}, the remaining assertions, \eqref{bounds on eigenmatrix F} and \eqref{bound on spectral gap}, of Lemma~\ref{lmm:Spectrum of cal F},  are consequences of the following general result that is proven in the appendix.
It generalises to a non-commutative setting the basic fact (cf. Lemma \ref{lmm:basic sp-gap for matrices}) that symmetric matrices with strictly positive entries 
have a positive spectral gap. The proof of Lemma~\ref{lmm:Spectral gap} is given in the appendix.

\begin{lemma}[Spectral gap]
\label{lmm:Spectral gap}
Let $\cal{T}: \C^{N \times N} \to \C^{N \times N}$ be a linear self-adjoint operator preserving the cone $\overline{\scr{C}}_+$ of positive semidefinite matrices. Suppose $\cal{T}$ is normalized, $\norm{\cal{T}}_{\rm{sp}}=1$, and %
\bels{bounds on operator T}{
\gamma\,\avg{\bs{R}}\2\bs{1} \,\leq\, \cal{T}[\bs{R}]
\,\leq\, 
\Gamma \,\avg{\bs{R}}\2\bs{1}\,, \qquad \bs{R} \in \overline{\scr{C}}_+\,,
}
for some positive constants $ \gamma $ and $ \Gamma $.
Then $ \cal{T} $ has a spectral gap of size $ \theta := \frac{\gamma^{\16}}{\12\2\Gamma^{\14}} $, i.e.,
\bels{Spectral gap for T}{
\spec{\cal{T}} \subseteq [-1+\theta,\21-\theta\2] \cup \sett{1} 
\,. 
}
Furthermore, the eigenvalue $1$ is non-degenerate and the corresponding normalized, $\norm{\bs{T}}_{\rm{hs}}=1$, eigenmatrix $\bs{T} \in \scr{C}_+$ satisfies
\vspace{-0.6cm}
\bels{bounds on eigenmatrix T}{
\tsfrac{\gamma}{\!\sqrt{\Gamma\2}}\,
\bs{1}\,\leq\, \bs{T} \,\leq\,\Gamma\, \bs{1}\,.
}
\end{lemma}

Lemma \ref{lmm:Spectral gap} shows the uniqueness of the eigenmatrix $\bs{F}$ as well. In the regime $\abs{\zeta}\ge 3\1(1+\kappa)$ the constants hidden in the comparison relation of \eqref{connectedness of F} will depend on $\abs{\zeta}$, but otherwise the upcoming arguments are not affected. 
In particular the qualitative property of having a unique eigenmatrix $\bs{F}$ remains true even  for large values of $\abs{\zeta}$. 
\medskip 
\\
\textit{Proof of} \eqref{connectedness of F}: 
The bounds in \eqref{connectedness of F} are a consequence of Assumption {\bf A1} and the bounds \eqref{bounds on W} on $\bs{W}$ and \eqref{im M lower and upper bound} on $\im \bs{M}$, respectively. Indeed, from {\bf A1} we have $\cal{S}[\bs{R}]\sim\avg{\bs{R}}\1\bs{1}$ for positive semidefinite matrices $\bs{R}$. By the definition \eqref{def of F} of $\cal{F}$ this immediately yields 
\[
\cal{F}[\bs{R}]
\;\sim\;
\avgb{\2\cal{C}_{\!\sqrt{\im \bs{M}}}\;\cal{C}_{\bs{W}}[\bs{R}]} 
\,
\cal{C}_{\1\bs{W}}\2\cal{C}_{\!\sqrt{\im \bs{M}}}\1[\1\bs{1}\1]
\;=\; 
\scalar{\2\cal{C}_{\bs{W}}[\im \bs{M}]}{\bs{R}\2}\; \cal{C}_{\bs{W}}[\im \bs{M}]
\,.
\]
Since \eqref{bounds on W} and \eqref{im M lower and upper bound} imply
$
\norm{\bs{M}}^{-2}\1\bs{1}\lesssim\cal{C}_{\bs{W}}[\im \bs{M}]\lesssim \norm{\bs{M}}^3\1\bs{1}
$, 
we conclude that \eqref{connectedness of F} holds. 
\end{Proof}
\begin{Proof}[Proof of Proposition \ref{prp:Linear stability}] To show \eqref{Linear stability} we consider the regime of large and small values of $\abs{\zeta}$ separately. We start with the simpler regime, $\abs{\zeta}\ge 3\1(1+\kappa)$. In this case we apply the bound 
$%
\norm{\bs{M}(\1\zeta\1)}
\le 
(\1\abs{\zeta}-\kappa\2)^{-1}
$, %
which is an immediate consequence of the Stieltjes transform representation \eqref{M as Stieltjes transform} of $\bs{M}$. In particular, 
\bels{Bound on C M S}{
\norm{\1\cal{C}_{\bs{M}(\zeta)}\cal{S}\1}_{\rm{sp}}\,\le\, \frac{\norm{\cal{S}}_{\rm{sp}} }{(\abs{\zeta}-\kappa)^2}
\,\le\, \frac{\norm{\cal{S}} }{4\1(1+\kappa)^2}\,\le\,\frac{1}{4}\,,
}
where we used $\kappa\ge \norm{\cal{S}}^{1/2}$ in the last and second to last inequality. We also used that $\norm{\cal{T}}_{\rm{sp}}\le\norm{\cal{T}} $ for any self-adjoint $\cal{T} \in \C^{N \times N}$. The claim \eqref{Linear stability} hence follows in the regime of large $\abs{\zeta}$.

Now we consider the regime $\abs{\zeta}\le 3(1+\kappa)$. Here we will use the spectral properties of the  saturated self-energy  operator $\cal{F}$, established in Lemma~\ref{lmm:Spectrum of cal F}. First we rewrite $\rm{Id}-\cal{C}_{\bs{M}(\zeta)}\cal{S}$ in terms of $\cal{F}$. For this purpose we recall the definition of $\bs{U}$ from \eqref{definition of U}. With the identity \eqref{M in terms of U}
we find 
\bels{cal C M in terms of cal C U}{
\cal{C}_{\2\bs{M}}
\,=\,\cal{C}_{\!\sqrt{\im \bs{M}}}
\,\cal{C}_{\bs{W}}\2\cal{C}_{\bs{U}^*}\2\cal{C}_{\bs{W}}\2\cal{C}_{\!\sqrt{\im \bs{M}}}
\,.
}
Combining \eqref{cal C M in terms of cal C U} with the definition of $\cal{F}$ from \eqref{def of F} we verify
\bels{expressing stability derivative in terms of cal F}{
\rm{Id}-\cal{C}_{\1\bs{M}}\cal{S}
\,=\, 
\cal{C}_{\!\sqrt{\im \bs{M}}}\,\cal{C}_{\bs{W}}\2\cal{C}_{\1\bs{U}^*}(\1\cal{C}_{\bs{U}}-\cal{F}\2) \2 \cal{C}_{\bs{W}}^{-1}\2\cal{C}_{\!\sqrt{\im \bs{M}}}^{\1-1}\,.
}

The bounds \eqref{bounds on W} on $\bs{W}$ and \eqref{im M lower and upper bound} on $\im \bs{M}$ imply bounds on $\cal{C}_{\bs{W}}$ and $\cal{C}_{\!\sqrt{\im \bs{M}}}$, respectively. In fact,  in the regime of bounded $|\zeta|$, we have 
\bels{upper and lower bounds on C W and C im M}{
\norm{\1\cal{C}_{\bs{W}}}\,\lesssim\, \frac{\norm{\bs{M}}}{\rho}
\,,
\quad 
\norm{\2\cal{C}_{\bs{W}}^{\1-1}}\,\lesssim\,\rho \2\norm{\bs{M}}^2
\,,
\quad 
\norm{\2\cal{C}_{\!\sqrt{\im \bs{M}}}\2}\,\lesssim\, \rho\2\norm{\bs{M}}^2
\,,
\quad 
\norm{\2\cal{C}_{\!\sqrt{\im \bs{M}}}^{\1-1}\2}\,\lesssim\,\frac{1}{\rho}
\,.
}
Therefore, taking the inverse and then the norm $\norm{\2\cdot\2}_{\rm{sp}}$ on both sides of \eqref{expressing stability derivative in terms of cal F} and using \eqref{upper and lower bounds on C W and C im M} as well as $\norm{\cal{C}_{\bs{T}}}_{\rm{sp}}\le\norm{\cal{C}_{\bs{T}}} $ for self-adjoint $\bs{T} \in \C^{N \times N}$
yields
\bels{stability derivative bound in terms of B operator}{
\norm{(\rm{Id}-\cal{C}_{\bs{M}}\cal{S})^{-1}}_{\rm{sp}}\,\lesssim\, \norm{\bs{M}}^5\norm{(\cal{C}_{\bs{U}}-\cal{F})^{-1}}_{\rm{sp}}\,.
}
Note that $\cal{C}_{\bs{U}}$ and $\cal{C}_{\bs{U}^*}$ are unitary operators on $\C^{N \times N}$ and thus $\norm{\cal{C}_{\bs{U}}}_{\rm{sp}}=\norm{\cal{C}_{\bs{U}^*}}_{\rm{sp}}=1$. We estimate the norm of the inverse of $\cal{C}_{\bs{U}}-\cal{F}$.
In case $\norm{\cal{F}}_{\rm{sp}}< 1/2$ we will simply use the bound $\norm{(\cal{C}_{\bs{U}}-\cal{F})^{-1}}_{\rm{sp}}\le 2$ in \eqref{stability derivative bound in terms of B operator} and \eqref{M upper bound} for estimating $\norm{\bs{M}}$, thus verifying \eqref{Linear stability} in this   case.

If $\norm{\cal{F}}_{\rm{sp}}\ge 1/2$, we apply the following lemma, which was stated as  Lemma~5.8  in \cite{AEK1cpam}.

\begin{lemma}[Rotation-Inversion Lemma]
\label{lmm:Rotation-Inversion Lemma}
Let $ \cal{T} $ be a self-adjoint and $\cal{U}$ a unitary operator on $ \C^{N \times N} $. 
Suppose that $ \cal{T} $ has a spectral gap, i.e., there is a constant  $ \theta > 0 $  %
such that
\[
\spec{\cal{T}}
\,\subseteq\, 
\bigl[ -\norm{\cal{T}}_{\rm{sp}} \!+ \theta\2,%
\norm{\cal{T}}_{\rm{sp}} \!-\theta\,\bigr] 
\cup \setb{\norm{\cal{T}}_{\rm{sp}}}
\,,
\]
with a non-degenerate largest eigenvalue $ \norm{\cal{T}}_{\rm{sp}} \leq 1 $.
Then there exists a universal positive constant $C$ such that
\[
\norm{\1(\1\cal{U}-\cal{T}\1)^{-1}\1}_{\rm{sp}} 
\,\leq\, 
\frac{C}{\theta}\,\abs{\21 - \norm{\cal{T}}_{\rm{sp}}\2\scalar{\bs{T}}{\1 \cal{U}[\1\bs{T}\1]\1}\1}^{-1}
\,,
\]
where $ \bs{T} $ is the normalized, $\norm{\bs{T}}_{\rm{hs}}=1$, eigenmatrix of $\cal{T}$, corresponding to $\norm{\cal{T}}_{\rm{sp}}$.
\end{lemma}

With the lower bound \eqref{bound on spectral gap} on the spectral gap of $\cal{F}$, we find
\bels{Upper bound on C U - cal F inverse}{
\norm{\1(\1\cal{C}_{\bs{U}}-\cal{F}\1)^{-1}}_{\rm{sp}}
\,&\lesssim\,
\frac{\norm{\bs{M}}^{42 }}{\max\setb{1-\norm{\cal{F}}_{\rm{sp}}\1,\abs{\11-\scalar{\1\bs{F}}{\cal{C}_{\bs{U}}[\1\bs{F}\1]\1}}\2}}
\,.
}
Plugging \eqref{Upper bound on C U - cal F inverse} into \eqref{stability derivative bound in terms of B operator} and using  \eqref{M upper bound} to estimate $\norm{\bs{M}}$, shows \eqref{Linear stability}, provided the denominator on the right hand side of \eqref{Upper bound on C U - cal F inverse} satisfies 
\bels{Lower bound on denominator}{
\max\setb{\21-\norm{\cal{F}}_{\rm{sp}}\1,\2\abs{\21-\scalar{\1\bs{F}}{\cal{C}_{\bs{U}}[\1\bs{F}\1]\1}}\2}
\;\gtrsim\, 
(\2\rho(\1\zeta\1)+\dist(\zeta,\supp\rho\1)\1)^{ C}
,
}
for some universal constant $ C>0$.

In the remainder of this proof we will verify \eqref{Lower bound on denominator}. We establish lower bounds on both arguments of the maximum in \eqref{Lower bound on denominator} and combine them afterwards. We start with a lower bound on $1-\norm{\cal{F}}_{\rm{sp}}$. Estimating the numerator of the fraction on the right hand side of \eqref{spectral radius of F identity} from below 
\[
\scalar{\1\bs{F}}{\cal{C}_{\bs{W}}[\im\bs{M}]\1}
\,\gtrsim\, 
\rho \,
\scalar{\1\bs{F}}{\!\bs{W}^2}\,\gtrsim \,\norm{\bs{M}}^{-2}\avg{\bs{F}}
\,,
\]
and its denominator from above,
$
\scalar{\1\bs{F}}{\!\bs{W}^{-2}}\lesssim \rho\2\norm{\bs{M}}^2\avg{\bs{F}},
$
by applying  the bounds from \eqref{bounds on W} and \eqref{im M lower and upper bound},
we see that 
\bels{first lower bound on 1 - norm F}{
1-\norm{\cal{F}(\1\zeta\1)}_{\rm{sp}}\,\gtrsim\, \frac{\im \zeta}{\rho(\1\zeta\1)\2\norm{\bs{M}(\1\zeta\1)}^4}\,.
}
Since $\rho(\1\zeta\1)$ is the harmonic extension of a probability density (namely the self-consistent density of states $\rho$), we have the trivial upper bound $\rho(\1\zeta\1)\lesssim \im \zeta /\dist(\zeta, \supp \rho)^2$. Continuing from \eqref{first lower bound on 1 - norm F} we find the lower bound 
\bels{lower bound on 1 - norm F}{
1-\norm{\cal{F}}_{\rm{sp}}
\,\gtrsim\,
\norm{\bs{M}}^{-4}\dist (\zeta, \supp \rho)^2
\,\gtrsim\,
(\1\rho+\dist(\zeta,\supp\rho)\1)^{4}\dist (\zeta, \supp \rho)^2
\,,
}
where we used \eqref{M upper bound} in the second inequality.

Now we estimate $\abs{\21-\scalar{\1\bs{F}}{\cal{C}_{\bs{U}}[\1\bs{F}\1]\1}}$ from below. 
We begin with 
\bels{First lower bound on 1- scalar F C U F}{
\!\abs{\21-\scalar{\1\bs{F}}{\cal{C}_{\bs{U}}[\bs{F}]\1}}
\ge
\re\bigl[ \21-\scalar{\1\bs{F}}{\cal{C}_{\bs{U}}[\bs{F}]\1}\2\bigr]
\,=\,  
1-\avgb{\1\bs{F},(\2\cal{C}_{\re \bs{U}}-\cal{C}_{\im \bs{U}})[\bs{F}]\1}
\ge 
\scalar{\2\bs{F}}{\cal{C}_{\im \bs{U}}[\bs{F}]}
\,,
}
where $1- \scalar{\bs{F}}{\cal{C}_{\re \bs{U}}\bs{F}}\ge 0$ in the last inequality, because $\bs{U}$ is unitary and 
$\norm{\bs{F}}_{\rm{hs}}=1$. 
Since $\im \bs{U} =- \bs{W}^{-2}$ (cf. \eqref{definition of U} and \eqref{definition of cal F}) 
and because of \eqref{bounds on W} we have 
$
-\im \bs{U}\gtrsim \norm{\bs{M}}^{-2}\rho
\,.
$
Continuing from \eqref{First lower bound on 1- scalar F C U F}, using the normalization $\norm{\bs{F}}_{\rm{hs}}=1 $  and  \eqref{M upper bound},  we get the lower bound 
\[
\abs{\21\1-\scalar{\1\bs{F}}{\cal{C}_{\bs{U}}[\1\bs{F}\1]\1}}\,\gtrsim\, 
\rho^{\22}\2\norm{\bs{M}}^{-4}
\,\gtrsim\;
(\2\rho+\dist(\zeta,\supp\rho)\1)^{4}
\rho^{\12}
\,.
\]
Combining this with \eqref{lower bound on 1 - norm F} shows \eqref{Lower bound on denominator} and thus finishes the proof of Proposition~\ref{prp:Linear stability}.
\end{Proof}
\begin{Proof}[Proof of Proposition \ref{prp:Density of states}]
We show that the harmonic extension $\rho(\1\zeta\1)$ of the self-consistent density of states (cf. \eqref{definition extended DOS}) is uniformly  $c$-H\"older  continuous on the entire complex upper half plane. Thus its unique continuous extension to the real line, the self-consistent density of states, inherits this regularity. 

We differentiate both sides of the MDE with respect to $\zeta$ and find the equation
\bels{M derivative}{
(\1\rm{Id}-\cal{C}_{\bs{M}}\cal{S}\1)[\partial_\zeta\bs{M}]
=
\bs{M}^{\12}
\,.
}
Inverting the operator $\rm{Id}-\cal{C}_{\bs{M}}\cal{S}$ and taking the normalized Hilbert-Schmidt norm reveals a bound on the derivative of the solution to the MDE,
\bels{bound on partial M}{
\norm{\partial_\zeta\bs{M}\1}_{\rm{hs}}
\,\leq\,
\norm{\1 (\1\rm{Id}-\cal{C}_{\bs{M}}\cal{S}\1)^{-1}}_{\rm{sp}} \norm{\bs{M}}^2
.
}
Since $\zeta \to \avg{\bs{M}(\1\zeta\1)}$ is an analytic function on $\Cp$, we have the basic  identity
$2\pi\ii\1\partial_\zeta\rho=2\1\ii\1\partial_\zeta\im\avg{\bs{M}}\,=\, \partial_\zeta\avg{\bs{M}}$. 
Therefore, making use of \eqref{bound on partial M}, we get
\bels{bound on partial rho}{
\abs{\1\partial_\zeta\rho\2}
\,=\, 
{\textstyle \frac{1}{\12\1\pi\1}}\abs{\avg{\partial_\zeta\bs{M}\1}}
\,\leq\,  
{\textstyle  \frac{1}{\12\1}}\norm{\1(\1\rm{Id}-\cal{C}_{\bs{M}}\cal{S}\1)^{-1}}_{\rm{sp}} 
\norm{\bs{M}}^2
\,\lesssim\, 
\rho^{-(C+2)}
\,.
}
For the last inequality in \eqref{bound on partial rho} we employed the bound \eqref{M upper bound} and the linear stability, Proposition~\ref{prp:Linear stability}. The universal constant $C$ stems from its statement \eqref{Linear stability}.
From \eqref{bound on partial rho} we read off that the harmonic extension $\rho$ of the self-consistent density of states is  $\frac{1}{C+3}$-H\"older continuous.
 
It remains to prove that $ \rho $ is real analytic at any $\tau_0$ with $\rho(\tau_0)>0$. Since $\rho$ is continuous, it is bounded away from zero in a neighborhood of $\tau_0$. Using \eqref{M derivative}, \eqref{M upper bound} and \eqref{Linear stability} we conclude that $\bs{M}$ is uniformly continuous in the intersection of a small neighborhood of $\tau_0$ in $\C$ with the complex upper half plane. In particular, $\bs{M}$ has a unique continuous extension $\bs{M}(\tau_0)$ to $\tau_0$. Furthermore, by differentiating \eqref{MDE} with respect to $\zeta$ and by the uniqueness of the solution to \eqref{MDE} with positive imaginary part one verifies that $\bs{M}$ coincides with the solution $\bs{Q}$ to the holomorphic initial value problem
\[
\partial_\omega \bs{Q}\,=\,
(\rm{Id}-\cal{C}_{\bs{Q}}\1\cal{S})^{-1} [\bs{Q}^2]\,,\qquad \bs{Q}(0)\,=\, \bs{M}(\tau_0)\,,
\]
i.e. $\bs{M}(\tau_0+\omega)=\bs{Q}(\omega)$ for any $\omega \in \Cp$ with sufficiently small absolute value. Since the solution $\bs{Q}$ is analytic in a small neighborhood of zero, we conclude that $\bs{M}$ can be holomorphically extended to a neighborhood of $\tau_0$ in $\C$. By continuity \eqref{definition extended DOS} remains true for $\zeta \in \R$ close to $\tau_0$ and thus $\rho$ is real analytic there. 
\end{Proof}

In the proof of Theorem~\ref{thr:Stability} we will often consider $\cal T:(\C^{N \times N}, \norm{\2\cdot\2}_A)\to (\C^{N \times N}, \norm{\2\cdot\2}_B)$, i.e., $\cal T$ is a linear operator on $\C^{N \times N}$ equipped with two different norms. We indicate the norms in the notation of the corresponding induced operator norm $\norm{\cal{T}}_{A \to B}$. We will use $A,B=\rm{hs}, \norm{\2\cdot\2}, 1, \infty,  \rm{max}$, etc.
We still keep our convention that 
$
\norm{\cal{T}}_{\rm{sp}}=\norm{\cal{T}}_{\rm{hs}\to \rm{hs}} $ and $ \norm{\cal{T}}=\norm{\cal{T}}_{\norm{\1\cdot\1}\to\norm{\1\cdot\1}}
$.
Furthermore, we introduce the norms 
\bels{ell-1, ell-infty and max-ell-1-infty matrix norms}{
\norm{\bs{R}}_1 := \max_y{ \textstyle \sum_x}\,\abs{\1r_{x y}}
\,,\quad
\norm{\bs{R}}_\infty := \max_x {\textstyle \sum_y}\, \abs{\1r_{xy}}
\,,\quad
\norm{\bs{R}}_{1\vee\infty}
\,:=\, 
\max\setb{\2\norm{\bs{R}}_1,\norm{\bs{R}}_\infty }
\,.
}
Some of the norms on matrices $ \bs{R} \in \C^{N \times N}$ are ordered, e.g. $\max\{\2\norm{\bs{R}}_\rm{max},\norm{\bs{R}}_{\rm{hs}}\}\le \norm{\bs{R}}\le\norm{\bs{R}}_{1\vee\infty}$. 
Note that if $\norm{\2\cdot\2}_{\wt A}\le \norm{\2\cdot\2}_{ A}$ and $\norm{\2\cdot\2}_{B}\le \norm{\2\cdot\2}_{ \wt B}$, then $\norm{\2\cdot\2}_{A \to B}\le \norm{\2\cdot\2}_{\wt A \to \wt B}$. In particular, for $\cal T:\C^{N \times N}\to \C^{N \times N}$ we have e.g. 
$ 
\norm{\cal{T}}_{\rm{max}\to \rm{hs}} 
\le 
\norm{\cal{T}}_{\rm{max}\to \norm{\1\cdot\1}} 
\le 
\norm{\cal{T}}_{\rm{max}\to 1 \vee \infty}
$.

In order to  show the existence and properties of the map $ \bs{D} \mapsto \bs{\frak{G}}(\bs{D}) $ from Theorem~\ref{thr:Stability}  we rely on an implicit function theorem,  which we state here for reference purposes. 

\begin{lemma}[Quantitative implicit function theorem]
\label{lmm:Implicit function theorem}
Let $ T: \C^{A}\times \C^{D} \to\C^{A}$ be a continuously differentiable function with invertible derivative $\nabla^{(1)}T(0,0)$ at the origin with respect to the first argument and ${T}(0,0)=0$. 
Suppose $\C^{A}$ and $\C^{D}$ are equipped with norms that we both denote by $\norm{\2\cdot\2}$, and let the linear operators on these spaces be equipped with the corresponding induced operator norms.
Let $\delta > 0 $ and $ C_1,C_2 < \infty $ be constants, such that 
\begin{itemize}
\titem{a}\label{a:IFT lower bound on D_1T at 0} 
$ \norm{\1(\nabla^{(1)}T(0,0))^{-1}} \leq C_1 $;
\titem{b}\label{b:IFT 1-DD bound}
$ \normb{\,\rm{Id}_{\C^{A}}-(\nabla^{(1)}T(0,0))^{-1}\nabla^{(1)}T(a,d\1)\2} \leq \frac{1}{2} $, for every $ (a,d\1) \in B^A_\delta \times B^D_\delta $;
\titem{c}\label{c:IFT upper bound on nablaD of T} 
$ \norm{\nabla^{(2)}T(a,d\1)} \le C_2 $, for every $ (a,d\1) \in B^A_\delta \times B^D_\delta $.
\end{itemize}
Here $B_\delta^\#$ is the $\delta$-ball around $0$ with respect to $\norm{\2\cdot\2}$ in $\C^\#$, and $\nabla^{(2)}$ denotes the derivative with respect to the second variable.   

Then there exists a constant $ \eps > 0 $, depending only on $\delta$, $C_1$ and $C_2$, and a unique continuously differentiable function $ f : B^D_\eps \to B^A_\delta $, such that $ T(f(d),d\1) = 0 $, for every $ d \in B^D_\eps $.
Furthermore, if $T$ is analytic, then so is $f$. 
\end{lemma}
The proof of this result is elementary and left to the reader.

\begin{Proof}[Proof of Theorem~\ref{thr:Stability}] 
To apply Lemma~\ref{lmm:Implicit function theorem} we define  $\cal{J}:\C^{N \times N}\times \C^{N \times N}\to \C^{N \times N}$ by
\[
\cal{J}[\bs{\frak{G}},\bs{D}]\,:=\,  \bs{1} \2+\2 (\2\zeta\2\bs{1}-\bs{A}+\cal{S}[\bs{\frak{G}}]\1)\1\bs{\frak{G}}+\bs{D}
\,.
\]
With this definition the perturbed MDE \eqref{perturbed MDE2} takes the form
\[
\cal{J}[\bs{\frak{G}}(\bs{D}),\bs{D}]\,=\,  \bs{0}  %
\,.
\]
 In particular,  the unperturbed MDE \eqref{MDE} is $\cal{J}[\bs{M},\bs{0}\1]=  \bs{0}$, with  $\bs{M}=\bs{\frak{G}}(\bs{0})$. 

For the application of the implicit function theorem we control the derivatives of $\cal{J}$ with respect to $\bs{\frak{G}}$ and $\bs{D}$. With the short hand notation,
\bels{definition of cal W}{
\cal{W}_{\bs{R}}[\bs{T}]
\,:=\,
\bs{M}\1(\1\cal{S}[\bs{T}]\1\bs{R}\2+\cal{S}[\bs{R}]\1\bs{T}\2)
\,,
}
we compute the directional derivative  of $\cal{J}$  with respect to $\bs{\frak{G}}$ in the direction $\bs{R} \in \C^{N \times N}$, 
\bels{derivatives wrt G}{
\nabla^{(\bs{\frak{G}})}_{\!\bs{R}}\!\cal{J}[\bs{\frak{G}},\bs{D}]
\,&=\, (\1\zeta\2\bs{1}-\bs{A}+\cal{S}[\bs{\frak{G}}]\1)\1\bs{R}+\cal{S}[\bs{R}]\1\bs{\frak{G}}
\\
&=\, -\2\bs{M}^{-1}(\1\rm{Id}-\cal{C}_{\bs{M}}\cal{S}-\cal{W}_{\bs{\frak{G}}-\bs{M}})[\bs{R}]
\,.
}
For the second identity in \eqref{derivatives wrt G} we used \eqref{MDE}. 

The derivative with respect to $\bs{D}$ is simply the identity operator,
$
\nabla^{(\bs{D})}\!\cal{J}[\bs{\frak{G}},\bs{D}]= \rm{Id}
$.
Therefore, estimating $\nabla^{(\bs{D})}\cal{J}$  for the hypothesis \hyperref[c:IFT upper bound on nablaD of T]{(c)} of Lemma~\ref{lmm:Implicit function theorem} is  trivial.

We consider $ \C^{N \times N} \cong \C^{N^2}$ with the entrywise maximum norm $\norm{\2\cdot\2}_\rm{max}$ and use the short hand notation $\norm{\cal{T}}_\rm{max}:=\norm{\cal{T}}_{\rm{max}\to \rm{max}}$ for the induced operator norm of any linear $\cal{T}: \C^{N \times N}\to \C^{N \times N}$. 
To apply Lemma~\ref{lmm:Implicit function theorem}  in this setup  we need the following two estimates: 
\begin{itemize}
\item[(i)] 
The operator norm of  $ (\2\rm{Id}-\cal{C}_{\bs{M}}\cal{S}\2)^{-1} $ on $(\C^{N \times N}, \norm{\2\cdot\2}_\rm{max})$ is controlled by its spectral norm,
\bels{0 to 0 bound on stability derivative}{
\norm{(\1\rm{Id}-\cal{C}_{\bs{M}}\cal{S}\1)^{-1}}_{\rm{max}}
\,\lesssim\; 
1 \,+\norm{\bs{M}}^2+\norm{\bs{M}}^4\norm{(\1\rm{Id}-\cal{C}_{\bs{M}}\cal{S}\1)^{-1}}_{\rm{sp}}
\,.
}
\item[(ii)] 
The operator norm of $\cal{W}_{\bs{\frak{G}}-\bs{M}}$ is small, provided $\bs{\frak{G}}$ is close to $\bs{M}$,
\bels{0 to 0 bound on W G - M}{
\norm{\1\cal{W}_{\bs{\frak{G}}-\bs{M}}}_{\rm{max}}\lesssim \, \norm{\bs{M}}_{1\vee \infty} \norm{\bs{\frak{G}}-\bs{M}}_{\rm{max}}\,.
}
\end{itemize} 
We will prove these estimates after we have used them to show that the hypotheses of the quantitative inverse function theorem hold. 

Let us first bound the operator $ \bs{R}\mapsto \nabla^{(\bs{\frak{G}})}_{\!\bs{R}}\!\cal{J}[\bs{M},\bs{0}\1] $. 
To this end, using \eqref{derivatives wrt G} we have
\bels{stability proof: inverse bound}{
\norm{\1
(\nabla^{(\bs{\frak{G}})}\!\cal{J}[\bs{M},\bs{0}\1]\1)^{-1}\msp{-1}[\bs{R}]\1}_\rm{max} 
&\leq\,\norm{(\1\rm{Id}-\cal{C}_{\bs{M}}\cal{S}\1)^{-1}}_{\rm{max}}\,
\norm{\bs{M}}_{1\vee \infty}
\norm{\bs{R}}_{\rm{max}}
\,,
}
for an arbitrary $ \bs{R} $. For the last line we have used $ \norm{\bs{M}\1\bs{R}}_{\rm{max}} \leq \norm{\bs{M}}_{1\vee \infty} \norm{\bs{R}}_{\rm{max}}$. 
By Theorem~\ref{thr:Arbitrarily high polynomial decay of solution} there is a sequence $\ul{\gamma}$, depending only on $\delta$ and $\scr{P}$, such that 
\bels{bound on 1 infty norm of M}{
\norm{\bs{M}}_{1\vee \infty}
\le\, 
\max_x \sum_y \frac{\gamma(\nu)}{(1+d(x,y))^{\nu}} +\gamma(0)
\,, 
\qquad \nu \in \N\,.
}
Here and in the following unrestricted summations $\sum_x$ are understood to run over the entire index set from $1$ to $N$.  %
Since the sizes of  the  balls with respect to $d$ grow only polynomially in their radii (cf. \eqref{polynomial ball growth}), the right hand side of \eqref{bound on 1 infty norm of M} is bounded by a constant that only depends on $\delta$ and $\scr{P}$ for a sufficiently large choice of $\nu$.
Using this estimate together with the bound (i) for the inverse of $ \rm{Id}-\cal{C}_{\bs{M}}\cal{S} $ in \eqref{stability proof: inverse bound} yields the bound $ \norm{\1 (\nabla^{(\bs{\frak{G}})}_{\!\bs{R}}\!\cal{J}[\bs{M},\bs{0}\1]\1)^{-1}}_\rm{max} \lesssim 1 $  for \hyperref[a:IFT lower bound on D_1T at 0]{(a)} of Lemma~\ref{lmm:Implicit function theorem}. 

Next, to verify  the assumption \hyperref[b:IFT 1-DD bound]{(b)} 
of Lemma~\ref{lmm:Implicit function theorem} we write
\bels{stability-1}{
\rm{Id}-(\nabla^{(\bs{\frak{G}})}\!\cal{J}[\bs{M},\bs{0}])^{-1}\2\nabla^{(\bs{\frak{G}})}\!\cal{J}[\bs{\frak{G}},\bs{D}]
\,=\, 
(\1\rm{Id}-\cal{C}_{\bs{M}}\cal{S}\1)^{-1}\cal{W}_{\bs{\frak{G}}-\bs{M}}
\,.
}
Using \eqref{0 to 0 bound on W G - M} and \eqref{0 to 0 bound on stability derivative}, in conjunction with \eqref{M upper bound} and \eqref{Linear stability}, we see that 
\bels{the key bound for max-stability}{
\normb{\2\rm{Id}\2-\2(\nabla^{(\bs{\frak{G}})}\!\cal{J}[\bs{M},\bs{0}\1]\1)^{-1}\2\nabla^{(\bs{\frak{G}})}\!\cal{J}[\bs{\frak{G}},\bs{D}]\2}_{\rm{max} }
\le\, 
\frac{1}{2}
\,, 
}
for all $ (\bs{\frak{G}},\bs{D}) \in B^{\rm{max}}_{c_2}(\bs{M}) \times B^{\rm{max}}_{c_1}(\bs{0}) $, provided $ c_1,c_2 \sim 1 $ are sufficiently small.
The first part of Theorem~\ref{thr:Stability}, the existence and uniqueness of the analytic function $\bs{\frak{G}}$, now follows from the implicit function theorem Lemma~\ref{lmm:Implicit function theorem}. 
In particular, \eqref{MDE stability2} follows from the analyticity. 
\medskip
\\
\textit{Proof of} (i): 
First we remark that \eqref{S operator decay}  for a large enough $\nu$ and together with \eqref{polynomial ball growth} imply 
\bels{max to 1 infty norm of S}{
\norm{\cal{S}}_{\rm{max} \to 1\vee \infty}\,\lesssim\, 1\,.
}

We expand the geometric series corresponding to the operator $(\rm{Id}-\cal{C}_{\bs{M}}\cal{S})^{-1}$ to second order,
\bels{Expand geometric series in C1}{
(\1\rm{Id}-\cal{C}_{\bs{M}}\cal{S}\1)^{-1}
\,=\, {\rm{Id}+\cal{C}_{\bs{M}}\cal{S}+\frac{(\cal{C}_{\bs{M}}\cal{S})^2}{\rm{Id}-\cal{C}_{\bs{M}}\cal{S}}}\,.
}
We consider each of the three terms on the right hand side separately and estimate
 their norms as operators from $\C^{N \times N}$ with the entrywise maximum norm to itself. 

The easiest is $\norm{\rm{Id}}_\rm{max}=1$.  For the second term we use the estimate
\bels{Second term estimate in C1}{
\norm{\1\cal{C}_{\bs{M}}\cal{S}}_{\rm{max}}\,\le\,\norm{\1\cal{C}_{\bs{M}}\cal{S}}_{\rm{max} \to \norm{\1\cdot\1}}
\,\le\, \norm{\cal{C}_{\bs{M}}}\norm{\cal{S}}_{\rm{max} \to \norm{\1\cdot\1}}\,.
}
 For the third term on the right hand side of \eqref{Expand geometric series in C1} we apply
\bels{Third term estimate in C1}{
\normB{\frac{(\cal{C}_{\bs{M}}\cal{S})^2}{\rm{Id}-\cal{C}_{\bs{M}}\cal{S}}}_{\rm{max}}
\,\le\, 
\norm{\1\cal{C}_{\bs{M}}\cal{S}}_{\rm{hs} \to \rm{max}}\2 
\norm{(\1\rm{Id}-\cal{C}_{\bs{M}}\cal{S}\1)^{-1}}_{\rm{sp}} 
\2\norm{\1\cal{C}_{\bs{M}}\cal{S}}_{\rm{max} \to \rm{hs}}\,.
}
The last factor on the right hand side of \eqref{Third term estimate in C1} is bounded by 
\bels{Last factor bound in third term in C1}{
\norm{\1\cal{C}_{\bs{M}}\cal{S}}_{ \rm{max}\to \rm{hs} }
\,\le\, 
\norm{\1\cal{C}_{\bs{M}}\cal{S}}_{ \rm{max}\to \norm{\1\cdot\1} }
\,\le\,
\norm{\1\cal{C}_{\bs{M}}}\2\norm{\1\cal{S}}_{\rm{max}\to \norm{\1\cdot\1}}\,.
}
For the first factor we use 
$
\norm{\cal{C}_{\bs{M}}\cal{S}}_{ \rm{hs} \to\rm{max}} 
\le \norm{\cal{C}_{\bs{M}}\cal{S}}_{ \rm{hs} \to \norm{\1\cdot\1}}
\le \norm{\cal{C}_{\bs{M}}} \norm{\cal{S}}_{\rm{hs} \to \norm{\1\cdot\1}} $.
We plug this 
and \eqref{Last factor bound in third term in C1} into \eqref{Third term estimate in C1}. Then we use the resulting inequality in combination with \eqref{Second term estimate in C1} in \eqref{Expand geometric series in C1} and find
\[
\norm{(\rm{Id}-\cal{C}_{\bs{M}}\cal{S})^{-1}}_{\rm{max}}\,\lesssim\, 1+\norm{\cal{C}_{\bs{M}}}
+\norm{\cal{C}_{\bs{M}}}^2\norm{(\rm{Id}-\cal{C}_{\bs{M}}\cal{S})^{-1}}_{\rm{sp}}\,,
\]
where we also used $\norm{\cal{S}}_{\rm{max} \to \norm{\1\cdot\1}}\le \norm{\cal{S}}_{\rm{max} \to 1\vee \infty}$ and \eqref{max to 1 infty norm of S}.
Since $\norm{\cal{C}_{\bs{M}}}\le \norm{\bs{M}}^2$ the claim \eqref{0 to 0 bound on stability derivative} follows.
\medskip
\\
\textit{Proof of} (ii):  Recall the definition of $\cal{W}_{\bs{R}}$ in \eqref{definition of cal W}. We estimate
\begin{equation}
\begin{split}
\norm{\1\cal{W}_{\1\bs{R}}[\bs{T}]\1}_\rm{max}&\le\, 2\norm{\bs{M}}_{1\vee \infty}\norm{\cal{S}}_{\rm{max}\to 1\vee \infty}\norm{\bs{R}}_\rm{max}\norm{\bs{T}}_\rm{max}\,.
\end{split}
\end{equation}
From the bound \eqref{max to 1 infty norm of S} we infer \eqref{0 to 0 bound on W G - M}.
\medskip
\\
\textit{Proof of} \eqref{derivative representation with Z} \textit{and} \eqref{bounds on Z}:
Now we are left with showing the second part of  Theorem~\ref{thr:Stability}, namely that the derivative of $\bs{\frak{G}}$ at $\bs{D}=\bs{0}$ can be written in the form \eqref{derivative representation with Z}  with the operator $\cal{Z}$ satisfying  \eqref{bounds on Z}.

Since we have shown the analyticity of $\bs{\frak{G}}$, the calculation leading up to \eqref{Derivative equation at D eq 0} is now justified and we see that
\[
\nabla_{\bs{R}}\bs{\frak{G}}(\bs{0})\,=\, (\rm{Id}-\cal{C}_{\bs{M}}\cal{S})^{-1}[\bs{M}\bs{R}]
\,=\,
\bs{M}\bs{R}
+\cal{Z}[\bs{R}]
\,,
\]
for all $\bs{R} \in \C^{N \times N}$.
Here, the linear operator $\cal{Z}$ is given by
\bels{definition of cal Z}{
\cal{Z}[\bs{R}] \,&:=\, \frac{\cal{C}_{\bs{M}}\cal{S}}{\1\rm{Id}-\cal{C}_{\bs{M}}\cal{S}\1}[\1\bs{M}\1\bs{R}]
\,=\, 
\pB{\2\cal{C}_{\bs{M}}\cal{S}+(\cal{C}_{\bs{M}}\cal{S}\1)^2+(\cal{C}_{\bs{M}}\cal{S}\1)^2(\1\rm{Id}-\cal{C}_{\bs{M}}\2\cal{S}\2)^{-1}\1\cal{C}_{\bs{M}}\cal{S}\2}[\1\bs{M}\1\bs{R}\1]
\,.
}
We will estimate the entries of the three summands separately. 

We show that $\norm{\cal{Z}[\bs{R}]}_{\ul{\gamma}}\le \frac{1}{2}$  for any $\bs{R} \in \C^{N \times N}$ with $\norm{\bs{R}}_{\rm{max}}\le 1$, where $\ul{\gamma}$ depends only on $\delta$ and $\scr{P}$.  
We begin with a few easy observations:
For two matrices $\bs{R}, \bs{T} \in \C^{N \times N}$ that have faster than power law decay, $\norm{\bs{R}}_{\ul{\gamma}_R}\le 1$ and $\norm{\bs{T}}_{\ul{\gamma}_T}\le 1$, their sum and product have faster than power law decay as well, i.e., $\norm{\bs{R}+\bs{T}}_{\ul{\gamma}_{R+T}}\le 1$ and $\norm{\bs{R}\bs{T}}_{\ul{\gamma}_{RT}}\le 1$. Here, $\ul{\gamma}_{R+T}$ and $\ul{\gamma}_{RT}$  depend only on $\ul{\gamma}_R,\ul{\gamma}_T$ and $P$ (cf. \eqref{polynomial ball growth}).
Furthermore, we see that by \eqref{S operator decay} the matrix $\cal{S}[\bs{R}]$ has faster than power law decay for any $\bs{R}\in \C^{N \times N}$ with $\norm{\bs{R}}_{\rm{max}}\le 1$. 

By the following argument we estimate the first summand on the right hand side of \eqref{definition of cal Z}. Using \eqref{S operator decay}, $\norm{\bs{M}\bs{R}}_{\rm{max}}\le \norm{\bs{M}}_{1\vee \infty}\norm{\bs{R}}_{\rm{max}}$ and the estimate \eqref{bound on 1 infty norm of M}, the matrix $\cal{S}[\bs{M}\bs{R}]$ has faster than power law decay. 
Since $\cal{C}_{\bs{M}}$ multiplies with $\bs{M}$ on both sides (cf. \eqref{sandwiching operator}) and $\bs{M}$ has faster than power law decay (cf. Theorem~\ref{thr:Arbitrarily high polynomial decay of solution}), we conclude that so has $\cal{C}_{\bs{M}}\cal{S}[\bs{M}\bs{R}]$. 

Now we turn  to the second summand on the right hand side of \eqref{definition of cal Z}. Since $\cal{C}_{\bs{M}}\cal{S}[\bs{M}\bs{R}]$ has faster than power law decay, its entries are bounded.
  Using again \eqref{S operator decay} as above, we see that $\cal{C}_{\bs{M}}\cal{S}$ applied to $\cal{C}_{\bs{M}}\cal{S}[\bs{M}\bs{R}]$ has
 faster than power law decay as well.  

Finally, we estimate the third summand from \eqref{definition of cal Z}. 
Since the matrix $\cal{C}_{\bs{M}}\cal{S}[\bs{M}\bs{R}]$ has faster than power law decay, its $\norm{\2\cdot\2}_{\rm{hs}}$-norm is bounded. 
By the linear stability \eqref{Linear stability}  and $\zeta\in \D_\delta$, we conclude 
$
\norm{\2(\1\rm{Id}-\cal{C}_{\bs{M}}\cal{S}\1)^{-1}[\2\cal{C}_{\bs{M}}\cal{S}[\bs{M}\bs{R}]\1]\1}_{\rm{hs}}\,\le\, C(\delta)
$. 
Thus, we get
\[
\norm{\,\cal{C}_{\bs{M}}\1\cal{S}\,(\1\rm{Id}-\cal{C}_{\bs{M}}\cal{S}\1)^{-1}[\2\cal{C}_{\bs{M}}\cal{S}[\bs{M}\bs{R}]\1]}_{\rm{max}}
\le\,
C(\delta)\2\norm{\1\cal{S}\1}_{\rm{hs}\to \rm{max}}\norm{\1\cal{C}_{\bs{M}}}_{\rm{max}}
\le\,
C(\delta)\2\norm{\1\cal{S}\1}_{\rm{hs}\to 1\vee \infty}
\norm{\bs{M}}_{1\vee \infty}^2
\,,
\]
which is bounded by \eqref{max to 1 infty norm of S} and \eqref{bound on 1 infty norm of M}. Therefore, the third term on the right hand side of \eqref{definition of cal Z} is an application of $\cal{C}_{\bs{M}}\cal{S}$ to a matrix with bounded entries, which results in a matrix with faster than power law decay. Altogether we have established that \eqref{bounds on Z} hold with only $\norm{\cal{Z}[\bs{R}]}_{\ul{\gamma}}$ on the left hand side. 

It remains to show that also $\cal{Z}^*[\bs{R}]$ satisfies this bound. Since $\cal{Z}^*$ has a structure that resembles the structure \eqref{definition of cal Z} of $\cal{Z}$, namely
\[
\cal{Z}^*[\bs{R}]\,=\, 
\bs{M}^* \pB{\cal{S}\2\cal{C}_{\bs{M}^*}+(\cal{S}\2\cal{C}_{\bs{M}^*})^2+(\cal{S}\2\cal{C}_{\bs{M}^*})^2\big(\2\rm{Id}-(\1\cal{C}_{\bs{M}}\cal{S}\1)^*\big)^{-1}\cal{S}\2\cal{C}_{\bs{M}^*}\msp{-2}}[\1\bs{R}\1]
\,,
\]
we can follow the same line of reasoning as for the entries of $\cal{Z}[\bs{R}]$.
This finishes the proof of  \eqref{bounds on Z} and with it the proof of Theorem~\ref{thr:Stability}.
\end{Proof}

\section{Estimating the error term}
\label{sec:Estimating the error term}

In this section we prove the  key estimates, stated precisely in Lemmas~\ref{lmm:Smallness of error matrix} and \ref{lmm:Smallness of error matrix away from support}, for the error matrix $ \bs{D} $ that appears as the perturbation in the equation \eqref{D-perturbed MDE and D} for the resolvent $ \bs{G}$. 
We start  by estimating $\bs{D}(\zeta) $ in terms of the auxiliary quantity $\Lambda(\zeta) $ (cf.  \eqref{definition of Lambda}) when $ \zeta $ is away from the convex hull of $\supp \rho $.
To this end, we recall the two endpoints of this convex hull (cf. Proposition~\ref{prp:fluctuation averaging}):
\bels{definition kappa-+}{
\kappa_-\,:=\,\min \supp \rho\,, \qquad \kappa_+\,:=\,\max \supp \rho
\,.
}

\begin{lemma}
\label{lmm:Smallness of error matrix away from support} Let $\delta>0$  and $\eps>0$. 
Then  the error matrix $\bs{D}$, defined in \eqref{definition of error matrix D}, satisfies
\bels{error bound away from convex hull of support}{
\norm{\bs{D}(\1\zeta\1)}_{\rm{max}}\,\bbm{1}(\1\Lambda(\zeta\1)\le  N^{-\eps})\,\prec\,\frac{1}{\sqrt{N }}+\pB{\frac{\Lambda(\1\zeta\1)}{N \im \zeta}}^{1/2},
}
for all $\zeta \in \Cp$ with $ \delta \le\dist(\1\zeta, [\kappa_-,\kappa_+])\le \delta^{-1}$ and $\im \zeta \ge N^{-1+\eps}$.
\end{lemma}

\begin{convention}
Throughout this section we will use Convention~\ref{con:Comparison relation} with the set of model parameters $\scr{P}$ replaced by the set $\scr{K}$ from \eqref{definition of model parameters K}. If the constant $C$, hidden in the comparison relation, depends on additional parameters $\scr{L}$, then we write $\alpha \lesssim_\scr{L} \beta$.
\end{convention}
We rewrite the entries $d_{xy}$ of $\bs{D}$ in a different form, that allows  us  to see their smallness, by expanding the term $(\bs{H}-\bs{A})\bs{G}$ (cf. \eqref{definition of error matrix D}) in neighborhoods of $x$ and $y$. 
For any $B \subseteq \{1, \dots,N\}$ we introduce the matrix 
\bels{def of H^B}{
\bs{H}^B\,=\,(h^B_{xy})_{x,y=1}^N\,,\qquad  
h^B_{xy}\,:=\,h_{xy} \2\bbm{1}(x,y \not \in B )\,,
} 
obtained from $\bs{H}$ by setting the rows and the columns labeled by the elements of $ B $ equal to zero.
The corresponding resolvent is 
\bels{definition GB}{
\bs{G}^B\msp{-1}(\zeta)  \,:=\, (\2\bs{H}^B-\1\zeta\1\bs{1}\2)^{-1}.
}
With this definition, we  have the   resolvent expansion formula
$
\bs{G}= \bs{G}^B - \bs{G}^{B}(\bs{H}-\bs{H}^B)\bs{G}
$.
In particular, for any $y \in B $ the rows of $\bs{G}$ outside $B$ have the expansion
\bels{Expansion of G in B}{
G_{u y}\,=\, - \sum_{v}^B \sum_{z\in B} G_{u v}^B \2 h_{vz} \2 G_{z y}\,,\qquad  u \not \in B\,.
}
Here we introduced, for any two index sets $A,B \subseteq \bb{X}$,  the short hand notation 
\[
\sum_{x \in A}^{B} \;:=\, \sum_{x \in A\setminus B}
\,.
\]
In case $A=\bb{X}$ we simply write $\sum_x^B$ and $ \sum_x = \sum_x^\emptyset$, i.e., the superscript over the summation means exclusion of these indices from the sum. 
Recall that $\bs{H}$ is written as a sum of its expectation matrix $\bs{A}$ and its fluctuation $\frac{1}{\sqrt{N}}\bs{W}$ (cf. \eqref{H structure}) and therefore
\[
\bs{D}\,=\,-  
N^{-1/2}\bs{W}\bs{G}-\cal{S}[\bs{G}]\bs{G}\,.
\]

We use the expansion formula \eqref{Expansion of G in B} on the resolvent elements in $(\bs{W}\bs{G})_{xy}=\sum_u w_{xu}G_{uy}$ and find that the entries of $\bs{D}$ can be written in the form
\bels{First expansion formula for dxy}{
d_{xy}\,&=\, - \frac{1}{\sqrt{N}}\sum_{u \in B}w_{x u}\2G_{u y}\,+\,\frac{1}{\!\sqrt{N\2}}\sum_{u,v}^{B} \sum_{z \in B}w_{x u}\2G_{u v}^{B} \2h_{vz}\2 G_{z y}\msp{0}-\, \frac{1}{N}\sum_{u, v, z}G_{u v}\2(\1\E \2w_{x u}w_{v z})\2G_{z y}\,.
}
Note that the set $ B $ with $y \in B$ here is arbitrary, e.g.,  it may depend on $x$ and $y$. In fact, we will  choose   it to be a neighborhood of $\{x,y\}$, momentarily. 

Let  $A \subseteq B $  be another index set. We split the sum over $z \in B$ in the second term 
 on the right hand side of \eqref{First expansion formula for dxy} into a sum over $w \in A$ and $w \in B \setminus A$ and use \eqref{H structure} again,
\[
\sum_{z \in B}h_{vz}\2 G_{z y}\,=\, \sum_{z \in A}a_{vz}\2 G_{z y}+\sum_{z \in B }^Ah_{vz}\2 G_{z y}+\frac{1}{\sqrt{N}}\sum_{z \in A}w_{vz}\2 G_{z y}
\,.
\]

We end up with the following decomposition of the error matrix  $ \bs{D}= \sum_{\1k\1=\11}^{\25}\bs{D}^{(k)} $, where the entries $d^{(k)}_{xy}$ of the individual matrices $\bs{D}^{(k)}$ are given by
\begin{subequations}
\label{definition of Dk}
\begin{align}
d_{xy}^{(1)}\,&=\, - \frac{1}{\sqrt{N}}\msp{-2}\sum_{u \in B_{xy}}\msp{-5}w_{x u}\2G_{u y}\,,
\\\label{definition of D2}
d_{xy}^{(2)}\,&=\, \frac{1}{\sqrt{N}}\sum_{u,v}^{B_{xy}} \sum_{z \in A_{xy}}\msp{-4}w_{x u}\2G_{u v}^{B_{xy}} a_{vz}\2 G_{z y}\,,
\\
d_{xy}^{(3)}\,&=\, \frac{1}{\sqrt{N}}\sum_{u,v}^{B_{xy}} \sum_{z \in B_{xy}}^{A_{xy}}\msp{-3}w_{x u}\2G_{u v}^{B_{xy}} h_{vz}\2 G_{z y}\,,
\\
\label{definition of D4}
d_{xy}^{(4)}\,&=\, \frac{1}{{N}}\sum_{u,v}^{B_{xy}} \sum_{z \in A_{xy}}\msp{-3}G_{u v}^{B_{xy}} (w_{x u}w_{vz}-\E \2w_{x u}w_{vz})\2 G_{z y}\,,
\\
d_{xy}^{(5)}\,&=\, \frac{1}{{N}}\sum_{u,v}^{B_{xy}} \sum_{z \in A_{xy}}\msp{-3}G_{u v}^{B_{xy}} (\E \2w_{x u}w_{vz})\2 G_{z y}
- \frac{1}{N}\sum_{u, v, z}G_{u v}\2(\E \2w_{x u}w_{v z})\2G_{z y}\,,
\end{align}
\end{subequations}
and
\vspace{-0.3cm}
\bels{definition of Axy and Bxy}{
B_{xy}\,:=\, B_{2N^{\eps_1}}(x) \cup B_{2N^{\eps_1}}(y) \,,\qquad  A_{xy}\,:=\, B_{N^{\eps_1}}(x) \cup B_{N^{\eps_1}}(y) \,,
}
for some $\eps_1>0$. Note that although $\bs{D}$ itself does not depend on the choice of $\eps_1$, its decomposition into $\bs{D}^{(k)}$ does. We will estimate each error matrix $\bs{D}^{(k)}$ separately, where the estimates may still depend on $\eps_1$.   
Since $\eps_1>0$ is arbitrarily small, it is eliminated from the final bounds on $\bs{D}$ using the following property of the stochastic domination (Definition~\ref{def:Stochastic domination}): If some positive random variables $ X,Y$ satisfy $ X \prec N^\eps Y $ for every $ \eps > 0$, then $ X \prec Y $. 

 The following lemma provides entrywise estimates on the individual error matrices. 
\begin{lemma}
\label{lmm:Estimates on error matrices} Let $C>0$ a constant and $\zeta \in \Cp$  with $\dist(\zeta, \spec(\bs{H}^B))^{-1}\prec N^{C}$ for all $B \subsetneq\bb{X}$. 
The entries of the error matrices $\bs{D}^{(k)}=\bs{D}^{(k)}(\1\zeta\1)$, defined in \eqref{definition of Dk}, satisfy the bounds
\begin{subequations}
\label{estimates for Dk}
\begin{align} 
\label{estimate for D1}
\abs{\1d_{xy}^{(1)}}\,
&\prec\,
\frac{\abs{B_{xy}}}{\!\sqrt{N}}\,\norm{\bs{G}}_{\rm{max}} \,,
\\ 
\label{estimate for D2}
\abs{\1d_{xy}^{(2)}}\,
&\prec\,
N\abs{A_{xy}}\2\pbb{\max_{u\not \in B_{xy}}\frac{ \im G_{uu}^{B_{xy}}}{N \im \zeta}}^{\!1/2}
\pbb{\max_{z \in A_{xy}}\sum_{v}^{B_{xy}}\abs{a_{vz}}^2}^{\!1/2}\norm{\bs{G}}_{\rm{max}}\,,
\\ 
\label{estimate for D3}
\abs{\1d_{xy}^{(3)}}\,
&\prec\, 
\frac{\abs{B_{xy}}}{\sqrt{N \im \zeta}} \max_{z \in B_{xy}}\frac{(\2\im G_{zz}^{B_{xy}\setminus\{z\}})^{1/2}
}{
\abs{G_{zz}^{B_{xy}\setminus\{z\}}}}\,\norm{\bs{G}}_{\rm{max}}\,,
\\  
\label{estimate for D4}
\abs{\1d^{(4)}_{xy}}\,
&\prec\, 
\abs{A_{xy}}\,
\biggl(\frac{N^{-1}\!\sum_{u}^{B_{xy}} \im G_{uu}^{B_{xy}}}{N \im \zeta}\biggr)^{\msp{-7}1/2}
\norm{\bs{G}}_{\rm{max}}
\,,
\\
\label{estimate for D5}
\abs{\1d^{(5)}_{xy}}\,
&\prec\, 
\abs{A_{xy}}\abs{B_{xy}}
\max_{k=0}^{\abs{B_{xy}}-1}\biggl(\,\frac{\im G_{x_kx_k}^{B_k}}{\abs{G_{x_kx_k}^{B_k}}{N \im \zeta}}\,+\,\pbb{\frac{\im G_{x_kx_k}^{B_k}}{{N^2}\im \zeta}}^{\msp{-5}1/2}\,\biggr)\,\norm{\bs{G}}_{\rm{max}} 
\\ 
\notag 
&\msp{25}
+\,\pbb{\frac{\im G_{yy}}{N \im \zeta}}^{\msp{-5}1/2}\norm{\bs{G}}_{\rm{max}}\,,
\end{align}
\end{subequations}
where  the  $(B_k)_{k=0}^{\abs{B_{xy}}}$ in \eqref{estimate for D5}  are an arbitrary increasing sequence of 
 subsets of $B_{xy}$ with $B_{k+1}=B_{k}\cup\{x_k\}$ for some $x_k \in B_{xy}$. In particular,
$\emptyset = B_0 \subsetneq B_1 \subsetneq \dots \subsetneq B_{\abs{B_{xy}}-1} \subsetneq  B_{\abs{B_{xy}}}=B_{xy}$.
\end{lemma}

\begin{Proof} We show the estimates \eqref{estimate for D1} to \eqref{estimate for D5} one by one. 
 The bound \eqref{estimate for D1} is trivial since by the bounded moment assumption \eqref{bounded moments} the entries of $\bs{W}$ satisfy $\abs{w_{xy}}\prec 1$. 
 For the proof of \eqref{estimate for D2} we simply use first Cauchy-Schwarz in the $v$-summation of \eqref{definition of D2} and then the \emph{Ward-identity},
\bels{Ward-identity}{
\sum_{u}^B\,\abs{\1G^B_{xu}(\1\zeta\1)}^2\,=\, \frac{\im G^B_{xx}(\1\zeta\1)}{\im \zeta}\,,\quad B \subsetneq \bb{X}\,,\;x \not \in B\,.
}

For \eqref{estimate for D3} we rewrite the entries of $\bs{D}^{(3)}$ in the form 
\bels{rewriting of D3}{
d_{xy}^{(3)}\,=\, -\frac{1}{\sqrt{N}} \msp{-3}\sum_{z \in B_{xy}}^{A_{xy}}\msp{-3}\sum_{u}^{B_{xy}}w_{x u}\2\frac{G_{u z}^{B_{xy}\setminus\{z\}} }{G_{zz}^{B_{xy}\setminus\{z\}} }\, G_{z y}\,,
}
where we used the Schur complement formula in the form of the general resolvent expansion identity 
\[
G_{uz}^B\,=\, -G_{zz}^B\sum^B_v G^{B \cup\{z\}}_{u v} h_{v z}\,,\qquad B \subsetneq \bb{X}\,,\;u,z \not \in B\,.
\]
To the $u$-summation in \eqref{rewriting of D3} we apply the  large deviation estimate \eqref{Single sum large deviation} of Lemma~\ref{lmm:Linear large deviation} with the choices $X_u:=w_{xy}$ and $b_u:=G_{uz}^{B_{xy}\setminus \{z\}}$, i.e.
\bels{LDA for D3}{
\absB{\sum_{u}^{B_{xy}}w_{x u}G_{u z}^{B_{xy}\setminus\{z\}}}
\,\prec\, 
\pbb{\sum_{u}^{\,B_{xy}}\absb{G_{u z}^{B_{xy}\setminus\{z\}}}^2}^{\!1/2}.
}
The assumption \eqref{correlation decay in large deviation} of Lemma~\ref{lmm:Linear large deviation} is an immediate consequence of the decay of correlation \eqref{decay of correlation}. In order to verify \eqref{small correlation between X and b} we use both \eqref{decay of correlation} and the $N$-dependent smoothness 
\bels{N dependent smoothness}{
\norm{\nabla_{\bs{R}}\bs{G}^B}\,=\, 
N^{-1/2}\norm{\bs{G}^B\bs{R}\bs{G}^B}\,\le\, N^{2C}\norm{\bs{R}}\,,
}
of the resolvent, where $\nabla_{\bs{R}}$ denotes the directional derivative with respect to $\bs{W}^B$ in the direction $\bs{R}=\bs{R}^* \in \C^{(N-\abs{B}) \times (N-\abs{B})}$. For the inequality in \eqref{N dependent smoothness} we used the assumption $\dist(\zeta, \spec(\bs{H}^B))\ge N^{-C}$ with high probability.
By the Ward-identity \eqref{Ward-identity} the bound \eqref{estimate for D3} follows from \eqref{LDA for D3} and \eqref{rewriting of D3}.

To show \eqref{estimate for D4} we employ the  quadratic large deviation result Lemma~\ref{lmm:Quadratic large deviation}  with the choices 
\[
X\,:=\,(w_{xu})_{u  \in \bb{X} \setminus B_{xy}}\,,\qquad Y\,:=\,(w_{vy})_{v  \in \bb{X} \setminus B_{xy}}\,,\qquad b_{uv}\,:=\,G_{u v}^{B_{xy}}\,.
\]
The assumptions \eqref{correlation decay in quadratic large deviation} and \eqref{b uncorrelated to X ynd Y} are again easily verified using \eqref{decay of correlation} and \eqref{N dependent smoothness}. Applying  \eqref{Double sum large deviation}  on the $(u,v)$-summation in \eqref{definition of D4} we find 
\[
\absbb{\sum_{u,v}^{B_{xy}} G_{u v}^{B_{xy}} (w_{x u}w_{vz}-\E \2w_{x u}w_{vz})\,}
\;\prec\; 
\pbb{\,\sum_{u,v}^{B_{xy}} \abs{G_{u v}^{B_{xy}}}^2}^{\!1/2}
\!=\, 
\pbb{\,\sum_{u}^{B_{xy}} \frac{\im{G_{u u}^{B_{xy}}}}{\im \zeta}}^{\!1/2}
,
\]
where we used \eqref{Ward-identity} again.

Finally, we turn to the proof of \eqref{estimate for D5}. Let $B_k$ be as in the statement of Lemma~\ref{lmm:Estimates on error matrices}. 
We set
\[
\alpha^{(k)}_{xz}\,:=\,  \frac{1}{{N}}\sum_{u,v}^{B_k} G_{u v}^{B_k}\,\E \2w_{x u}w_{vz} 
\,,
\]
and use a telescopic sum to write $d^{(5)}_{xy}$ as 
\bels{rewriting D5}{
{d^{(5)}_{xy}}\,= 
\sum_{z \in A_{xy}}\msp{-3}\sum_{k=0}^{\abs{B_{xy}}-1}\msp{-3}(\1\alpha^{(k+1)}_{xz}\!-\1\alpha^{(k)}_{xz}\1)\2
G_{zy}
-
\frac{1}{N}\sum_{u, v}\sum_z^{A_{xy}}G_{u v}\2(\1\E \2w_{x u}w_{v z})\2G_{z y}
\,.
}
We estimate the rightmost term in \eqref{rewriting D5} simply by 
\begin{equation*}
\begin{split}
\absbb{\frac{1}{N}\sum_{u, v}\sum_z^{A_{xy}}G_{u v}\2(\1\E \2w_{x u}w_{v z})\2G_{z y}}^2
\,&\le\, 
\frac{\norm{\bs{G}}_{\rm{max}}^2}{N^2}\,
\sum_z^{A_{xy}}
\pbb{\,\sum_{u, v}\,\abs{\1\E \2w_{x u}w_{v z}}}^{\!2}\,
\sum_z^{A_{xy}}\1\abs{\1G_{z y}}^2
\lesssim\,
\norm{\bs{G}}_{\rm{max}}^2\frac{\im G_{yy}}{N\im \zeta}
\,,
\end{split}
\end{equation*}
where the sum over $u$ and $v$ on the right hand side of the first inequality is bounded by a constant because of the  decay of covariances  \eqref{decay of covariances}  and we used \eqref{Ward-identity} in the second ineqality. Thus, \eqref{estimate for D5} follows from \eqref{rewriting D5} and the bound 
\bels{bound on alpha difference}{
\abs{\1\alpha^{(k+1)}_{xz}-\alpha^{(k)}_{xz}}
\;\prec\;
\frac{1}{N \im \zeta\1}\2\frac{\im G_{x_k x_k}^{B_k}}{\abs{G_{x_k x_k}^{B_k}}} 
+\frac{1}{\sqrt{N}}\pbb{\frac{\im G_{x_kx_k}^{B_k}}{{N}\im \zeta}\!}^{\msp{-5}1/2}
.
}
To show \eqref{bound on alpha difference} we first see that
\bels{difference of alphas}{
\msp{-15}
\alpha^{(k+1)}_{xz}\!-\alpha^{(k)}_{xz}
\,=
-\frac{1}{{N}}\msp{-6}\sum_{u,v}^{B_{k+1}} \msp{-5}\frac{G_{\!u x_k}^{B_k}G_{\!x_k \msp{-2}v}^{B_k}\!}{G^{B_k}_{x_kx_k}} \, \E \2w_{x u}w_{vz}
\!-\frac{1}{{N}}\!\sum_{u}^{B_{k}}G_{\!u x_k}^{B_{k}}\2 \E \2w_{x u}w_{x_kz}
- 
\frac{1}{{N}}\msp{-10}\sum_{v}^{\;B_{k+1}}
\msp{-8}G_{\!x_k v}^{B_{k}} \E \2w_{x x_k}w_{vz}
\,,
}
where we used the general resolvent identity
\[
G^{B}_{xy}\,=\,G^{B\cup\{u\}}_{xy}+ \frac{G^{B}_{xu}G^{B}_{uy}}{G^{B}_{uu}}\,,\qquad B \subsetneq\{1, \dots,N\}\,,\; x,y,u \not \in B\,,\; x \ne u\,,\; y \ne u\,.
\]
The last two terms on the right hand side of \eqref{difference of alphas} are estimated by the second term on the right hand side of \eqref{bound on alpha difference} using first Cauchy-Schwarz, the decay of covariances  \eqref{decay of covariances},   and then the Ward-identity \eqref{Ward-identity}. For the first term in \eqref{difference of alphas} we use the same argument as in \eqref{Bound S by Qnu} to see that \eqref{decay of covariances} implies 
\bels{S applied to rank one matrix}{
\sum_{u,\1v} \,\abs{\1r_u t_v}\, \abs{\1\E \2w_{x u}w_{vz}}\,\lesssim\, \norm{\bs{r}}\2\norm{\1\bs{t}}
\,,
}
for any two vectors $\bs{r},\bs{t} \in \C^{N \times N}$. We obtain \eqref{bound on alpha difference} by applying \eqref{S applied to rank one matrix} with the choice $r_u:=G_{u x_k}^{B_k}$, $t_v:=G_{x_k v}^{B_k}$ and using the Ward-identity afterwards. 
   In this way \eqref{estimate for D5} follows and Lemma~\ref{lmm:Estimates on error matrices} is proven.
\end{Proof}

The following definition is motivated by the formula that expresses  the matrix   elements of $\bs{G}^{B}$   in terms of
 the matrix   elements of $\bs{G}$. For $\bs{R} \in \C^{N \times N}$ and $A,B \subsetneq \bb{X}$ we denote by $\bs{R}_{AB}:=(r_{xy})_{x \in A, y \in B}$ its submatrix. In case $A=B$ we write $\bs{R}_{AB}=\bs{R}_{A}$ for short.  
 Then we have
\bels{GB identity}{
\bs{G}^B_{\bb{X}\setminus B}\,=\, (\2\bs{H}_{\bb{X}\setminus B}-\zeta\2\bs{1})^{-1}\,=\, ((\bs{G}^{-1})_{\bb{X}\setminus B})^{-1}
.
}
In particular, $ (\bs{G}^B)_{\bb{X}\backslash B} = \bs{G}^B_{\bb{X}\backslash B}$. 

\begin{definition}
For $B \subsetneq \bb{X}$ we define the $\C^{(N-\abs{B})\times (N-\abs{B})}$-matrix
\bels{definition of mBxx}{
\bs{M}^B\,:=\,((\bs{M}^{-1})_{\bb{X}\setminus B})^{-1}.
}
\end{definition}

\begin{lemma} 
\label{lmm:Removed resolvent bounds}
Let $\delta>0$ and $\zeta\in \Cp$  be such that
$
\delta \le \dist(\zeta,[\kappa_-,\kappa_+])+\rho(\1\zeta\1)\le \delta^{-1}
$.
Then for all $  B \subsetneq \bb{X}$ the matrix $\bs{M}^B$, defined in \eqref{definition of mBxx},  satisfies
\bels{decay of MB}{
\norm{\1\bs{M}^B}_{\ul{\gamma}}\,\lesssim_{\1\delta} 1\,,
}
for some sequence $\ul{\gamma}$, depending only on $\delta$ and the model parameters. 
For every $x \not \in B$ we have
\bels{Bounds on mB}{
\abs{m^B_{xx}(\1\zeta\1)}\,\sim_\delta 1\,,\qquad  \im m^B_{xx}(\1\zeta\1)\,\sim_\delta \rho(\1\zeta\1)
\,.
}
Furthermore, there is a positive constant $c$, depending only on $\scr{K}$ and $\delta$, such that 
\bels{bound on GB-mB}{
\max_{x,\1y\1\notin\1 B}\2\abs{\1G^B_{xy}(\1\zeta\1)-m^B_{xy}(\1\zeta\1)}\,\bbm{1}\msp{-1}\pB{\1\Lambda(\1\zeta\1) \le \tsfrac{c}{\21\2+\abs{B}}\,}
\;\lesssim_{\2\delta}\;
(\11+ \abs{B}\1)\,\Lambda(\1\zeta\1)
\,.
}
\end{lemma}
\begin{Proof} We begin by establishing upper and lower bounds on the singular values of $\bs{M}^B$,
\bels{singular value bounds for MB}{
\norm{\1\bs{M}^B} \,\sim_\delta 1\,,\qquad \norm{(\bs{M}^B)^{-1}}\,\sim_\delta 1\,.
}

We will make use of the following general fact: If $\bs{R}\in \C^{N \times N}$ satisfies $\norm{\bs{R}}\lesssim_\delta 1$, as well as 
\bels{lower bound on im or re}{
\im \bs{R}\,\gtrsim_\delta \bs{1}\,,\qquad \text{or}\qquad \re \bs{R}\,\gtrsim_\delta\bs{1}\,,\qquad \text{or}\qquad -\re \bs{R}\gtrsim_\delta\bs{1}¨
\,,
}
then any submatrix $\bs{R}_{A}$ of $\bs{R}$ satisfies
\bels{Conclusion for submatrix of R}{
\norm{\1\bs{R}_A}\,\sim_\delta1\,,\qquad\norm{(\bs{R}_A)^{-1}}\,\sim_\delta\, 1\,, \qquad A \subseteq \bb{X}\,.
}

We verify \eqref{lower bound on im or re} for $\bs{R}=\bs{M}$ in two separate regimes and thus show \eqref{singular value bounds for MB}. First  let $\zeta$ be such that $\rho(\zeta)\ge \delta/2$. Then the lower bound in the imaginary part in \eqref{lower bound on im or re} follows from \eqref{im M lower and upper bound} and \eqref{M upper bound}.

Now let $\zeta$ be such that $\delta/2 \le \dist(\zeta,[\kappa_-,\kappa_+])\le \delta^{-1}$. 
 Then we may also assume that we have $\dist(\re \zeta, [\kappa_-,\kappa_+])\ge \delta/4 $, because otherwise $\im \zeta \ge \delta/4$ and thus $\rho(\zeta)\gtrsim_\delta 1$. In this situation the claim follows from the case that we already considered, namely $\rho(\zeta)\ge \delta/2$ because there $\delta$ was arbitrary.
 Since $\bs{M}$ is the Stieltjes transform of a $\scr{C}_+$-valued measure with support in $[\kappa_-,\kappa_+]$ (cf. \eqref{M as Stieltjes transform}), its real part is positive definite to the left of $\kappa_-$ and negative definite to the right of $\kappa_+$. In both cases we also have the effective bound $\abs{\re \bs{M}}\gtrsim_\delta \bs{1}$ because $\dist(\re \zeta, [\kappa_-,\kappa_+])\ge \frac{\delta}{4} $. 
 
 Now we apply \eqref{singular value bounds for MB} to see \eqref{decay of MB}.
 By \eqref{decay of expectation} and \eqref{S operator decay} the right hand side of \eqref{RM MDE} and with it $\bs{M}^{-1}$  has faster than power law decay. The same is true for its submatrix with indices in $\bb{X} \setminus B$. Thus \eqref{decay of MB} follows directly from the definition \eqref{definition of mBxx} of $\bs{M}^B$, the upper bound on its singular values from \eqref{singular value bounds for MB} and  the Combes-Thomas estimate in  Lemma~\ref{lmm:Perturbed Combes-Thomas}.

To prove \eqref{Bounds on mB} we use 
\[
\im \bs{M}^{B}
\,=\, -\2(\bs{M}^{B})^* \im (\1\bs{M}^{-1})_{\bb{X}\setminus B}\2\bs{M}^{B}
\;\sim_{\1\delta}\, 
- \im (\1\bs{M}^{-1})_{\bb{X}\setminus B}
\;\sim_\delta\, 
\rho\2\bs{1}
\,,
\]
where we applied \eqref{singular value bounds for MB} for the first comparison relation and used $-\im \bs{M}^{-1}\sim_\delta \rho\1\bs{1}$ (cf. \eqref{im M lower and upper bound} and \eqref{M lower bound}) for the second. The bound on $\im m_{xx}^B$ in  \eqref{Bounds on mB} follows and the bound on $\abs{m_{xx}^B}$  follows at least in the regime $\rho(\zeta)\ge \delta/2$. 
We are left with showing $\abs{m_{xx}^B}\gtrsim_\delta1 $ in the case $\delta/2 \le \dist(\zeta,[\kappa_-,\kappa_+])\le \delta^{-1}$. As we did above, we may assume that $\dist(\re \zeta, [\kappa_-,\kappa_+])\ge \frac{\delta}{4} $.  We restrict to $\re \zeta \le \kappa_--\frac{\delta}{2}$.
The case $\re \zeta \ge \kappa_++\frac{\delta}{2}$ is treated analogously. In this regime 
\[
\re \bs{M}^{B}
\,=\, 
(\1\bs{M}^{B})^* \re (\1\bs{M}^{-1})_{\bb{X}\setminus B}\2\bs{M}^{B}
\,\sim_\delta\,
\re (\1\bs{M}^{-1})_{\bb{X}\setminus B}\,\sim_\delta \bs{1}
\,,
\]
where we used  $\re \bs{M}^{-1}= (\1\bs{M}^{-1})^*(\1\re \bs{M}\1)\2\bs{M}^{-1} \!\sim_\delta \re \bs{M}\sim_\delta\bs{1}$  for the last comparison relation. Thus, \eqref{Bounds on mB} follows.

Now we show \eqref{bound on GB-mB}. By the Schur complement formula we have for any $\bs{T} \in \C^{N \times N}$ the identity 
\bels{submatrix Schur}{
\bigl((\bs{T}^B)_{\{x,y\}}\bigr)^{-1}
=\, 
\bigl(\2\bs{T}_{\{x,y\}}- \bs{T}_{\!\{x,y\}B}(\bs{T}_{B})^{-1}\bs{T}_{\!B\{x,y\}}\bigr)^{-1}
\,=\, (\1(\bs{T}_{\!B \cup\{x,y\}})^{-1})_{\msp{-1}\{x,y\}}
\,,
}
for $x,y \not \in B$ and $\bs{T}^B:=((\bs{T}^{-1})_{\bb{X}\setminus B})^{-1}$, provided all inverses exist. 
We will use this identity for $\bs{T}= \bs{M},\bs{G}$. Note that this definition $\bs{T}^B$ with $\bs{T}=\bs{G}$ is consistent with the definition \eqref{definition GB} on the index set $\bb{X}\setminus B$ because of \eqref{GB identity}.
Recalling  that $\bs{G}_{B\cup \{x,y\}}=(G_{u,v})_{u,v \in B\cup \{x,y\}}$ and  $\bs{M}_{B\cup \{x,y\}}$ are
 matrices of dimension $|B|+2$,  we have
\[
\norm{\bs{G}_{B\cup \{x,y\}}-\bs{M}_{B\cup \{x,y\}}}
\,\le\, 
(\abs{B}+2)\,\normb{\bs{G}_{B\cup \{x,y\}}-\bs{M}_{B\cup \{x,y\}}}_{\rm{max}}
\,\le\, 
(\abs{B}+2)\Lambda\,.
\]
Therefore, as long as  $(\abs{B}+2)\2\Lambda\2\norm{(\bs{M}_{B\cup \{x,y\}})^{-1}}\le \frac{1}{2}$  we get 
\bea{
\normb{(\bs{G}_{B\cup \{x,y\}})^{-1}-(\bs{M}_{B\cup \{x,y\}})^{-1}}
\,&\le\, 
2\,\normb{(\bs{M}_{B\cup \{x,y\}})^{-1}}^2\2
\normb{\bs{G}_{B\cup \{x,y\}}-\bs{M}_{B\cup \{x,y\}}}
\,\lesssim_{\1\delta}\, (\21+\abs{B})\2\Lambda
\,,
}
where we used in the last step that $\norm{(\bs{M}_{B\cup \{x,y\}})^{-1}}\sim_\delta 1$, which follows from using \eqref{lower bound on im or re} and \eqref{Conclusion for submatrix of R} for the choice $\bs{R}=\bs{M}$
 in the regimes $\rho\gtrsim_\delta 1$ and $\dist(\re \zeta, [\kappa_-,\kappa_+])\gtrsim_\delta 1$, respectively.   

Again using the definite signs of the imaginary and real part of $\bs{M}$ as well as that of $(\bs{M}_{B \cup\{x,y\}})^{-1}$ in these two regimes, we infer that 
\[
\normb{\bigl(((\bs{M}_{B \cup\{x,y\}})^{-1})_{\{x,y\}}\bigr)^{-1}} 
\,\sim_\delta\, 
1
\,,
\]
as well. We conclude that there is a constant $c$, depending only on $\delta$ and $\scr{K}$, such that 
\[
\normB{
\bigl(((\bs{G}_{B \cup\{x,y\}})^{-1})_{\{x,y\}}\bigr)^{-1}
\!-\2
\bigl(((\bs{M}_{B \cup\{x,y\}})^{-1})_{\{x,y\}}\bigr)^{-1}}\,
\bbm{1}\pB{\Lambda \le \frac{c}{1+\abs{B}}}
\,\lesssim_\delta\, 
(1+\abs{B})\2\Lambda
\,. 
\]
With the identity \eqref{submatrix Schur} the claim \eqref{bound on GB-mB} follows and Lemma~\ref{lmm:Removed resolvent bounds} is proven.
\end{Proof}

\begin{Proof}[Proof of Lemmas~\ref{lmm:Smallness of error matrix} and \ref{lmm:Smallness of error matrix away from support}]
We begin with the proof of \eqref{bound on error D away from real line}. We continue the estimates on all the error matrices listed in Lemma~\ref{lmm:Estimates on error matrices}. Therefore, we fix $\zeta =\tau +\ii\eta$ with $\abs{\tau}\le C$ and $\eta \in [1,C]$. 
Since $\im \zeta \ge1$, we have the trivial resolvent bound  and a lower bound on diagonal elements,
\bels{trivial resolvent bound AND trivial bound on GBxx inverse}{
\norm{\bs{G}^B}\,\le\, 1\,,
\qquad\text{and}\qquad
\frac{1}{\abs{G_{xx}^B}}\,\prec\, 1\,,\qquad x \not \in B  \subseteq\bb{X}
\,.
}
Indeed, to get the lower bound we apply the Schur complement formula applied to the $(x,x)$-element of the resolvent $\bs{G}^B=(\bs{H}^B-\zeta\2\bs{1})^{-1}$ to obtain
\[
-\frac{1}{G_{xx}^B}\,=\, \zeta-a_{xx}+\sum_{u,v}^Bh_{x u}G_{u v}^{B \cup \{x\}}h_{vx}\,.
\]
We take absolute value on both sides and estimate trivially,
\[
\frac{1}{\abs{G_{xx}^B}}\,\le\, \abs{\zeta\1}+\abs{\1a_{xx}}+\norm{\bs{G}^{B \cup \{x\}}}\,\sum_u\,\abs{h_{xu}}^2
\,\prec\; 
1\,.
\]
Here we used the first bound of \eqref{trivial resolvent bound AND trivial bound on GBxx inverse} to control the norm of the resolvent and the assumptions \eqref{decay of expectation} and \eqref{bounded moments} to bound $\sum_u\abs{h_{xu}}^2$. 
Combining \eqref{definition of Axy and Bxy} and and the assumption \eqref{polynomial ball growth} we get
\bels{bound on size of Axy and Bxy}{
\abs{A_{xy}}\,\le\, \abs{B_{xy}}\,\prec\, N^{\eps_1\msp{-2}P}.
}

Using \eqref{bound on size of Axy and Bxy} and \eqref{trivial resolvent bound AND trivial bound on GBxx inverse} in the main estimates \eqref{estimates for Dk} for $ \abs{\1d^{(k)}_{xy}} $'s yields
\bels{dxy bound for large im zeta}{
\abs{\1d_{xy}}\,\prec\, \frac{N^{2\2\eps_1\msp{-2}P}}{\sqrt{N}}+N^{\eps_1\msp{-2}P}N \frac{\kappa_2(\nu)}{N^{\eps_1\nu}}
\,,\qquad\text{for all }\nu \in \N
\,.
}
Here we also used that by assumption \eqref{decay of expectation} for any $\nu \in \N$ the expectation matrix satisfies
\bels{N bound on expectation far from diagonal}{
\abs{\1a_{vz}}
\,\le\, 
\kappa_2(\nu) N^{-\eps_1 \nu}
\,,\qquad z \in A_{xy}\,,\; v \not \in B_{xy}\,,
}
to obtain the second summand on the right hand side of \eqref{dxy bound for large im zeta} from estimating $\abs{d^{(2)}_{xy}}$.
Since $\eps_1>0$ was arbitrary \eqref{dxy bound for large im zeta} implies \eqref{bound on error D away from real line}.

Now we prove \eqref{bound on error D in the bulk} and \eqref{error bound away from convex hull of support} in tandem. Let $\delta>0$ and $\zeta \in \Cp$ such that  $\delta\le \dist(\zeta,[\kappa_-,\kappa_+])+\rho(\zeta)\le \delta^{-1}$ and $\im \zeta \ge N^{-1+\eps}$. We show that
\bels{general entrywise bound on D}{
\norm{\bs{D}}_{\rm{max}}\,\bbm{1}(\2\Lambda\le  N^{-\eps})
\,\prec\,
\pB{\frac{\rho+\Lambda}{N \im \zeta}}^{\!1/2}.
}
From \eqref{general entrywise bound on D} the bound \eqref{bound on error D in the bulk} follows immediately in the regime where $\rho \ge \delta$. Also \eqref{error bound away from convex hull of support} follows from \eqref{general entrywise bound on D}. Indeed, in the regime of spectral parameters $\zeta \in \Cp$ with $\delta\le \dist(\zeta,[\kappa_-,\kappa_+])\le \delta^{-1}$ we have $\rho\sim_\delta \im \zeta$ because $\rho$ is the harmonic extension of a probability density supported inside $[\kappa_-,\kappa_+]$. 

For the proof of \eqref{general entrywise bound on D} we use \eqref{bound on GB-mB}, \eqref{Bounds on mB},
 \eqref{N bound on expectation far from diagonal}, \eqref{bound on size of Axy and Bxy} and $\norm{\bs{G}}_{\rm{max}}\lesssim 1+\Lambda$ (cf. \eqref{M upper bound}) to estimate the right hand side of each inequality in \eqref{estimates for Dk}.
In this way we get
\bels{entrywise estimate of D with eps1}{
\norm{\1\bs{D}\1}_{\rm{max}}\,\bbm{1}(\Lambda\le  N^{-\eps})
\,\prec\, 
N^{2\1\eps_1\msp{-2}P}N^{3/2}\pB{\frac{ \rho+\Lambda}{N \im \zeta}}^{\!1/2} 
\frac{\kappa_2(\nu)}{N^{\eps_1\nu}} 
+ 
N^{2\1\eps_1\msp{-2}P}
\pB{\frac{\rho+\Lambda}{N \im \zeta}}^{\!1/2},
}
for any $\nu \in \N$, 
provided $\eps > \eps_1 P$  to ensure $N^{-\eps} \le c/|B_{xy}|$, i.e. that the constraint $\Lambda\le N^{-\eps}$ makes 
\eqref{bound on GB-mB} applicable. 
Here, we also used $\rho \gtrsim_\delta\im \zeta$ to see that $\frac{\rho}{N \im \zeta}\gtrsim_\delta \frac{1}{N}$. Since \eqref{entrywise estimate of D with eps1} holds for arbitrarily small $\eps_1>0$, the claim \eqref{general entrywise bound on D} and with it Lemmas~\ref{lmm:Smallness of error matrix} and \ref{lmm:Smallness of error matrix away from support} are proven. 
\end{Proof}

\section{Fluctuation averaging}
\label{sec:Fluctuation averaging}
In this section we prove Proposition~\ref{prp:fluctuation averaging} by which a error bound $\Psi$ for the entrywise local law can be used to improve the bound on the error matrix $\bs{D}$ to $\Psi^2$, once $\bs{D}$ is averaged against a non-random matrix $\bs{R}$ with faster than power law decay. 

\begin{Proof}[Proof of Proposition~\ref{prp:fluctuation averaging}] Let $\bs{R} \in \C^{N \times N}$ with $\norm{\bs{R}}_\beta\le 1$ for some positive sequence $\beta$. 
Within this proof we use Convention~\ref{con:Comparison relation} such that $\varphi\lesssim \psi$ means  $\varphi\le C \psi$ for a constant $C$, depending only on $\wt{\scr{P}}:=(\scr{K}, \delta, \eps_1,\beta,C)$, where $C$ and $\delta$ are the constants from the statement of the proposition, $\scr{K}$ are the model parameters (cf \eqref{definition of model parameters K}) and $\eps_1$ enters in the splitting of the error matrix $\bs{D}$ into $\bs{D}^{(k)}$ (cf. \eqref{definition of Axy and Bxy}). Note that since $\eps_1$ is arbitrary it suffices to show \eqref{general FA} up to factors of $N^{\eps_1}$. We will also use the notation $\ord_\prec(\Phi)$ for a random variable that is stochastically dominated by some nonnegative $\Phi$. 

We split the expression $\avg{\bs{R},\bs{D}}$ from \eqref{general FA} according to the definition \eqref{definition of Dk} of the matrices $\bs{D}^{(k)}$. Then we estimate  $\avg{\bs{R},\bs{D}^{(k)}}$ separately for every $k$. We do this in three steps. 
First we estimate $\norm{\bs{D}^{(k)}}_{\rm{max}}$ for $k=2,3,5$ directly without using the averaging effect of $\bs{R}$. Afterwards we show the bounds on $\avg{\bs{R},\bs{D}^{(1)}}$ and $\avg{\bs{R},\bs{D}^{(4)}}$, respectively. In the upcoming arguments the following observation will be useful. 
The local law \eqref{Assumption of FA on G-M} together with \eqref{bound on GB-mB} implies that for every $B \subseteq \bb{X}$ with $\abs{B}\le N^{\eps/2}$ we have 
\bels{local law for GB}{
\norm{\1\bs{M}^B\!-\bs{G}^B}_{\rm{max}}\,\prec\, (1+\abs{B})\Psi
\,.
}
Here, until the end of this proof, we consider $\bs{G}^B$ as the $\C^{(N-\abs{B}) \times (N-\abs{B})}$-matrix $\bs{G}^B=(G_{xy}^B)_{x,y \not \in B}$ as opposed to the general convention \eqref{def of H^B}. 
\medskip
\\
\emph{Estimating }$\norm{\bs{D}^{(k)}}_{\rm{max}}$: Here, we show that under the assumption \eqref{Assumption of FA on G-M}  the error matrices with indices $k=2,3,5$ satisfy the improved entrywise bound
\bels{Psi2 bound on Dk}{
\norm{\bs{D}^{(k)}}_{\rm{max}}\,\prec\, N^{3\1\eps_1 P}\Psi^2\,,\qquad k=2,3,5\,,
}
where $\eps_1$ stems from \eqref{definition of Axy and Bxy} and $P$ is the model parameter from \eqref{polynomial ball growth}.

We start by estimating the entries of $\bs{D}^{(2)}$. Directly from its definition in \eqref{definition of D2} we infer
\[
\abs{\2d_{xy}^{(2)}}
\;\prec\; 
N^{3/2} \abs{A_{xy}}\, \norm{\bs{G}}_{\rm{max}} \norm{\bs{G}^{B_{xy}}}_{\rm{max}}\,
{\textstyle \max_{\1z \1\in\1 A_{xy}}\!\sum_{\,v}^{B_{xy}}} \abs{\1a_{vz}}
\,.
\]
The maximum norm on the entries of the resolvents $\bs{G}$ and $\bs{G}^{B_{xy}}$ are bounded by \eqref{local law for GB} and \eqref{decay of MB}. The decay \eqref{Expectation decay} of the entries of the bare matrix and that $d(v,z)\ge N^{\eps_1}$ in the last sum then imply $\norm{\bs{D}^{(2)}}_{\rm{max}}\prec N^{-\nu}$ for any $\nu \in \N$. 

To show \eqref{Psi2 bound on Dk} for $k=3$ we use the representation \eqref{rewriting of D3} and the large deviation estimate \eqref{LDA for D3} just as we did in the proof of Lemma~\ref{lmm:Estimates on error matrices}. In this way we get 
\[
\abs{d_{xy}^{(3)}}
\,\prec\,
\frac{\abs{B_{xy}}}{\sqrt{N}}\,\biggl(\,\max_{z \not \in A_{xy}}\sum_{u}^{B_{xy}}\abs{G_{u z}^{B_{xy}\setminus\{z\}}}^2\biggr)^{\!1/2}\msp{-7}\frac{\abs{G_{zy}}}{\absb{G^{B_{xy}\setminus\{z\}}_{zz}}}\,.
\]
Now we use  \eqref{local law for GB}, \eqref{decay of MB}, \eqref{Bounds on mB} and \eqref{bound on size of Axy and Bxy} to conclude 
\bels{improved bound on D3}{
\norm{\bs{D}^{(3)}}_{\rm{max}}
\,\prec\, 
\frac{N^{\eps_1 P}\!}{\sqrt{N}}\biggl(\Psi+\max_{ z \not \in A_{xy}}\abs{m_{zy}}\biggr)\,.
}
The faster than power law decay of $\bs{M}$ from \eqref{Arbitrarily high polynomial decay of solution} together with the definition of $A_{xy}$ in \eqref{definition of Axy and Bxy} implies 
$
\max_{z \not \in A_{xy}}\abs{m_{zy}} \leq C(\nu) N^{-\eps_1\nu}
$ 
for any $\nu\in \N$. Since $\Psi\ge N^{-1/2}$ we infer \eqref{Psi2 bound on Dk} for $k=3$ from \eqref{improved bound on D3}. 

Finally we consider the case $k=5$. We follow the proof of Lemma~\ref{lmm:Estimates on error matrices} and use the representation \eqref{rewriting D5}. We estimate the two summands on the right hand  side of \eqref{rewriting D5}, starting with the second term. We rewrite this term in the form $\sum_z^{A_{xy}}\cal{S}_{xz}[\bs{G}]\2G_{zy}$ and use \eqref{S operator decay} as well as $\norm{\bs{G}}_{\rm{max}}\le \norm{\bs{M}}_{\rm{max}}+\Psi$ together with the upper bound on $\bs{M}$ in \eqref{Arbitrarily high polynomial decay of solution}. 

To bound the first term on the right hand  side of \eqref{rewriting D5} we use \eqref{difference of alphas}. Each of the three terms on the right hand side of \eqref{difference of alphas} has to be bounded by $N^{2\1\eps_1 P}\Psi^2$. The second and third term are bounded by $\frac{1}{N}$ by the  decay of covariances \eqref{decay of covariances}.  For the first term we use \eqref{S applied to rank one matrix}, \eqref{local law for GB} and \eqref{decay of MB}.
\medskip
\\
\emph{Estimating }$\avg{\2\bs{R}\1,\bs{D}^{(1)}}$: 
Here we will show that
\bels{Estimate on RD1}{
\abs{\avg{\2\bs{R}\1,\bs{D}^{(1)}}}\,\prec\,N^{CP\eps_1} \Psi^2\,,
}
for some numerical constant $C>0$.

We split the error matrix $\bs{D}^{(1)}$ into two pieces $ \bs{D}^{(1)}= \bs{D}^{\rm{(1a)}}+\bs{D}^{\rm{(1b)}} $, defined by 
\[
d_{xy}^{\rm{(1a)}}\,:=\, - \frac{1}{\sqrt{N}}\msp{-2}\sum_{u \in B_{xy}}\msp{-5}w_{x u}\2m_{u y}
\,,
\qquad\text{and}\qquad 
d_{xy}^{\rm{(1b)}}\,:=\, - \frac{1}{\sqrt{N}}\msp{-2}\sum_{u \in B_{xy}}\msp{-5}w_{x u}\2(G_{uy}-m_{u y})
\,,
\]
where $ B_{xy} $ is a $ 2N^{\eps_1}$-environment of the set $ \sett{x,y} $ (cf. \eqref{definition of Axy and Bxy}).  
The second part is trivially bounded, $ \norm{\bs{D}^{\rm{(1b)}}}_{\rm{max}} \prec N^{\eps_1 P} \Psi^2 $, using the local law \eqref{local law for GB}, with $ B = \emptyset $. 

For the bound on $\avg{\bs{R},\bs{D}^{\rm{(1a)}}}$ we write
\bels{Splitting D1a into X+Y+Z}{
\avg{\2\bs{R}\1,\bs{D}^{(1a)}}\,=\, X+Y+Z
\,,
} 
where the term on the right hand side are sums of $ \sigma_{xuy}:=N^{-3/2}\,\ol{r}_{xy}w_{x u}\2m_{u y} $ over disjoint  domains
\begin{equation*}
\begin{split}
 X\,:=\,  \sum_{x}\sum_{y}^{B_x^1}\sum_{u \in B_{x}^2}^{B_y^2}\msp{-5}\sigma_{xuy}
\,,
\quad 
Y\,:=\,  \sum_{x}\sum_{y}^{B_x^1}\sum_{u \in B_{y}^2}\msp{-5}\sigma_{xuy}
\,,
\quad
Z\,:=\,  \sum_{x}\sum_{y\in B_x^1}\sum_{u \in B_{xy}}\msp{-5}\sigma_{xuy}\,,
\end{split}
\end{equation*}
expressed in terms of the following metric balls:
\[ 
B_x^k\,:=\, B_{kN^{\eps_1}}(x)\,.
\]
The fast decay of off-diagonal entries for $\bs{R}$ and $\bs{M}$, $\abs{r_{xy}}+\abs{m_{xy}}\lesssim \frac{1}{N}$ for $d(x,y)\ge N^{\eps_1}$ (cf. \eqref{Arbitrarily high polynomial decay of solution}), yields immediately $ \abs{X} \prec_\mu N^{2P\eps_1\mu}N^{-3\mu} $. This suffices for \eqref{Estimate on RD1}. 
The off-diagonal decay also yields
\begin{subequations}
\label{bounds on X Y Z for D1a}
\begin{align}
\label{estimate Y for D1a}
\E\2\abs{Y}^{2\mu}\,
&\lesssim_\mu 
\frac{N^{2P\eps_1\mu}}{N^{5\mu}}\sum_{\bs{x},\bs{u}}
\absB{\,\E\prod_{i=1}^\mu w_{x_iu_i}\ol{w}_{x_{\mu+i}u_{\mu+i}}}\,,
\\
\label{estimate Z for D1a}
\E\2\abs{Z}^{2\mu}\,
&\lesssim_\mu 
\frac{N^{2P\eps_1\mu}}{N^{3\mu}}\sum_{\bs{x}}\sum_{\bs{u} \in B_{\bs{x}}^3}
\absB{\,\E\prod_{i=1}^\mu w_{x_iu_i}\ol{w}_{x_{\mu+i}u_{\mu+i}}}
\,,
\end{align}
\end{subequations}
where the sums are over index tuples $\bs{x}=(x_1,\dots,x_{2\mu}) \in \bb{X}^{2\mu}$ and $B_{\bs{x}}^k:=B_{kN^{\eps_1}}(\bs{x})$ is the ball around $\bs{x}$ with respect to the product metric
\[
d(\bs{x},\bs{y}):=\max_{i=1}^{2\mu} d(x_i,y_i)\,.
\]
In \eqref{estimate Z for D1a} we have used the triangle inequality to conclude that $ d(\bs{u},\bs{x}) \leq 3N^{\eps_1} $.
For $Y$ and $Z$ we continue the estimates in \eqref{bounds on X Y Z for D1a} by using the decay of correlations \eqref{decay of correlation} and the ensuing lumping of index pairs $(x_i,u_i)$: 
\bels{Lumping bound for D1a}{
\absB{\E\prod_{i=1}^\mu w_{x_iu_i}\ol{w}_{x_{\mu+i}u_{\mu+i}}}
\,\lesssim_{\mu,\nu}
\begin{cases}
N^{-\nu}&\exists\,i \text{ s.t. }  d_{\rm{sym}} ((x_i,u_i),\{(x_j,u_j): j\ne i\})\ge N^{\eps_1};
\\
1& \text{ otherwise, } 
\end{cases}
}
where $d_{\rm{sym}}((x_1,x_2),(y_1,y_2)):=d((x_1,x_2), \{(y_1,y_2),(y_2,y_1)\})$ is the symmetrized distance on $\bb{X}^2$, induced by $d$. Inserting \eqref{Lumping bound for D1a} into the moment bounds on $Y$ and $Z$ effectively reduces the combinatorics of the sum in \eqref{estimate Y for D1a} from $N^{4\mu}$ to $N^{2\mu}$ and in \eqref{estimate Z for D1a} from $N^{2\mu}$ to $N^{\mu}$. 
We conclude that  $\abs{\avg{\bs{R},\bs{D}^{(1\rm{a})}}} \prec N^{CP\eps_1}N^{-1} $. 
Moreover, together with $\Psi\ge N^{-1/2}$ and our earlier estimate for $ \avg{\bs{R},\bs{D}^{(1\rm{b})}} $ this yields \eqref{Estimate on RD1}.
\medskip
\\
\emph{Estimating }$\avg{\2\bs{R}\1,\bs{D}^{(4)}}$: Similarly to the strategy for estimating $\abs{\avg{\bs{R},\bs{D}^{(1)}}}$ we write
\[
\bs{D}^{(4)}
\,=\,
\bs{D}^{\rm{(4a)}}+\bs{D}^{\rm{(4b)}} 
\,,\qquad
d_{xy}^{\rm{(4a)}}\,:=\, \sum_{z \in A_{xy}}Z^{B_{xy}}_{xz}\2m_{zy}
\,,\qquad 
d_{xy}^{\rm{(4b)}}\,:=\, \sum_{z \in A_{xy}}Z^{B_{xy}}_{xz}\2(G_{zy}-m_{zy})
\,,
\]
where $ A_{xy} $ is from  \eqref{definition of Axy and Bxy}, and
we have introduced for any $B \subsetneq \bb{X} $ the short hand
\bels{definition of ZBxy}{
Z^B_{xz}\,:=\,\frac{1}{{N}}\sum_{u,v}^{B}G_{u v}^{B}\, (\1w_{x u}w_{vz}-\E \2w_{x u}w_{vz})
\,.
}
From the decay of correlations \eqref{decay of correlation} and $d(\{x,z\},\bb{X}\setminus B_{xy})\ge N^{\eps_1}$ for any $z \in A_{xy}$ as well as the $N$-dependent smoothness of the resolvent as a function of the matrix entries of $\bs{W}$ for $\dist(\zeta, \spec(\bs{H}))\ge N^{-C}$ we see that Lemma~\ref{lmm:Quadratic large deviation} can be applied for a large deviation estimate on the $(u,v)$-sum in the definition \eqref{definition of ZBxy} of $Z^B_{xz}$ for $B=B_{xy}$, i.e. 
\bels{a priori bound on ZBxy}{
\abs{Z_{xz}^{B_{xy}}}
\,\prec\, 
\biggl(\frac{1}{N^2}\sum^{B_{xy}}_{u,v}\abs{G_{uv}^{B_{xy}}}^2\biggr)^{\!1/2}
\!\prec\, N^{\eps_1P}\Psi\,.
}
Here we also used \eqref{local law for GB} and \eqref{decay of MB} for the second stochastic domination bound. 
Combining \eqref{a priori bound on ZBxy} with \eqref{local law for GB} we see that
\vspace{-0.2cm}
\bels{bound on D4b}{
\norm{\bs{D}^{\rm{(4b)}}}_{\rm{max}}\,\prec\, N^{2\eps_1P}\Psi^2\,.
}

The rest of the proof of Proposition~\ref{prp:fluctuation averaging} is dedicated to showing the high moment bound
\bels{high moment bound on RD4a}{
\E\2\abs{\avg{\bs{R},\bs{D}^{\rm{(4a)}}}}^{2\mu}\,\lesssim_\mu\,N^{C(\mu)\eps_1}\Psi^{2\mu}.
}
Together with \eqref{Psi2 bound on Dk}, \eqref{Estimate on RD1}  and \eqref{bound on D4b} this bound implies \eqref{general FA} since $\eps_1$ can be chosen arbitrarily small. 
In analogy to \eqref{Splitting D1a into X+Y+Z} we write  $ \avg{\2\bs{R}\1,\bs{D}^{\rm{(4a)}}} = X+Y+Z $,
where the three terms on the right hand side are obtained by summing $ \sigma_{xzy}:=N^{-1}\2\ol{r}_{xy}Z_{x z}^{B_{xy}}\2m_{z y} $ over disjoint sets of indices:
\begin{equation*}
\begin{split}
X\,:=\,  \sum_{x}\sum_{y}^{B_x^1}\sum_{z \in B_{x}^1}^{B_y^1}\msp{-5}\sigma_{xzy}
\,,
\quad 
Y\,:=\,  \sum_{x}\sum_{y}^{B_x^1}\sum_{z \in B_{y}^1}\msp{-5}\sigma_{xzy}
\,,
\quad
Z\,:=\,  \sum_{x}\sum_{y\in B_x^1}\sum_{z \in A_{xy}}\msp{-5}\sigma_{xzy}\,.
\end{split}
\end{equation*}
Similar to \eqref{bounds on X Y Z for D1a} the fast decay of off-diagonal entries of both $\bs{R}$ and $\bs{M}$, and the a priori bound \eqref{a priori bound on ZBxy} immediately yield $ \abs{X} \prec_\mu N^{4\eps_1P\mu}  N^{-2\mu}\Psi^{2\mu} $. 
Since this is already sufficient for \eqref{high moment bound on RD4a}, we focus on the terms $ Y $ and $ Z $ in \eqref{bounds on X Y Z for D4a}. Using again the decay of off-diagonal entries yields:
\begin{subequations}
\label{bounds on X Y Z for D4a}
\begin{align}
\label{estimate Y for D4a}
\E\2\abs{Y}^{2\mu}\,&\lesssim_\mu 
\frac{1}{N^{4\mu}}\sum_{\bs{x},\bs{y}}\sum_{\bs{z} \in B_{\bs{y}}^1}
\absbb{\,\E\prod_{i=1}^\mu Z^{B_{x_iy_i}}_{x_iz_i}\ol{Z}^{B_{x_{\mu+i}y_{\mu+i}}}_{x_{\mu+i}z_{\mu+i}}}\,,
\\
\label{estimate Z for D4a}
\E\2\abs{Z}^{2\mu}\,&\lesssim_\mu 
\frac{1}{N^{2\mu}}\sum_{\bs{x}}\sum_{\bs{y} \in B_{\bs{x}}^1}\sum_{\bs{z} \in B_{\bs{x}}^1\cup B_{\bs{y}}^1}\absbb{\,\E\prod_{i=1}^\mu Z^{B_{x_iy_i}}_{x_iz_i}\ol{Z}^{B_{x_{\mu+i}y_{\mu+i}}}_{x_{\mu+i}z_{\mu+i}}}\,.
\end{align}
\end{subequations}
We call the subscripts $ i $ of the indices $ x_i $ and $ z_i $ \emph{labels}.
In order to further estimate the moments of $Y$ and $Z$ we introduce the set of \emph{lone labels} of $ (\bs{x},\bs{z}) $:
\bels{def of L}{
L(\bs{x},\bs{z}) \,:=\, \cB{i : d\pb{\{x_i,z_i\},{\textstyle \bigcup_{j \ne i}}\{x_j,z_j\} } \ge 3\1N^{\eps_1}
}
\,.
}
The corresponding index pair  $(x_i, z_i)$ for $ i\in L(\bs{x},\bs{z}) $, is called \emph{lone index pair}. We partition the sums in \eqref{estimate Y for D4a} and \eqref{estimate Z for D4a} according to the number of lone labels, i.e. we insert the partition of unity 
$
1 = \sum_{\ell=0}^{2\mu}\bbm{1}(\1\abs{L(\bs{x},\bs{z})}=\ell\2)
$.  
A simple counting argument reveals that fixing the number of lone labels reduces the combinatorics of the sums in  \eqref{estimate Y for D4a} and \eqref{estimate Z for D4a}. More precisely,
\bels{reduces combinatorics through lone labels}{
\sum_{\bs{x},\bs{y}}\sum_{\bs{z} \in B_{\bs{y}}^1}
\bbm{1}(\2\abs{L(\bs{x},\bs{z})}=\ell\,)
\,&\le\, 
N^{C \mu\2\eps_1}N^{2\mu+\ell}, 
\\
\sum_{\bs{x}}\sum_{\bs{y} \in B_{\bs{x}}^1}
\msp{-15}\sum_{\quad \bs{z} \in B_{\bs{x}}^1\cup B_{\bs{y}}^1}
\msp{-20}
\bbm{1}(\2\abs{L(\bs{x},\bs{z})}=\ell\,)
\,&\le\, 
N^{C \mu\2\eps_1}N^{\mu+\ell/2}
.
}
The expectation in \eqref{estimate Y for D4a} and \eqref{estimate Z for D4a} is bounded using the following technical result.

\begin{lemma}[Key estimate for averaged local law]
\label{lmm:Key estimate for averaged local law}
Assume the hypotheses of Proposition~\ref{prp:fluctuation averaging} hold, let $ \mu \in \N  $ and  $ \bs{x},\bs{y}  \in \bb{X}^{2\mu}$.
Suppose there are $ 2\mu $ subsets $ Q_1,\dots, Q_{2\mu} $ of $ \bb{X} $, such that $ B_{N^{\eps_1}}(x_i,y_i) \subseteq Q_i\subseteq B_{3N^{\eps_1}}(x_i,y_i)$ for each $ i $. Then
\bels{FA:resolvent product bound}{
\absbb{\,
\E\,
\prod_{i=1}^{\mu}
Z^{(Q_i)}_{x_iy_i}
\ol{Z}^{(Q_{\mu+i})}_{x_{\mu+i}y_{\mu+i}}
\msp{6}
}
\,&\lesssim_{\mu}\, 
N^{C(\mu)\eps_1}
\,
\Psi^{\12\mu\,+\2\abs{\1L(\bs{x},\1\bs{y})\1}}
.
}
\end{lemma}
Using \eqref{reduces combinatorics through lone labels} and Lemma~\ref{lmm:Key estimate for averaged local law} on the right hand sides of  \eqref{estimate Y for D4a} and \eqref{estimate Z for D4a} after partitioning according to the number of lone labels, yields
\bels{final moment bound on Y and Z for D1a}{
\E\2\abs{Y}^{2\mu}\,\lesssim_\mu \frac{N^{C(\mu)\eps_1}}{N^{4\mu}}\sum_{\ell=0}^{2\mu}\Psi^{2\mu+\ell}N^{2\mu+\ell},\qquad
\E\2\abs{Z}^{2\mu}\,\lesssim_\mu \frac{N^{C(\mu)\eps_1}}{N^{2\mu}}\sum_{\ell=0}^{2\mu}\Psi^{2\mu+\ell}N^{\mu+\ell/2}.
}
Since $\Psi\ge N^{-1/2}$ the high moment bounds in \eqref{final moment bound on Y and Z for D1a} together  with the simple estimate for $ X $  imply \eqref{high moment bound on RD4a}. This finishes the proof of Proposition~\ref{prp:fluctuation averaging} up to verifying Lemma~\ref{lmm:Key estimate for averaged local law} which will occupy the rest of the section.
\end{Proof}

\begin{Proof}[Proof of Lemma~\ref{lmm:Key estimate for averaged local law}]
Let us consider the data 
$
\xi := (\1\bs{x},\bs{y},(Q_i)_{i=1}^{2\mu}\1)
$ 
fixed. We start by writing the product on the left hand side of \eqref{FA:resolvent product bound} in the form.
\bels{Z in terms of w and Gamma}{
\prod_{i=1}^{\mu}
Z^{(Q_i)}_{x_iy_i}
\ol{Z}^{(Q_{\mu+i})}_{x_{\mu+i}y_{\mu+i}}
\,=\, %
\sum_{\bs{u},\bs{v}}
w_{\bs{x},\bs{y}}(\bs{u},\bs{v})\,
\Gamma_{\!\xi}
(\bs{u},\bs{v})
\,,
}
where the two auxiliary functions $ \Gamma_{\!\xi}, w_{\bs{x},\bs{y}} : \bb{X}^{2\mu} \times \bb{X}^{2\mu} \to \C $, are defined by
\begin{subequations}
\label{defs of w and Gamma}
\begin{align}
\label{def of Gamma at level k=0}
\Gamma_{\!\xi}
(\bs{u},\bs{v}) 
\,&:=\, 
\bbm{1}\cb{u_i,v_i\notin Q_i,\forall\2i=1,\dots,2\mu} 
\prod_{i=1}^{\mu}
G^{(Q_i)}_{u_i v_i}
\ol{G}^{(Q_{\mu+i})}_{u_{\mu+i} v_{\mu+i}}\,,
\\
\label{def of product w_xy}
w_{\bs{x},\bs{y}}(\bs{u},\bs{v}) 
\,&:=
\prod_{i=1}^{\mu}
w_{x_iy_i}\msp{-2}(u_i,v_i)
\,\ol{w}_{x_{\mu+i}y_{\mu+i}}(u_{\mu+i},v_{\mu+i})
\,,
\end{align}
\end{subequations}
and 
\vspace{-0.5cm}
\bels{def of w_xy(i,j)}{
w_{xy}(u,v) \,:=\,\frac{1}{N}(w_{xu}w_{vy}-\E\,w_{xu}w_{vy})
\,.
}

In order to estimate \eqref{Z in terms of w and Gamma} we partition the sum over the indices $ u_i $ and $ v_i $  depending on their distance from the set of {lone index} pairs,  $ (x_i,y_i) $ with $ i \in L $, where $ L = L(\bs{x},\bs{y}) $.
To this end we introduce the partition $ \c{B_i:i \in \c{0}\cup L} $ of $ \bb{X} $, %
\bels{def of the parts B_s}{
B_i \,:=\,
\begin{cases}
\;B_{N^{\eps_1}}(x_i)\cup B_{N^{\eps_1}}(y_i)
&\text{when }i \in L\,,
\\
\bb{X}\2\backslash\bigcup_{j\in L} (B_{N^{\eps_1}}(x_j)\cup B_{N^{\eps_1}}(y_j))
&\text{when }i = 0\,, 
\end{cases} 
}
and the shorthand
\bels{def of BB(xi,sigma)}{
\mathbb{B}(\xi,\sigma) 
\,:=\, \cB{(\bs{u},\bs{v}) \in \bb{X}^{4\mu}: u_i \in B_{\sigma'_i}\!\backslash Q_i\2,\;v_i \in B_{\sigma''_i}\!\backslash Q_i\2,\,i=1,\dots,2\mu} 
\,,
}
where the components $ \sigma_i = (\sigma'_i,\sigma''_i) \in (\1\c{0} \1\cup\2 L\1)^2 $ of $ \sigma = (\sigma_i)_{i=1}^{2\mu} $ specify whether $u_i$ and  $v_i$ are close to a lone index pair or not; e.g. $\sigma'_i$ determines
which lone index $u_i$ is close to, if any.  For any fixed $\xi$, as $\sigma$ runs through all 
possible elements of  $(\1\c{0} \1\cup\2 L\1)^{4\mu}$,   the sets $\mathbb{B}(\xi,\sigma) $
form a partition  of the summation set on the right hand side of  \eqref{Z in terms of w and Gamma} (taking into account
the  restriction $u_i, v_i\not\in Q_i$). Therefore it will be sufficient to estimate
\bels{fixed sigma}{
   \sum_{\quad(\bs{u},\bs{v}) \in  \mathbb{B}(\xi,\sigma) }
\msp{-20}
w_{\bs{x},\bs{y}}(\bs{u},\bs{v})\,
\Gamma_{\!\xi}
(\bs{u},\bs{v})
}
for every fixed $\sigma\in(\1\c{0} \1\cup\2 L\1)^{4\mu}$. 
Since $ x_i $ and $ y_i $ are fixed, while $ u_i $ and $ v_i $ are free variables, with their domains depending on $ (\xi,\sigma) $, we say that the former are \emph{external} indices and the latter  are \emph{internal} indices. 

Let us define the set of \emph{isolated labels},
\bels{def of wht-L}{ 
\wh{L}(\bs{x},\bs{y},\sigma) \,=\, L(\bs{x},\bs{y})\1 \backslash\1 \c{\1\sigma'_1,\dots,\sigma'_{2\mu},\sigma''_1,\dots,\sigma''_{2\mu}\1}
\,,
}
so that if an external index has an isolated label as subscript, then it is isolated from all the other indices in the following sense: 
\[
\qquad
d\pB{\1\c{x_i,y_i}\2,\,{\textstyle \bigcup_{j=1}^{\12\mu}}\c{u_j,v_j}  \cup\, {\textstyle \bigcup_{j\ne i}}\c{x_j,y_j}}
\,\ge\, 
N^{\eps_1}\,,
\quad
(\bs{u},\bs{v})  \in \mathbb{B}(\xi,\sigma)
\,,
\; i \in \wh{L}(\bs{x},\bs{y},\sigma)
\,.
\]
Notice that isolated labels indicate not only separation from all other external indices, as lone labels do,  but also from 
all internal indices. 
Given a resolvent entry $ G^{B}_{uv} $ we will refer to $ u,v $ as \emph{lower indices} and the set $ B $ as an \emph{upper index set}.

The next lemma, whose proof we postpone until the end of this section, yields an algebraic representation for \eqref{fixed sigma} provided the internal indices are properly restricted. 

\begin{lemma}[Monomial representation]
\label{lmm:Monomial representation}
Let $ \xi $ and $ \sigma $ be fixed. 
Then the restriction $ \Gamma_{\!\xi}|_{\mathbb{B}(\xi,\sigma)} $ of the function \eqref{def of Gamma at level k=0} to the subset $ \mathbb{B}(\xi,\sigma) $ of $ \bb{X}^{4\mu}$ has a representation
\bels{monomial exp for Gamma on BB}{
\Gamma_{\!\xi}|_{\mathbb{B}(\xi,\sigma)} 
\,=\msp{-4} 
\sum_{\alpha\1=\11}^{\; M(\xi,\sigma)} \Gamma_{\!\xi,\1\sigma,\alpha}
\,,
}
in terms of 
\vspace{-0.5cm}
\bels{FA:bound on M}{
M(\xi,\sigma) \,\lesssim_\mu\,  N^{C(\mu)\eps_1},%
}
(signed) monomials $ \Gamma_{\!\xi,\sigma,\alpha} :\mathbb{B}(\xi,\sigma) \to \C $, such that $ \Gamma_{\!\xi,\sigma,\alpha}(\bs{u},\bs{v}) $ for each $\alpha$ is of the form:
\bels{typical monomial V2}{
(-1)^\#\prod_{t=1}^{n}
(\2G^{E_{t}}_{ a_{ t} b_{ t}})^\#
\,\prod_{r=1}^{q}
\frac{\!1}{(\2G^{F_{r}}_{w_{ r} w_{r}})^\#\!}
\prod_{\msp{10} r\in R^{(1)}}\msp{-5}
(\2G^{U_{r}}_{u_rv_r})^\#\msp{-5}
\prod_{\msp{10} t\in R^{(2)}}\msp{-5}
(\2G^{U'_{t}}_{u_{t}u'_{ t}}G^{V'_{t}}_{v'_{t}v_t})^\#
\,. 
}
Here the notations $ (-1)^\# $ and $ (\,\cdot\,)^\# $ indicate possible signs and complex conjugations that may depend only on $ (\xi,\sigma,\alpha) $, respectively, and that will be irrelevant for our estimates.
The dependence on $ (\xi,\sigma,\alpha) $ has been suppressed in the notations, e.g.,  $ n = n(\xi,\sigma,\alpha) $, $U_r=U_r(\xi,\sigma,\alpha)$, etc.

The numbers $ n $ and $ q $ of factors in \eqref{typical monomial V2} are bounded, $ n+q \lesssim_\mu 1 $.
Furthermore, for any fixed $\alpha$ the two subsets $ R^{(k)}$, $ k=1,2$, form a partition of $\{1,\dots,2\mu\}$, and the monomials \eqref{typical monomial V2} have the following  three  properties:
\begin{enumerate}
\item 
The lower indices $ a_{t} $, $ b_{t} $, $ u'_{t}$, $v'_{t}$, $ w_{t} $ are in $ \cup_{i\in \wh{L}} B_i $, and $  d(a_{t},b_{t}) \ge N^{\eps_1} $.
\item
The upper index sets $ E_{r}$, $ F_{r}$, $ U_{r}$, $ U'_{r}$, $ V'_{r} $ are bounded in size by $ N^{C(\mu)\eps_1} $, and $ B_r \subseteq U_{r},U'_{r},V'_{r} $.
The total number of these sets appearing in the expansion \eqref{monomial exp for Gamma on BB} is bounded by $ N^{C(\mu)\eps_1} $.  
\item
At least one of the following two statements is always true:
\vspace{-0.10cm}
\bea{
&\,\emph{(I)}\quad \exists\1 i \in \wh{L}\,,
\;\text{ s.t. }\;
B_i \;\subseteq\, 
{\textstyle \bigcap_{\1t=1}^{n}} E_{t} 
\,\cap\, 
{\textstyle \bigcap_{\1r=1}^{q}} F_{r} 
\,\cap\, 
{\textstyle \bigcap_{\1r\in R_1}}\!
U_{r} \,\cap\,
{\textstyle \bigcap_{\1 t\in R_2}}
(\,
U'_{t} \cap V'_{t})
\;;
\\ 
\label{FA: lower bound on regular resolvent entries}
&\emph{(II)}\quad n+ \abs{R^{(1)}} + 2\1\abs{R^{(2)}} 
\;\ge\; 2\mu + \abs{\1\wh{L\1}}
\,.
}
\end{enumerate}
\end{lemma}

Since  Lemma~\ref{lmm:Key estimate for averaged local law}  relies heavily on this representation, we make a few remarks:
(i) Monomials with different values of $ \alpha $ may be equal. The indices $ a_{t} $, $ b_{t} $, $ u'_{t}$, $v'_{t}$, $ w_{t} $ may overlap, but they are always distinct from the internal indices since from \eqref{def of BB(xi,sigma)} and \eqref{def of wht-L} we see that 
\[
\c{u_r,v_r}_{r=1}^{2\mu} \subseteq \bb{X}  \backslash \bigl(\cup_{\!i\2\in\2 \wh{L}}B_i\bigr) 
\,,\qquad
(\bs{u},\bs{v}) \in \mathbb{B}(\xi,\sigma) 
\,.
\]
(ii) The reciprocals of the resolvent entries are not important for our analysis because the diagonal resolvent entries are comparable to  $ 1 $ in absolute value when  a local law holds (cf. \eqref{Bounds on mB}). 
(iii) Property 3 asserts that each monomial  is either a deterministic function of $ \bs{H}^{(B_i)} $ for some isolated label $ i $, and consequently almost independent of the rows/columns of $ \bs{H} $ labeled by $x_i,y_i $ (Case (I)), or the monomial contains at least $ \abs{\wh{L}} $ additional off-diagonal resolvent factors (Case (II)). 
In the second case, each of these  extra factors will provide  an additional factor $\Psi$ for typical internal indices  due to faster than power law decay of $ \bs{M} $ and  the local law  \eqref{local law for GB}.
Atypical internal indices, e.g. when $u_r$ and $v_r$ are close to each other, do not
give a factor $\Psi$ since $m_{u_r v_r}$ is not small, but there are  much fewer atypical indices  than  typical ones and
this entropy factor makes up for the lack of smallness.
These arguments will be made rigorous in Lemma~\ref{lmm:Three sources of smallness} below.

By using the monomial sum representation \eqref{monomial exp for Gamma on BB} in \eqref{fixed sigma}, and estimating each summand 
 separately, we obtain
\bels{mcl-Z bound 1}{
&\Biggl|
\,\E\!
\prod_{i=1}^{\mu}
Z^{(Q_i)}_{x_iy_i}
\ol{Z}^{(Q_{\mu+i})}_{x_{\mu+i}y_{\mu+i}}\Biggr|
\,\lesssim_\mu\;
N^{\1C(\mu)\eps_1}
\max_{\sigma}
\max_{\alpha\,=\,1}^{M(\xi,\sigma)}\,
\Biggl|\,
\E
\msp{-23}\sum_{\quad(\bs{u},\bs{v}) \in \mathbb{B}(\xi,\sigma)}
\msp{-25}w_{\bs{x},\bs{y}}(\bs{u},\bs{v})\,
\Gamma_{\!\xi,\sigma,\alpha}(\bs{u},\bs{v})
\,
\Biggr|
\,,
}
where the factor $N^{\1C(\mu)\eps_1} $ originates from \eqref{FA:bound on M}, and we have bounded the summation over by a $\mu$-dependent constant.
Thus \eqref{FA:resolvent product bound} holds if we show, uniformly in $ \alpha $, that
\bels{FA:bound on monomials needed}{
\Biggl|\,
\E
\msp{-23}\sum_{\quad(\bs{u},\bs{v}) \in \mathbb{B}(\xi,\sigma)}
\msp{-25}w_{\bs{x},\bs{y}}(\bs{u},\bs{v})\,
\Gamma_{\!\xi,\sigma,\alpha}(\bs{u},\bs{v})
\,
\Biggr|
\;&\leq\;
N^{\1C(\mu)\eps_1}
N^{-\frac{1}{2}(\2\abs{\1L(\bs{x},\bs{y})}\,-\,\abs{\1\wh{L}(\bs{x},\bs{y},\sigma)}\2)}
\,
\Psi^{\22\mu\,+\2\abs{\wh{L}(\bs{x},\bs{y},\sigma)\1}}
\,.
}
In order to prove this bound, we fix $ \alpha $, and sum over the internal indices to get
\bels{FA: smallness can be extracted}{
\msp{-10}\Biggl|\;
\E
\msp{-23}\sum_{\quad(\bs{u},\bs{v}) \in \mathbb{B}(\xi,\sigma)}
\msp{-25}w_{\bs{x},\bs{y}}(\bs{u},\bs{v})\,
\Gamma_{\!\xi,\sigma,\alpha}(\bs{u},\bs{v})
\,
\Biggr|
\,
\leq\;
\E\,
\prod_{t=1}^{n}\,
\abs{\2G^{E_{t}}_{ a_{t} b_{t}}}
\msp{10}\prod_{r=1}^{q}\;
\abs{\2G^{F_{ r}}_{w_{ r} w_{r}}}^{-1}
\msp{-10}
\prod_{\quad r\in R^{(1)}} \msp{-14}\Theta^{(1)}_r
\msp{-12}\prod_{\quad r\in R^{(2)}} \msp{-14}\Theta^{(2)}_r
\,,
}
where we have used  the formula \eqref{typical monomial V2} for the monomial $\Gamma_{\!\xi,\sigma,\alpha}$. 
The sums over the internal indices have been absorbed into the following factors:
\bels{defs of Theta^(1) and Theta^(2)}{
\qquad \Theta^{(1)}_r
\,&:=\,
\Biggl|\msp{-16}
\sum_{\quad u\in B_{\sigma'_r}\!\backslash Q_r}
\msp{-13}
\sum_{\quad v\in B_{\sigma''_r}\!\backslash Q_r}
\msp{-20}
w_{{x_r},y_r}\msp{-2}(u,v)\;
G^{U_{r}}_{uv}
\,\Biggr|\,,
\msp{70}
r \in R^{(1)},
\\
\qquad \Theta^{(2)}_r
\,&:=\,
\Biggl|\msp{-16}
\sum_{\quad u\in B_{\sigma'_r}\!\backslash Q_r}
\msp{-13}
\sum_{\quad v\in B_{\sigma''_r}\!\backslash Q_r}
\msp{-20}
w_{{x_r},y_r}\msp{-2}(u,v)\;
G^{U'_{r}}_{u\1u'_{r}}G^{V'_{r}}_{v'_{r}v}
\;
\Biggr|
\,,
\quad
r\in R^{(2)}
.
}

The right hand side of \eqref{FA: smallness can be extracted} will be bounded using the following three estimates which follow by combining the monomial representation with our previous stochastic estimates.

\begin{lemma}[Three sources of smallness]
\label{lmm:Three sources of smallness}
Consider an arbitrary monomial  $ \Gamma_{\!\xi,\sigma,\alpha} $, of the form \eqref{typical monomial V2}. 
Then, under the hypotheses of Proposition~\ref{prp:fluctuation averaging}, 
the following three estimates hold:
\begin{enumerate}
\item
The resolvent entries with no internal lower indices are small while the reciprocals of the resolvent entries are bounded, in the sense that 
\bels{LL for monomial factors}{ 
\abs{\2G^{E_{t}}_{ a_{t} b_{ t}}} 
\,&\prec\,
\Psi
\,,\qquad
\abs{\2 G^{F_{ r}}_{w_{r} w_{r}}}^{-1} 
\prec\;
1
\,.
}
\item
If $ \Gamma_{\!\xi,\sigma,\alpha} $ satisfies (I) of Property 3 of Lemma~\ref{lmm:Monomial representation}, then its contribution is very small in the sense that
\bels{too small monomial}{
\abs{\2\E\,w_{\bs{x},\bs{y}}(\bs{u},\bs{v})\,\Gamma_{\!\xi,\sigma,\alpha}(\bs{u},\bs{v})\1}
\,&\lesssim_{\mu,\nu}\;
N^{-\nu}  
\,,
\qquad 
(\bs{u},\bs{v})\2\in\2\mathbb{B}(\xi,\sigma)
\,.
}
\item
Sums over the internal indices around external indices with lone  labels yield extra smallness: 
\vspace{-0.3cm}
\bels{Theta^(k)_r estimates}{
\qquad\quad
{\Theta^{(k)}_r}
\;&\prec\,
N^{C(\mu)\eps_1}
N^{\1-\1\frac{\11\1}{2}\1\abs{\1\sigma_r}_\ast}
\Psi^{\1k}
,
\qquad 1 \leq r \leq 2\mu\,,\; k=1,2 
\,,
}
where $ \abs{\sigma_r}_\ast \,:=\, \abs{\c{\10\1,\sigma'_r,\sigma''_r}} - 1 $ counts how many, if any, of the two indices $ u_r $ and $ v_r $, are restricted to vicinity of distinct external indices.
\end{enumerate}
\end{lemma}

We postpone the proof of Lemma~\ref{lmm:Three sources of smallness} and first see how it is used to finish the proof of Lemma~\ref{lmm:Key estimate for averaged local law}.
The bound \eqref{FA:bound on monomials needed} follows by combining Lemma~\ref{lmm:Monomial representation} and Lemma~\ref{lmm:Three sources of smallness} to estimate the right hand side of \eqref{FA: smallness can be extracted}.
If (I) of Property 3 of Lemma~\ref{lmm:Monomial representation}  holds, then applying \eqref{too small monomial} and \eqref{LL for monomial factors}
in \eqref{FA: smallness can be extracted} yields \eqref{FA:bound on monomials needed}.
On the other hand, if  (I) of Property 3 of Lemma~\ref{lmm:Monomial representation} is not true, then we use \eqref{LL for monomial factors} and \eqref{Theta^(k)_r estimates} to get
\bels{FA final form before completion}{
\prod_{t=1}^{n}\,
\abs{\2G^{E_{t}}_{ a_{t} b_{ t}}}
\msp{10}\prod_{r=1}^{q}\,
\abs{\2G^{F_{r}}_{w_{ r} w_{r}}}^{-1}
\msp{-13}
\prod_{\quad r\in R^{(1)}} \msp{-10}\Theta^{(1)}_r
\msp{-12}\prod_{\quad r\in R^{(2)}} \msp{-10}\Theta^{(2)}_r
\;\prec\;N^{C(\mu)\eps_1}
\Psi^{\1n\2+\2\abs{R^{(1)}}\2+\22\1\abs{R^{(2)}}}\,N^{\1-\1\frac{1}{2}\sum_r \abs{\1\sigma_r}_\ast}
\,.
}
By Part 3 of Lemma~\ref{lmm:Monomial representation} we know that (II) holds. 
Thus the power of $ \Psi $ on the right hand side of \eqref{FA final form before completion} is at least $ 2\mu + \abs{\wh{L}} $. 
On the other hand, from \eqref{def of wht-L} we see that
\[
\abs{\1L\1}-\abs{\1\wh{L}\1}
\,\leq\, 
\absb{\,{\textstyle \bigcup_{\1r=1}^{\12\1\mu}} \c{\sigma'_r, \sigma''_r}\backslash\c{0}}
\,\leq\, 
\sum_{r=1}^{2\mu} \absb{\2\c{\sigma'_r, \sigma''_r}\backslash\c{0}} 
\,\leq\, 
\sum_{r=1}^{2\mu}\abs{\sigma_r}_*
\,.
\]
Hence the power of $ N^{-1/2} $ on the right hand side of \eqref{FA final form before completion} is at least $ \abs{L}-\abs{\wh{L}} $.
Using these bounds together with $\Psi \ge N^{-1/2}$  %
 in \eqref{FA final form before completion}, and then taking expectations yields \eqref{FA:bound on monomials needed}.
Plugging \eqref{FA:bound on monomials needed} into \eqref{mcl-Z bound 1} completes the proof of \eqref{FA:resolvent product bound}.   
\end{Proof}

\begin{Proof}[Proof of Lemma~\ref{lmm:Three sources of smallness}]
 Combining \eqref{local law for GB} and \eqref{decay of MB} we see that for some sequence $ \alpha $ 
\bels{scaling of resolvent entries with non-trivial upper index sets}{
\abs{G^{E}_{uv}}  
\,\prec\;
N^{C\eps_1}\Psi 
+ \frac{\alpha(\nu)}{(1+d(u,v))^\nu}
\,,\quad\text{whenever}
\quad
u,v \notin E\,,\;\text{ and }\;
\abs{E} \leq N^{C\eps_1}
\,.
}
By the bound on the size of $ E_{t}$, $ F_{r} $ in Property~2 of Lemma~\ref{lmm:Monomial representation}, \eqref{scaling of resolvent entries with non-trivial upper index sets} is applicable for these upper index sets. Then  \eqref{LL for monomial factors} follows from the second bound of Property 1 of Lemma~\ref{lmm:Monomial representation} and the decay of the entries of $\bs{M}^{E}$ from \eqref{decay of MB}.

In order to prove Part 2, let $ i \in \wh{L} $ be the label from (I) of Property 3 of Lemma~\ref{lmm:Monomial representation}. 
We have %
\bea{
\E\,w_{\bs{x},\bs{y}}(\bs{u},\bs{v})\,\Gamma_{\!\xi,\sigma,\alpha}(\bs{u},\bs{v})
\;&=\;
\E\Bigl[\2w_{x_iy_i}\msp{-2}(\1u_i,v_i)
^\#\,\Bigr]\,\cdot\, \E\biggl[\;
\Gamma_{\!\xi,\sigma,\alpha}(\bs{u},\bs{v})
\;\prod_{r\1\neq\1 i}w_{x_ry_r}\msp{-2}(u_r,v_r)^\#%
\,\biggr]
\\
&\msp{24}+\,
\rm{Cov}\biggl(\, 
w_{x_iy_i}\msp{-2}(u_i,v_i)^\#,\;\Gamma_{\!\xi,\sigma,\alpha}(\bs{u},\bs{v})\,
\prod_{r\1\neq\1 i}w_{x_ry_r}\msp{-2}(u_r,v_r)^\#\biggl)
\,,
}
where the first term on the right hand side vanishes because $ w_{xy}(u,v) $'s are centred random variables by  \eqref{def of w_xy(i,j)}. 
Now the covariance is smaller than any inverse power of $N$, since $ w_{x_iy_i}\!(u_i,v_i) $ depends only on 
the $ x_i $-th and $y_i$-th row/column of $ \bs{H} $, while 
$ \Gamma_{\!\xi,\sigma,\alpha} $  is a deterministic function of $ \bs{H}^{B_i} $ by  (I) of Property~3 of Lemma~\ref{lmm:Monomial representation}.
Indeed, the faster than power law decay of correlations \eqref{decay of correlation}  yields \eqref{too small monomial}, because the derivative of $\Gamma_{\!\xi,\sigma,\alpha}(\bs{u},\bs{v})$ with respect to the entries of $\bs{H}$ are bounded in absolute value by $ N^{C(\mu)} $ by the $N$-dependent smoothness of the resolvents $\bs{G}^E(\zeta)$ as a function of $\bs{H}$ for spectral parameters $\zeta$ with $\dist(\zeta, \spec(\bs{H}^{E}))\ge N^{-C}$. For more details we refer to the proof of Lemmas~\ref{lmm:Quadratic large deviation} and \ref{lmm:Linear large deviation}, where a similar argument was used. 

Now we will prove Part 3.
To this end, fix an arbitrary label $ r =1,2,\dots, 2\mu $. 
Let us denote $ B_L := \bigcup_{s\in L} B_s $ and $ B_{\wh{L}} := \bigcup_{t \in \wh{L}} B_t$.

Let us first consider $ \Theta^{(1)}_r $. If $ \sigma'_r = s $ and $ \sigma''_r = t $, then we need to estimate
\bels{estimating Theta^(1)}{
\sum_{\;u \in B_s\!\backslash Q_r}\sum_{\;v \in B_t\!\backslash Q_r} \msp{-6}w_{x_ry_r}\msp{-3}(u,v) \,G^{U_r\msp{-1}}_{uv}
\,,\quad\text{where }\,
s,t \in L\backslash \wh{L}
\,,\text{ and }
 Q_r \subseteq U_r \subseteq Q_r\cup B_{\wh{L}}
\,.
}
Since $ B_s\!\backslash Q_r, B_t\!\backslash Q_r \subseteq \bb{X}\1\backslash U_r $, the indices $ u,v $ do not overlap the upper index set $ U_r $.
Hence, in the case $ k=1$ and $ s = t = 0 $ the estimate \eqref{Theta^(k)_r estimates} follows from \eqref{Double sum large deviation} of Lemma~\ref{lmm:Quadratic large deviation}.

If $ s,t \in L $, then taking modulus of \eqref{estimating Theta^(1)} and using \eqref{scaling of resolvent entries with non-trivial upper index sets} yields \eqref{Theta^(k)_r estimates}:
\bels{Theta^1 estimate when s,t not 0}{
\Theta^{(1)}_r \,&\leq\,\abs{B_s\!\backslash Q_r}\,\abs{B_t\!\backslash Q_r}\,
\Bigl(\max_{\;u\1,\2v\2\in\2 \bb{X}}\, \abs{w_{x_ry_r}(u,v)}\Bigr) 
\max_{\;u \in B_s\!\backslash Q_r}\!\max_{\;v \in B_t\!\backslash Q_r}
\abs{G^{U_r}_{uv}}%
\\
&
\prec\;
\frac{\,N^{C\eps_1}\!}{N}\Bigl( 
N^{C\eps_1}\Psi + \frac{\alpha(\nu)}{(1+d(B_s,B_t))^\nu}\Bigr)
\,\leq\,
N^{C\eps_1} N^{-\frac{1}{2}\abs{\sigma_r}_\ast}\Psi 
\,,
}

where $ d(A,B) := \inf \sett{d(a,b):a\in A,\;b \in B} $ for any sets $ A $ and $ B $ of $ \bb{X} $.
Here we have also used the definition \eqref{def of L} of lone labels and $ \Psi \ge N^{-1/2} $.

Suppose now that exactly one component of $ \sigma_r $ equals $ 0 $ and one is in $ L $.  
In this case, we split $ w_{x_ry_r}(u,v) $ in \eqref{estimating Theta^(1)} into two parts corresponding to $ w_{x_ru}w_{v\1y_r} $ and its expectation, and estimate the corresponding sums separately.
First, using \eqref{Single sum large deviation} of Lemma~\ref{lmm:Linear large deviation} yields 
\bels{FA:non-average part for k=1}{
\msp{-14}\frac{1}{N}
\absbb{\msp{-6}\sum_{\;u \in B_s\!\backslash Q_r}\msp{-5}\sum_{\;v\in B_0\!\backslash Q_r} \msp{-6}
w_{x_ru}w_{v\1y_r}
G^{U_r}_{uv}}
&\prec
\frac{\abs{B_s\!\backslash Q_r}}{N}\,
\biggl(\max_{u\in \bb{X}}\, \abs{w_{{x_r}u}}\biggr)
\max_{u\2\notin\1 U_r}
\absbb{
\msp{-19}
\sum_{\quad v \in B_0\!\backslash Q_r} \msp{-16}
G^{U_r}_{uv}w_{v\1y_r}
}
\prec
\frac{N^{C\eps_1}\Psi}{N^{1/2}}
.
}
On the other hand, using \eqref{scaling of resolvent entries with non-trivial upper index sets} we estimate the expectation part:
\bels{FA:average part for k=1}{
\msp{-5}
\frac{1}{N}\absbb{\msp{-6}\sum_{\;u \in B_s\!\backslash Q_r}\sum_{\;v \in B_0\!\backslash Q_r} \msp{-6}
(\2\E\,w_{x_ru}w_{v\1y_r})\,
G^{U_r}_{uv}}
\!\prec
\frac{\abs{B_s\!\backslash Q_r}}{N}\, 
\max_{u\in \bb{X}}\msp{-10}
\sum_{\;v \in B_0\!\backslash Q_r} \msp{-6}
\abs{\E\,w_{x_ru}w_{v\1y_r}}\bigl(N^{C\eps_1}\Psi+\abs{m^{U_r}_{uv}}\bigr)
.
}
Similar to the part \eqref{FA:non-average part for k=1}, because of \eqref{decay of covariances},  we can estimate \eqref{FA:average part for k=1} by  $ \ord_{\!\prec}(N^{C\eps_1} N^{-1}) $. 
As $ \Psi \ge N^{-1/2} $, this  finishes the proof of \eqref{Theta^(k)_r estimates} in the case $ k = 1 $.

Now we prove \eqref{Theta^(k)_r estimates} for $\Theta^{(2)}_r $.
In this case, we need to bound,
\bels{estimating Theta^(2)}{
\sum_{\;u \in B_s\!\backslash Q_r}
\msp{-5}
\sum_{\; v \in B_t\!\backslash Q_r} \msp{-10}w_{x_ry_r}\msp{-3}(u,v) \,G^{U'_r}_{u\2u'_r}G^{V'_r}_{v'_r\1v}
\,,
}
where $ s = \sigma'_r$, $ t=\sigma''_r  $ have again values in $ \c{0} \cup L\backslash \wh{L} $. Here, 
$u'_r \in B_{\wh{L}} \!\backslash U'_r $, $ v'_r \in B_{\wh{L}}\backslash V'_r $, and 
$ Q_r \subseteq\, U'_r\2,V'_r \subseteq Q_r \msp{-2}\cup\msp{-2} B_{\wh{L}} $.

By definitions of the lone and isolated labels \eqref{def of L} and \eqref{def of wht-L}, respectively, we know that, if $  s \in L\backslash \wh{L}$, then $ d(\1u\1,u'_r) \ge N^{\eps_1} $, and similarly, if 
$ t \in L\backslash \wh{L}$, then $ d(v'_r,v\1) \ge N^{\eps_1} $.
Thus, if $ s,t \in L\backslash \wh{L} $, then estimating similarly as in \eqref{Theta^1 estimate when s,t not 0} with \eqref{scaling of resolvent entries with non-trivial upper index sets}, yields
\[
\Theta^{(2)}_r 
\,\prec\,
N^{C\eps_1} N^{-1}\Psi^2\,, 
\qquad
s,t \in L\backslash \wh{L}
\,.
\]

In the remaining cases, we split \eqref{estimating Theta^(2)} into two parts corresponding to the term $ w_{x_ru}w_{v\1y_r} $ and its expectation in the definition of \eqref{def of w_xy(i,j)} of $ w_{x_ry_r}(u,v)$, and estimate these two parts separately.

The average part is bounded similarly as in \eqref{FA:average part for k=1}, i.e., if $ s \in L\backslash \wh{L} $ and $ t=0$, then
\begin{align}
\label{FA:average part for k=2}
&\frac{1}{N}\absbb{\sum_{\;u\in B_s\!\backslash Q_r}\sum_{\;v\in B_0\!\backslash Q_r} \msp{-6}
(\2\E\,w_{x_ru}w_{v\1y_r})\,
G^{U'_r}_{u\2u'_r}G^{V'_r}_{v'_r\1v}}
\\ \nonumber
&\msp{10}
\prec\;
\frac{\abs{B_s\!\backslash Q_r}}{N}
\max_{u\in B_s}\msp{-10}
\sum_{\;v \in B_0\!\backslash Q_r} \msp{-6}
(\2\E\,w_{x_ru}w_{v\1y_r})
\biggl( 
N^{C\eps_1}\Psi + \frac{\alpha(\nu)}{(1+d(u,u'_r))^\nu}\biggr)
\biggl( 
N^{C\eps_1}\Psi + \frac{\alpha(\nu)}{(1+d(v'_r,v))^\nu}\biggr)
\,.
\end{align}
Here $ d(u,u'_r) \ge N^{\eps_1} $ since $ u \in B_s $, $ s \in L\backslash \wh{L} $, while $ u'_r \in B_{\wh{L}} $.
Taking $ \nu > C\eps_1^{-1}$ and using the \eqref{decay of covariances} to bound the sum over the covariances by a constant, we thus we see that the right hand side is $ \ord_{\!\prec}(N^{C\eps_1}N^{-1}\Psi\1) $. 
Since $ \Psi \ge N^{-1/2} $, this matches \eqref{Theta^(k)_r estimates} as $ \abs{\sigma_r}_\ast = \abs{\c{0,s,t}}-1 = 1 $.

Now, we are left to bound the size of terms of the form \eqref{estimating Theta^(2)}, where $w_{x_ry_t}(u,v) $ is replaced with  $ \frac{1}{N}w_{x_ru}w_{vy_r} $, and either $ s = 0 $ or $ t = 0 $. 
In these cases the sums over $ u $ and $ v $ factorize, i.e., we have
\vspace{-0.3cm}
\[
\frac{1}{N}
\biggl(
\msp{-12}
\sum_{\quad u \in B_s\!\backslash Q_r}
\msp{-15}
w_{x_ru}G^{U'_r}_{u\2u'_r}
\biggr)
\biggl(\msp{-12}
\sum_{\quad v \in B_t\!\backslash Q_r} \msp{-15}
G^{V'_r}_{v'_r\1v}w_{v\1y_r}
\biggr)
\,.
\]
When the sum is over a small set, i.e., over $ B_{s'} $ for some $ s' \in L\backslash \wh{L} $, then we estimate the sizes of the entries of $ \bs{W} $ and $ \bs{G}^{(\#)} $ by $ \ord_{\!\prec}(N^{-1/2})$ and $ \ord_{\!\prec}(\Psi) $, respectively. 
On the other hand, when $ u $ or $ v $ is summed over $ B_0\!\backslash Q_r $, we use \eqref{Single sum large deviation} of Lemma~\ref{lmm:Linear large deviation} to obtain a bound of size $ \ord_{\!\prec}(\Psi) $. In each case, we obtain an estimate that matches \eqref{Theta^(k)_r estimates}. 
\end{Proof}

\begin{Proof}[Proof of Lemma~\ref{lmm:Monomial representation}]
We consider the data $ (\xi,\sigma) $ fixed, and write $ \wh{L} = \wh{L}(\bs{x},\bs{y},\sigma) $, etc.
We start by enumerating the isolated labels (see \eqref{def of wht-L})
\bels{enumeration of wht-L}{ 
\c{s_1,\dots,s_{\wh{\ell}}\,} \,=\,\wh{L}
\,,\qquad
\wh{\ell} := \abs{\wh{L}}
\,,
}
and set $ \wh{B}(k) := \cup_{j=1}^k B_{s_j} $ for $ 1 \leq k \leq \wh{\ell}$
(recall the definition from \eqref{def of the parts B_s} and that $B_{s_j}$'s are disjoint).

The monomial expansion \eqref{monomial exp for Gamma on BB} is constructed iteratively in  $ \wh{\ell} $ steps. Indeed, we will define $ 1+ \wh{\ell} $ representations,
\vspace{-0.2cm}
\bels{level-k expansion}{
\msp{100}\Gamma_{\!\xi}|_{\mathbb{B}(\xi,\sigma)} %
\,=\,
\sum_{\alpha=1}^{M_k} \,\Gamma^{(k)}_\alpha
\,,\qquad
k=0,1,\dots,\wh{\ell}
\,.
}
where the $ M_k = M_k(\xi,\sigma) $ monomials $  \Gamma^{(k)}_{\!\alpha}  = \Gamma^{(k)}_{\!\xi,\sigma,\alpha} : \mathbb{B}(\xi,\sigma) \to \C $, evaluated at $ (\bs{u},\bs{v}) \in \mathbb{B}(\xi,\sigma)$, are of the form
\vspace{-0.3cm}
\bels{monomial in simple standard form}{
(-1)^{\#}\prod_{t=1}^m (\1G^{E_t}_{a_t b_t})^{\#} \prod_{r=1}^q \frac{1}{(\1G^{F_r}_{w_r w_r})^{\#}\msp{-6}}
\,,
}
with some indices $a_t,b_t \notin E_t $, $ w_r \notin F_r $. 
The numbers $m$ and $q$ as well as the sets $ E_t $, $ F_r $ may vary from monomial to monomial, i.e., they are functions of $k$ and $\alpha$. 
Furthermore, for each fixed $k$ and $\alpha$,  the lower indices and the upper index sets satisfy 
\begin{enumerate}
\item[(a)] 
$ a_t,b_t \2\in\2 \c{u_r,v_r}_{r=1}^p \!\cup \wh{B}(k)  $, and $ w_s \in \c{a_t,b_t}_{t=1}^m $; 
\item[(b)]
$ E_t \subseteq \wh{B}(k) \cup Q_{t'} $, for some $ 1 \leq t' \leq 2\mu $, 
and 
$ F_r \subseteq \wh{B}(k) \cup Q_{r'} $, for some $ 1 \leq r' \leq 2\mu $;
\item[(c)]
If $ a_t \in B_{s_i} $ and $ b_t \in B_{s_j} $, with $ 1\leq i,j\leq k $, then $ i \neq j $;
\item[(d)] For each $ s=1,\dots,2\mu $ there are two unique labels $ 1 \leq t'(s), t''(s) \leq m $, such that $a_{t'(s)} = u_s $, $ b_{t'(s)} \notin \c{\1v_r}_{r\neq s}$, and $ a_{t''(s)} \notin \c{\1u_r}_{r\neq s}$, $ b_{t''(s)} = v_s $ hold, respectively. 
\end{enumerate}
We will call the right hand side of \eqref{level-k expansion} the \emph{level-$ k $ expansion} in the following and we will define it by a  recursion on $ k $.

The level-$ 0 $ expansion is determined by the formula \eqref{def of Gamma at level k=0}:
\bels{def of level-0 expansion}{
\Gamma^{(0)}_1  := \Gamma_{\!\xi}|_{\mathbb{B}(\xi,\sigma)}
\,,\qquad M_0 := 1
\,.
}
This monomial clearly satisfies (a)--(d), with $ m = 2\mu$, $ q = 0 $,  $ E_t = Q_t $, and $ t'(s) = t''(s) = s $. 
The final goal, the representation \eqref{monomial exp for Gamma on BB}, is the last level-$ \wh{\ell}$ expansion, i.e.,
\bels{FA:final expansion identified}{
\Gamma_{\!\xi,\sigma,\alpha} :=\, \Gamma^{(\wh{\ell})}_\alpha
\,,\qquad
\alpha=1,2,\dots,M_{\wh{\ell}} \2=:\1 M(\xi,\sigma) 
\,.
}

Now we show how the level-$ k$  expansion is obtained given the level-$ (k-1)$ expansion.  
In order to do that, first we list the elements of each $B_{s_k}$  as  
$ 
\c{x_{ka}: 1 \leq a \leq \abs{B_{s_k}}} = B_{s_k}
$,
and we define 
\[
B_{k1} := \emptyset\,,\qquad
B_{ka} := \c{x_{kb}: 1 \leq b \leq a-1}
\,,
\qquad a =2,\dots,\abs{B_{s_k}}\,,
\quad
k =1,2,\dots,\wh{\ell}
\,,
\]
which is a one-by-one exhaustion of $B_{s_k}$;  namely  $B_{k1}\subseteq B_{k2} \subseteq \ldots \subseteq B_{k, \abs{B_{s_k}}}\subseteq B_{s_k}$. 
Note that $B_{k,a+1} = B_{ka}\cup\{ x_{ka}\}$. 

We now  consider a generic level-$ (k-1) $ monomial $ \Gamma^{(k-1)}_\alpha$, which is  of the form \eqref{monomial in simple standard form} and satisfies (a)--(d). 
Each monomial  $ \Gamma^{(k-1)}_\alpha$ will give rise  several level-$k$ monomials that are constructed independently for different $\alpha$'s as follows. 
Expanding each of the $ m $ factors in the first product of \eqref{monomial in simple standard form} using  the standard resolvent expansion identity
\begin{subequations}
\label{resolvent expansions for FA}
\bels{resolvent expansions for FA:regular}{
G^{E}_{ab} 
\,=\; 
G^{E\1\cup B_{s_k}}_{ab}
\,+\; 
\sum_{
a'\1=\11
}
^{\abs{B_{s_k}}}
\bbm{1}(\1x_{ka'}\notin E\2)\,
\frac{
G^{E\1\cup B_{ka'}}_{a\1x_{ka'}}G^{E\1\cup B_{ka'}}_{x_{ka'}b}
}{
G^{E\1\cup B_{ka'}}_{x_{ka'}x_{ka'}}
}
\,,
}
and each of  the $ q $ factors in the second product of \eqref{monomial in simple standard form} using 
\bels{resolvent expansions for FA:reciprocal}{
\frac{\!1}{G^{F}_{ww}} 
\,&=\, 
\frac{1}{\,G^{F \1\cup B_{s_k}}_{ww}\!}\; 
-
\sum_{a \1=\11}^{\abs{B_{s_k}}} 
\bbm{1}(\1x_{ka}\notin F\2)\,
\frac{
G^{F\1\cup B_{ka}}_{w\1x_{ka}}G^{F\1\cup B_{ka}}_{x_{sa}w}
}{
G^{F\1\cup B_{ka}}_{ww} G^{F\1\cup B_{k,a+1}}_{ww} G^{F\1\cup B_{ka}}_{x_{ka}x_{ka}}
}
\,,  
}
\end{subequations}
yields a product of sums of resolvent entries and their reciprocals.

Inserting these formulas into \eqref{monomial in simple standard form} and
expressing the resulting product as a single sum yields the representation
\bels{partial ordering sum: one step}{
\Gamma^{(k-1)}_\alpha
=\2 {\textstyle \sum_{\1\beta \1\in\1  \cal{A}_\alpha\msp{-1}(k)}}\2 \Gamma^{(k)}_\beta
\,,
}
where $ \cal{A}_\alpha(k) $  is some finite subset of integers and $\beta$ simply labels the resulting monomials in an arbitrary way.
From the resolvent identities \eqref{resolvent expansions for FA} it is easy to see that the monomials $ \Gamma^{(k)}_\beta $ inherit the properties (a)--(d) from the level-$ (k-1) $ monomials.
In particular, summing over $ \alpha =1,\dots,M_{k-1} $ in \eqref{partial ordering sum: one step} yields the level-$ k $ monomial expansion \eqref{level-k expansion}, with $ M_k := \sum_\alpha \abs{\1\cal{A}_\alpha(k)} $.
We will assume w.l.o.g. that the sets $ \cal{A}_\alpha(k) $, $ 1 \leq \alpha \leq M_{k-1}$, form a partition of the first $ M_k $ integers.  

This procedure defines the monomial representation recursively. 
Since  $  \Gamma^{(k)}_{\!\alpha}$ is a {\it function} of the $(\bs{u},\bs{v})$ indices, strictly speaking we should record which lower indices in the generic form \eqref{monomial in simple standard form} are considered \emph{independent variables}. 
Initially, at level $k=0$, all indices are variables, see \eqref{def of Gamma at level k=0}. Later, the expansion formulas \eqref{resolvent expansions for FA} bring in new lower indices, denoted generically by $x_{ka} $ from the set $ \cup_{s\in  \wh{L}} B_s $ which is disjoint from the range of the components $ u_r,v_r $ of the variables $ (\bs{u},\bs{v}) $ as $ \mathbb{B}(\xi,\sigma) $ is a subset of $ (\bb{X}\1\backslash\!\cup_{s\in  \wh{L}} B_s)^{2\mu} $.
However, the structure of \eqref{resolvent expansions for FA} clearly shows at which location the "old" $a, b$ indices from the left hand side of these formulas appear in the "new" formulas on the right hand side. 
Now the simple rule is that if any of these indices $a, b$ were variables on the left hand side, they are considered variables on the right hand side as well. In this way the concept of independent variables is naturally inherited along the recursion. With this simple rule we avoid the cumbersome notation of explicitly 
indicating which indices are variables in the formulas.

We note that the monomials of the final expansion \eqref{FA:final expansion identified} can be written in the form \eqref{typical monomial V2}. Indeed, the second products in \eqref{typical monomial V2} and \eqref{monomial in simple standard form} are the same, while the first product of \eqref{monomial in simple standard form} is split into the three other products in  \eqref{typical monomial V2} using (d).
Properties~1 and 2 in Lemma~\ref{lmm:Monomial representation} for the monomials in \eqref{FA:final expansion identified} follow easily from (a)--(d). 
Indeed,  (a) yields the first part of Property~1, while the second part of Property~1  follows from (c) and the basic property $ d(B_s,B_t) \ge N^{\eps_1} $ for  distinct lone labels $ s,t \in \wh{L} $.

For a given $\xi$, we define the family of subsets of $ \bb{X} $:
\[
\scr{E} 
\,:=\, 
\cB{
B_{1,a(1)} \cup B_{2,a(2)} \cup \cdots \cup B_{\wh{\ell},a(\wh{\ell})} \cup Q_r
: 
1 \leq a(k) \leq \abs{B_{s_k}}\,, \;1 \leq k\leq \wh{\ell}\,,\;1\leq r\leq 2\mu\;
}
\,.
\]
By construction (cf. \eqref{resolvent expansions for FA} and (b)) the upper index sets are members of this $ \xi $-dependent family.
Since $ \abs{\1Q_r},\abs{B_{s_k}} \leq N^{C_0\eps_1} $, for some $ C_0 \sim 1 $, we get $ \abs{\scr{E}} \lesssim_\mu \2N^{C_0\mu} $.
Property~2 follows directly from these observations.

Next we prove Property 3 of the monomials \eqref{FA:final expansion identified}.
To this end, we use the formula \eqref{partial ordering sum: one step} to define a partial ordering '$ < $' on the monomials by %
\bels{partial ordering}{
\Gamma^{(k-1)}_\alpha <\,\Gamma^{(k)}_\beta 
\quad\Longleftrightarrow\quad
\beta \in \cal{A}_\alpha(k)
\,.
}
It follows that for every $ \alpha =1,2,\dots, M = M_{\1\wh{\ell}} $, there exists a sequence $ (\alpha_k)_{k=1}^{\2\wh{\ell}-1} $, such that
\bels{FA:chain of monomials}{ 
\Gamma_{\!\xi}|_{\mathbb{B}(\xi,\sigma)}  =\; \Gamma_1^{(0)} \,<\; \Gamma^{(1)}_{\alpha_1} \,<\;\cdots\; <\; \Gamma^{(\1\wh{\ell}\1-1)}_{\!\alpha_{\wh{\ell}-1}} <\; \Gamma^{(\1\wh{\ell}\1)}_\alpha =\;\Gamma_{\!\xi,\sigma,\alpha} 
\,.
}

Let us fix an arbitrary label $ \alpha =1,\dots,M $ of the final expansion. 
Suppose that the $ k$-th monomial $ \Gamma^{(k)}_{\alpha_k} $, in the chain \eqref{FA:chain of monomials}, is of the form \eqref{monomial in simple standard form}, and define
\bels{defs of D_k and m_k}{
D_k \,:=\, 
\pB{\2\textstyle \bigcap_{\1t\1=\11}^{\1m} E_t\2} \,\cap\, \pB{\textstyle \2\bigcap_{\1r\1=\11}^{\1q} F_r\2}
,\qquad
m_k \,:=\, m
\,.
}
Here, $ D_k $ is the largest set $ A \subseteq \bb{X} $, such that $ \Gamma^{(k)}_{\alpha_k} $ depends only on the matrix elements of $\bs{H}^{(A)} $. 

Since both the upper index sets and the total number of resolvent elements of the form $ G^{(A)}_{ab}$ are both larger (or equal)  on the right hand side than on the left hand sides of the identities \eqref{resolvent expansions for FA}, and the added indices on the right hand side are from $ B_{s_k}$, we have
\[
D_{k-1} \subseteq D_k\,,\qquad
D_k\backslash D_{k-1} \subseteq B_{s_k}\,,
\qquad\text{and}\qquad
m_k \;\ge\; m_{k-1}
\,.
\]
We claim that
\bels{implications for the property 3}{
B_{s_k} \nsubseteq D_{\wh{\ell}}
\;&\implies\;
B_{s_k} \nsubseteq D_k
\;\implies\;
m_k \,\ge\,m_{k-1}\msp{-2}+1
\,. 
}
The first implication follows from the monotonicity of $ D_k$'s.
In order to get the second implication, suppose that $ \Gamma^{(k-1)}_{\alpha_{k-1}} $ equals \eqref{monomial in simple standard form}. 
Since $ D_k $ does not contain $ B_{s_k} $ the monomial $ \Gamma^{(k)}_{\alpha_k} $ can not be of the form 
 \eqref{monomial in simple standard form}, with the upper index sets $ E_t $ and $ F_t $ replaced with $ E_t \cup B_{s_k} $ and $ F_t \cup B_{s_k} $, respectively. 
The formulas \eqref{resolvent expansions for FA} hence show that $ \Gamma^{(k)}_{\alpha_k} $ contains at least one more resolvent entry of the form $ G^{(A)}_{ab} $ than $ \Gamma^{(k-1)}_{\alpha_{k-1}} $, and thus $ m_k \ge m_{k-1}+1 $.

Property 3 follows from \eqref{implications for the property 3}. Indeed, suppose that there are no isolated label $ s $ such that $ B_s \subseteq D_{\wh{\ell}} $. Then applying \eqref{implications for the property 3} for each $ k =1,\dots,\wh{\ell} $, yields $
m_{\1\wh{\ell}} 
\,\ge\,
m_0 +\2 \wh{\ell} 
\,.
$
Since $ m_0 = p $, using the notations from \eqref{typical monomial V2} we have
\[
m_{\1\wh{\ell}} \;=\, n + \abs{R^{(1)}} + 2\1\abs{R^{(2)}}
\,, 
\]
by Property (c) of the monomials. This completes the proof of Property 3.

Now only the bound \eqref{FA:bound on M} on the number of monomials $ M = M_{\wh{\ell}} $ remains to be proven, which is a simple counting. 
Let $p_k$ be the  largest number  of factors among the monomials  at the level-$ k $ expansion, i.e., writing a monomial $ \Gamma^{(k)}_\alpha = \Gamma^{(k)}_{\!\xi,\sigma,\alpha} $ in the form \eqref{typical monomial V2} we have
\[ 
p_k \,:= \msp{-10}\max_{\msp{10}1 \1\leq\1 \alpha \1\leq\1 M_k}\, \bigl(\2n(\alpha)+ \abs{R^{(1)}(\alpha)} + 2\1\abs{R^{(2)}(\alpha)} +q(\alpha)\,\bigr)
\,,
\]
where $ M_k = M_k(\xi,\sigma) $, $ n(\alpha) = n(\xi,\sigma,\alpha) $, $ R^{(1)}(\alpha) = R^{(1)}(\xi,\sigma,\alpha) $, etc.
Let us set $ b_\ast := 1+ \max_{x,y}\,\abs{B_{N^{\eps_1}}(x,y)} $.
Each of the factors in every monomial at the level $ k-1 $  is turned into a sum over monomials by the resolvent identities \eqref{resolvent expansions for FA}.
Since each such monomial contains at most five resolvent entries (cf. the last terms in \eqref{resolvent expansions for FA:reciprocal}), we obtain the first of the following two bounds: %
\bels{recursion for p_k and M_k}{
p_k \leq 5\1p_{k-1}
\qquad\text{and}\qquad
M_k \leq M_{k-1} \2b_\ast^{\,p_{k-1}}%
\,.
}
For the second bound we recall that each of the at most $ p_{k-1} $ factors in 
every level-$ (k-1)$ monomial 
is expanded by the resolvent identities \eqref{resolvent expansions for FA} into a sum of at most $ b_\ast $ terms. 
The product of these sums yields single sum of at most  $ b_\ast^{\2p_{k-1}} $ terms. 
From \eqref{def of level-0 expansion} and \eqref{def of Gamma at level k=0} we get: $ M_0 := 1 $,  $ p_0 = 2\mu $. 
Since $ k \leq \wh{\ell} \leq 2\mu $, we have $ \max_k p_k \leq  2\mu \,25^\mu $. 
Plugging this into the second bound of \eqref{recursion for p_k and M_k} yields $ M_k \leq (\1(b_\ast)^{2\mu \,25^\mu})^{2\mu} $. %
This proves \eqref{FA:bound on M} since $ b_\ast \leq N^{C\eps_1} $ by \eqref{polynomial ball growth}.
Finally, we obtain the bound on the number of factors in \eqref{typical monomial V2} using $ n+q \leq p_{\wh{\ell}} \lesssim_\mu 1 $. 
\end{Proof}

\section{Bulk universality and rigidity}
\label{sec:Bulk universality and rigidity}

 In this section we show how to use the strong local law, Theorem~\ref{thr:Local law for correlated random matrices}, to obtain the remaining results of Section \ref{subsec:Correlated random matrices} on random matrices with correlated entries.
\subsection{Rigidity}
\label{subsec:Rigidity}

\begin{proposition}[Local law away from {$[\kappa_-,\kappa_+]$}]
\label{prp:Local law away}
Let $\bs{G}$ be the resolvent of a random matrix $\bs{H}$ of the form \eqref{H structure} that satisfies {\bf B1}-{\bf B4}. Let $\kappa_-, \kappa_+$ be the endpoints of the convex hull of $\supp \rho$ as in \eqref{definition kappa-+}. 
For all $\delta,\eps>0$ and $\nu \in \N$ there exists a positive constant $C$ such that away from $[\kappa_-,\kappa_+]$,
\bels{entrywise local law away}{
\msp{-5}
\P\biggl[\,
\exists\, \zeta \in \Cp\, \text{ s.t. }\delta \le \dist(\1\zeta\1,[\kappa_-,\kappa_+])\le \frac{1}{\delta}\,,\, \max_{x,y=1}^N\abs{\1G_{xy}(\1\zeta\1)-m_{xy}(\1\zeta\1)}\ge \frac{N^\eps}{\sqrt{N }}
\,\biggr]
\le
\frac{C}{\2N^\nu\!}
\,.
}
The normalized trace converges with the improved rate 
\bels{trace local law away}{
\msp{-8}\P\biggl[\,
\exists\, \zeta \in \Cp \, \text{ s.t. }\delta \le \dist(\1\zeta\1,[\kappa_-,\kappa_+])\le \frac{1}{\delta}\,,\; \absB{\text{\small$\frac{1}{N}$}\msp{-1}\tr\bs{G}(\1\zeta\1)-\text{\small$\frac{1}{N}$}\msp{-1}\tr\bs{M}(\1\zeta\1)}\ge \frac{N^\eps}{N }
\,\biggr]
\le
\frac{C}{\2N^\nu\!}
\,.
}
The constant $C$ depends only on the model parameters $\scr{K}$ in addition to $\delta$, $\eps$ and $\nu$. 
\end{proposition}

\begin{remark}
\label{rmk:Local Law away from supp rho}
Theorem~\ref{thr:Local law for correlated random matrices} and Proposition~\ref{prp:Local law away} provide a local law with optimal convergence rate $\frac{1}{N \im \zeta}$ inside the bulk of the spectrum and convergence rate $\frac{1}{N}$ away from the convex hull of $\supp \rho$, respectively.
 In order to prove a local law inside spectral gaps and at the edges of the self-consistent spectrum, additional assumptions on $\bs{H}$ are needed to exclude a naturally appearing instability that may be caused by exceptional rows and columns of $\bs{H}$ and the outlying eigenvalues they create. This instability is already present in the case of independent entries as explained in Section~11.2 of \cite{AEK1}. 
\end{remark}

\begin{remark} The local law in Proposition~\ref{prp:Local law away} extends beyond the regime of bounded spectral parameters $\zeta$. The upper bound $\frac{1}{\delta}$ on the distance of $\zeta$ from $[\kappa_-,\kappa_+]$ can be dropped in both \eqref{entrywise local law away} and \eqref{trace local law away}. Furthermore, as was done e.g. for Wigner-type matrices in \cite{AEK2}, by following the $\abs{\zeta}$-dependence along the proof the estimates on the difference $\bs{G}-\bs{M}$ in  \eqref{entrywise local law away} and \eqref{trace local law away} can be improved to 
$\frac{N^\eps}{(1+\abs{\zeta}^2)\sqrt{N} }$ 
and 
$\frac{N^\eps}{(1+\abs{\zeta}^2){N}}$, 
respectively. Since this extra complication only extends the local law to a regime far outside the spectrum of $\bs{H}$ (cf. Lemma~\ref{lmm:No eigenvalues away} below) we refrain from carrying out this analysis.   
\end{remark}

\begin{Proof}[Proof of Proposition~\ref{prp:Local law away}] The proof has three steps. In the first step we will establish a weaker version of Proposition~\ref{prp:Local law away} where instead of the bound $\Lambda \prec N^{-1/2}$ we will only show $\Lambda \prec 
N^{-1/2} \!+ (N\im \zeta)^{-1}$.
Then we will use this version in the second step to prove that there are no eigenvalues outside a small neighborhood of $[\kappa_-,\kappa_+]$. Finally, in the third step we will show \eqref{entrywise local law away} and \eqref{trace local law away}.  
\medskip
\\
\emph{Step 1}:
The proof of this step follows the same strategy as the proof of Theorem~\ref{thr:Local law for correlated random matrices}. Only instead of using Lemma~\ref{lmm:Smallness of error matrix} to estimate the error matrix $\bs{D}$ we will use Lemma~\ref{lmm:Smallness of error matrix away from support}. In analogy to the proof of \eqref{entrywise local law in bulk} we begin by showing  the  entrywise bound
\bels{Step 1 Lambda bound away}{
\Lambda(\1\zeta\1)\,\prec\, \frac{1}{\sqrt{N}}+\frac{1}{N\im\zeta}\,,\qquad \zeta \in \Cp\,,\; \delta \le \dist(\1\zeta\1,[\kappa_-,\kappa_+])\le \frac{1}{\delta}\,,\;\im \zeta \ge N^{-1+\eps}.
}

In fact, following the same line of reasoning that was used to prove \eqref{Gap in the values of G - M max-norm}, but using \eqref{error bound away from convex hull of support} instead of \eqref{bound on error D in the bulk} to estimate $\norm{\bs{D}}_{\rm{max}}$ we see that 
\bels{Weak Lambda bound away}{
\Lambda(\1\zeta\1)\bbm{1}(\Lambda(\1\zeta\1)\le N^{-\eps/2})\,\prec\, \frac{1}{\sqrt{N }}+\pB{\frac{\Lambda(\1\zeta\1)}{N \im \zeta}}^{1/2}
\,\le\, \frac{1}{\sqrt{N }}+ \frac{N^{\eps}}{N \im \zeta} + 4\2 N^{-\eps}\Lambda(\1\zeta\1)\,,
}
 for any $\eps>0$. 
The last term on the right hand side can be absorbed into the left hand side and since $\eps$ was arbitrary \eqref{Weak Lambda bound away} yields
\bels{Weak Lambda bound away with indicator}{
\Lambda(\1\zeta\1)\bbm{1}(\Lambda(\1\zeta\1)\le N^{-\eps/2})\,\prec\,\frac{1}{\sqrt{N}}+\frac{1}{N\im\zeta}\,.
}
This inequality establishes a gap in the possible values that $\Lambda$ can take, provided $\eps< 1/2$ because $N^{-\eps} \ge N^{-1/2} \!+ (N\im \zeta)^{-1}$.
Exactly as we argued for \eqref{Gap in the values of G - M max-norm} we can get rid of the indicator function in \eqref{Weak Lambda bound away with indicator} by using a continuity argument together with a union bound in order to obtain \eqref{Step 1 Lambda bound away}. 

As in the proof of Theorem~\ref{thr:Local law for correlated random matrices} we now use the fluctuation averaging to get an improved convergence rate for the normalized trace,
\bels{Weak averaged bound away}{
\absB{\text{\small$\frac{1}{N}$}\msp{-1}\tr(\bs{G}(\1\zeta\1)-\bs{M}(\1\zeta\1))}
\,\prec\, \frac{1}{{N}}+\frac{1}{(N\im\zeta)^2}\,.
}
for all $\zeta \in \Cp$ with $ \delta \le \dist(\1\zeta\1,[\kappa_-,\kappa_+])\le \frac{1}{\delta}$ and $\im \zeta \ge N^{-1+\eps}$.
Indeed, \eqref{Weak averaged bound away} is an immediate consequence of \eqref{Step 1 Lambda bound away} and the  fluctuation averaging Proposition~\ref{prp:fluctuation averaging}.  
\medskip
\\
\emph{Step 2}: In this step we use \eqref{Weak averaged bound away} to prove the following lemma. 
\begin{lemma}[No eigenvalues away from {$[\kappa_-,\kappa_+]$}] 
\label{lmm:No eigenvalues away}
For any $\delta,\nu>0$ we have
\bels{No eigenvalues away}{
\P\sB{\2
\spec(\bs{H}) \cap (\2\R\setminus\1[\2\kappa_-\!-\delta\1,\kappa_+\!+\delta\2]\2) = \emptyset
 \2}
\,\ge\, 1- C N^{-\nu}
,
}
for a positive constant $C$, depending only on the model parameters $\scr{K}$ in addition to $\delta$ and $\nu$.
\end{lemma}

In order to show \eqref{No eigenvalues away} fix $\tau \in [-\delta^{-1}\!,\2\kappa_-\!-\delta\2] \cup [\1\kappa_+\!+\delta\1,\2 \delta^{-1}]$ and $ \eta \in[N^{-1+\eps}\!,1\2]$, and let $\sett{\lambda_i}_{i=1}^N $ be the eigenvalues of $ \bs{H}$. 
Employing \eqref{Weak averaged bound away} we get
\bels{single term in sum bound for no ev}{
\frac{\eta}{(\lambda_i-\tau)^2+\eta^2}
\,\le\,  
\im \tr \bs{G}(\tau +\ii\1 \eta)
\,\prec\, 
\im \tr \bs{M}(\tau +\ii\1 \eta) + 1 + \frac{1}{N\eta^2}
\,\lesssim_{\1\delta}\, 
N \eta + 1 + \frac{1}{N\eta^{\12}\!}
\,.
}
Here, we used in the last inequality that $\frac{1}{N} \tr \bs{M}$ is the Stieltjes transform of the self-consistent density of states $\rho$ with $\supp \rho \subseteq [\kappa_-, \kappa_+]$. Since the left hand side of \eqref{single term in sum bound for no ev} is a Lipschitz continuous function in $\tau$ with Lipschitz constant bounded by $N$ we can use a union bound to establish \eqref{single term in sum bound for no ev} first on a fine grid of $\tau$-values and then uniformly for all $\tau$ and for the choice $\eta=N^{-2/3}$,
\[
\sup_{\tau}\frac{1}{N^{4/3}(\lambda_i-\tau)^2+1}\,\prec\, \frac{1}{N^{1/3}}\,.
\]
In particular, the eigenvalue $\lambda_i$ cannot be at position $\tau$ with very high probability, i.e.
\bels{No eigenvalues on compact interval away}{
\P\1\sb{\,\exists \,i \text{ s.t. }
\delta \le \dist(\lambda_i, [\kappa_-, \kappa_+])\le \delta^{-1}  
 }
\,\le\,  C(\delta,\nu) N^{-\nu}.
}

Now we exclude that there are eigenvalues  far away from the self-consistent spectrum by using a continuity argument. 
Let $\wt{\bs{W}}$ be a standard GUE matrix with $\E\abs{\wt{w}_{xy}}^2=\frac{1}{N}$, $(\lambda_i^{(\alpha)})_i$ the eigenvalues of $\bs{H}^{(\alpha)}:=\alpha\2\bs{H}+(1-\alpha)\wt{\bs{W}}$ for $\alpha \in [0,1]$ and $\kappa:=\sup_{\alpha}\max\{\abs{\kappa^{(\alpha)}_+}, \abs{\kappa^{(\alpha)}_-}\}$, where $\kappa_\pm^{(\alpha)}$ are defined as in \eqref{definition kappa-+} for the matrix $\bs{H}^{(\alpha)}$. In particular, $\kappa_\pm^{(0)}=\pm 2$. 
Since the constant $C(\delta,\nu)$ in \eqref{No eigenvalues on compact interval away} is uniform for all random matrices with the same model parameters $\scr{K}$, we see that 
\[
\sup_{\alpha \1\in\1 [0,1]}\,
\P\sb{\,\exists\, i \text{ s.t. }\abs{\lambda_i^{(\alpha)}} \in [\kappa+\delta, \delta^{-1}]\,
}\,\le\,  C(\delta,\nu) N^{-\nu}
\]
The eigenvalues $\lambda_i^{(\alpha)}$ are Lipschitz continuous in $\alpha$. In fact, $\abs{\partial_\alpha \lambda_i^{(\alpha)}}\le \norm{\bs{H}-\wt{\bs{W}}}\prec\sqrt{N}$. Here, the simple bound on $\norm{\bs{H}-\wt{\bs{W}}}$ follows from 
$ \E\2\norm{\bs{H}-\wt{\bs{W}}}^{2\mu} = \E\2\s{\tr(\bs{H}-\wt{\bs{W}})^2}^\mu\,\le\,C(\mu)N^\mu $,  
for some positive constant $C(\mu)$, depending on $\mu$, the upper bound $\ul{\kappa}_1$ from \eqref{bounded moments} on the moments, the sequence $\ul{\kappa}_2$ from \eqref{decay of expectation} and $P$ from \eqref{polynomial ball growth}.  Thus we can use a union bound to establish
\bels{No eigenvalue can leave}{
\P\1\sb{\,\exists \2 \alpha, i\, \text{ s.t. } \abs{\lambda_i^{(\alpha)}}\in[\kappa+2\2\delta, \delta^{-1}-\delta]\,
}\,\le\,  C(\delta,\nu) N^{-\nu}.
}
Since for $\alpha=0$ all eigenvalues are in $[-\kappa-2\2\delta, \kappa+2\2\delta]$ with very high probability and with very high probability no eigenvalue can leave this interval by \eqref{No eigenvalue can leave}, we conclude that 
\[
\P\2\sb{\,\exists \, i \text{ s.t. }\abs{\lambda_i}\2\ge\2 \kappa+2\1\delta\,
}\,\le\,  C(\delta,\nu) N^{-\nu}.
\]
Together with \eqref{No eigenvalues on compact interval away} 
this finishes the proof of Lemma~\ref{lmm:No eigenvalues away}.
\medskip
\\
\emph{Step 3}: In this step we use \eqref{No eigenvalues away} to improve the bound on the error matrix $\bs{D}$ away from $[-\kappa_-, \kappa_+]$ and thus show \eqref{entrywise local law away} and \eqref{trace local law away} by following the same strategy that was used in Step~1 and in the proof of  Theorem~\ref{thr:Local law for correlated random matrices}. 

By Lemma~\ref{lmm:No eigenvalues away} there are with very high probability no eigenvalues in 
 $ \R\setminus\1 [\kappa_{\msp{-1}-}\msp{-7}-\msp{-2}\delta/2,\kappa_{\msp{-1}+}\msp{-5}+\msp{-2}\delta/2\1]\, $. 
Therefore, for any $B \subseteq \bb{X}$ also the submatrix $\bs{H}^B$ of $\bs{H}$ has no eigenvalues in this interval. 
In particular, for any $x \in \bb{X}\setminus B$ we have 
\bels{GBxx comparable to eta}{
\im G_{xx}^B(\1\zeta\1)\,\sim_\delta\, \im \zeta\,, \qquad \zeta \in \Cp\,,\; \delta \le \dist(\1\zeta\1,[\kappa_-,\kappa_+])\le \delta^{-1}\,,
}
in a high probability event. 
As in the proof of Lemma~\ref{lmm:Smallness of error matrix} we bound the entries of the error matrix $\bs{D}$ by estimating the right hand sides of the equations \eqref{estimate for D1} to \eqref{estimate for D5} further. But now we use \eqref{GBxx comparable to eta}, so that $\im \zeta$ in the denominators cancel and we end up with
\bels{optimal D bound away}{
\qquad \norm{\1\bs{D}(\1\zeta\1)}_{\rm{max}}\2\bbm{1}(\1\Lambda(\1\zeta\1)\le N^{-\eps})\,\prec\, N^{-1/2}
\,,
\quad\text{whenever}\quad \delta \le \dist(\1\zeta\1,[\kappa_-,\kappa_+])\le \delta^{-1}.
}
Following the strategy of proof from Step~1 we see that \eqref{optimal D bound away} implies \eqref{entrywise local law away} and \eqref{trace local law away}. This finishes the proof of Proposition~\ref{prp:Local law away}.
\end{Proof}

\begin{Proof}[Proof of Corollary~\ref{crl:Rigidity}]
The proof follows a standard argument that establishes rigidity from the local law, which we present here for the convenience of the reader. The argument uses a Cauchy-integral formula that was also applied in the construction of the Helffer-Sj\"ostrand functional calculus (cf. \cite{DFunc}) and it already appeared in different variants in \cite{EYYber}, \cite{ERSY} and \cite{EYY}.

Let $\tau \in \R$ such that $\rho(\tau)\ge \delta$ for some $\delta>0$. 
 We will now apply Lemma~5.1 of \cite{AEK2} which shows how to estimate the difference between two measures in terms of the difference of their Stieltjes transforms. With the same notation that was used in the statement of that lemma we make the choices 
\[
\nu_1(\dd \sigma)\,:=\, \rho(\sigma)\1\dd \sigma\,,\qquad \nu_2(\dd \sigma)\,:=\,
 \frac{1}{N}\sum_i \delta_{\lambda_i}\msp{-2}(\dd \sigma)\,,
\]
and $ \tau_1 := \kappa_- -\wt \delta$, $ \tau_2 := \tau $, $ \eta_1 := N^{-1/2} $, $ \eta_2 := N^{-1+\wt \eps}$, $ \eps :=  1$, for some fixed $\wt \delta, \wt \eps >0$. We estimate the error terms $J_1$, $J_2$ and $J_3$ from Lemma~5.1 of \cite{AEK2} by using \eqref{trace local law away} and \eqref{trace local law in bulk}.    In this way we find
\[
\absB{\,N
{\textstyle \,\int_{[\1\kappa_--\wt{\delta}\1,\1\tau\1]\2}}
\rho(\sigma) \2\dd \sigma\,-\, \abs{\2\spec(\bs{H})\!\cap\! [\kappa_-\!-\wt \delta,\tau\1]\2}\,}
\,\prec\, N^{\wt \eps}.
\]
Since $\wt \eps$ was arbitrary and there are no eigenvalues of $\bs{H}$ to the left of $\kappa_-\!-\wt \delta$ (cf. Lemma~\ref{lmm:No eigenvalues away}), we infer 
\bels{difference integrated DOS}{
\absB{\,N
{\textstyle \,\int_{[-\infty,\tau\1]\2}}
\rho(\sigma) \2\dd \sigma \,-\, \abs{\2\spec(\bs{H})\!\cap\! (-\infty,\tau\1]\2}\,}
\,\prec\, 
1\,,
}
for any $\tau \in \R$ with $\rho(\tau)\ge\delta$. Combining \eqref{difference integrated DOS} with the definition \eqref{energy index} of $i(\tau)$  yields the bound  $ \absb{ \int_\tau^{\lambda_{i(\tau)}} \!\rho(\sigma)\2\dd \sigma\,} \prec N^{-1} $.
This in turn implies \eqref{bulk rigidity} and Corollary~\ref{crl:Rigidity} is proven.
\end{Proof}

\subsection{Bulk universality}
\label{subsec:Bulk universality}
Given the local law (Theorem~\ref{thr:Local law for correlated random matrices}), the proof of bulk universality  (Corollaries~\ref{crl:Bulk universality}
and \ref{crl:Bulk Gap Universality})  follows standard arguments based upon the three step strategy explained in the introduction. We will only sketch the main differences due to the correlations. 
We start by introducing an Ornstein--Uhlenbeck (OU) process on random matrices $ \bs{H}_t $ that conserves the first two mixed moments of the matrix entries
\bels{SDE for matrix flow}{
\dd \bs{H}_{\1t} \,=\, -\frac{1}{2}(\2\bs{H}_{\1t}-\bs{A}\1)\2\dd t +  \Sigma^{1/2}[\1\dd \bs{B}_{\1t}]\,,
\qquad \bs{H}_0\,=\, \bs{H}\,,
} 
where the covariance operator $\Sigma: \C^{N \times N} \to \C^{N \times N}$ is given as
\[
\Sigma[\bs{R}]\,:=\, \E\,\scalar{\bs{W}}{\bs{R}}\2\bs{W}\,,
\]
and $\bs{B}_t$ is matrix of standard real (complex) independent Brownian motions with the appropriate symmetry $\bs{B}_t^*=\bs{B}_t$ for $\beta=1$ ($\beta =2$) whose distribution is invariant under the orthogonal (unitary) symmetry group.
We remark that a large Gaussian component, as created by the flow \eqref{SDE for matrix flow}, was first used in \cite{Johansson2000} to prove universality for the Hermitian symmetry class. 
  
Along the flow the matrix $ \bs{H}_t = \bs{A} +\frac{1}{\sqrt{N}} \bs{W}_t $ satisfies the condition {\bf B3} on the dependence of the matrix entries uniformly in $ t $. 
In particular, since $ \Sigma $ determines the operator $ \cal{S} $ we see that $ \bs{H}_t $ is associated to the same MDE as the original matrix $ \bs{H} $.
Also the condition {\bf B4} and {\bf B5} can be stated in terms of $ \Sigma$, and are hence both  conserved along the flow. 

For the following arguments we write $\bs{W}_t$ as a vector containing all degrees of freedom originating from the real and imaginary parts of the entries of $\bs{W}_t$. This vector has $N(N+1)/2$ real entries for $\beta=1$ and $N^2$ real entries  for $\beta=2$.
We partition $\bb{X}^2=\bb{I}_{\le}\dot{\cup}\bb{I}_>$ into its upper, $\bb{I}_{\le}:=\{(x,y): x \le y\}$, and lower, $\bb{I}_{>}=\{(x,y): x > y\}$, triangular part.
Then we identify
\[
\frac{1}{\sqrt{N}}\bs{W}_t\,=\, 
\begin{cases}
(w_t(\alpha))_{\alpha \in \bb{I}_{\le}}& \text{ if }\; \beta =1\,,
\\
(w_t(\alpha))_{\alpha \in \bb{X}^2}& \text{ if }\; \beta =2\,,
\end{cases}
\]
where $w_t((x,y)):=\frac{1}{\sqrt{N}}\2w_{xy}$ for $\beta =1$ and 
\[
w_t((x,y))\,:=\,
\begin{cases}
\frac{1}{\sqrt{N}}\re w_{xy}& \text{ for }\;  (x,y) \in \bb{I}_{\le}\,,
\\
\frac{1}{\sqrt{N}}\im w_{xy}& \text{ for }\; (x,y) \in \bb{I}_>\,,
\end{cases}
\]
for $\beta =2$. In terms of the vector $w_t$ the flow \eqref{SDE for matrix flow} takes the form
\bels{SDE flow}{
\dd w_t \,=\, -\1\frac{w_t}{2}
\2+\2 \Sigma^{1/2}\2\dd b_t\,,
}
where $ b_t=(b_t(\alpha))_\alpha$ is a vector of independent standard Brownian motions, and $ \Sigma^{1/2} $ is the square-root of the covariance matrix corresponding to $ \bs{H}=\bs{H}_0$:
\[
\Sigma(\alpha,\beta) \,:=\, \E\,w_0(\alpha)\1w_0(\beta) \,.
\]

Recall the notation $ B_\tau(x) = \c{y \in \bb{X} : d(x,y) \leq \tau } $ for any $x\in \bb{X}$,  and set
\begin{equation*}
\begin{split}
\cal{B}_k((x,y)) \,&:=\, (B_{k\1N^{\eps}}(x)\times B_{k\1N^{\eps}}(y))\cup (B_{k\1N^{\eps}}(y)\times B_{k\1N^{\eps}}(x))
\,,
\qquad
k=1,2
\,.
\end{split}
\end{equation*}
Using \eqref{polynomial ball growth} and {\bf B3} we see that for any $ \alpha $ 
\bels{size of mcl-B_2 and Sigma outside of it}{
\abs{\cal{B}_2(\alpha)} \,\leq\, N^{C\eps}
\qquad\text{and}\qquad
\abs{\Sigma(\alpha,\gamma)} &\leq C(\eps,\nu)N^{-\nu}\,,\qquad \gamma \notin \cal{B}_1(\alpha)
\,,
}
respectively.
For any fixed $\alpha  $, we denote by $ w^\alpha $ the vector obtained by removing all the entries of $ w $ which may become strongly dependent on the component $ w(\alpha) $ along the flow \eqref{SDE flow}, i.e., we define  
\bels{def of w^alpha}{
w^\alpha(\gamma) := w(\gamma)\1 \bbm{1}(\1\gamma \notin \cal{B}_2(\alpha)\1) 
\,.
} 
In the case that $ \bs{X}$ has independent entries it was proven in \cite{BY13EigenMoment} that the process \eqref{SDE flow} conserves the local eigenvalue  statistics of $ \bs{H} $ up to times $ t \ll N^{-1/2} $, provided bulk local law holds uniformly in $ t $ along the flow as well. 
We will now show that this insight extends for dependent random matrices as well.
The following result is a straightforward generalization of Lemma~A.1. from \cite{BY13EigenMoment}  to matrices with dependent entries.
A similar result was independently  given in \cite{Che2016}. 

\begin{lemma}[Continuity of the OU flow]
\label{lmm:Comparison flow}
For every $ \eps > 0 $, $\nu \in \N$ and smooth function $ f $ there is $ C(\eps,\nu) < \infty $, such that 
\bels{flow for matrix function}{
\abs{\2\E\,f(w_t)-\E\,f(w_0)\1} \,\leq\, C(\eps,\nu)\,\bigl(\,N^{1/2+\eps}\,\Xi\, + N^{-\nu}\1\wt{\Xi}\,\bigr)\,t
\,,
}
where
\vspace{-0.3cm}
\bels{def of Gamma and wti-Gamma}{
\Xi\,&:=\, 
\sup_{s\leq t}\max_{\alpha,\delta,\gamma}
\sup_{\theta\in[0,1]} \E\biggl[\,\Bigl(N^{1/2}\abs{w_s(\alpha)}+N^{3/2}\abs{\1w_s(\alpha)\1w_s(\delta)\1w_s(\gamma)}\Bigr)\,\absb{\2\partial^3_{\alpha\delta\gamma}f\bigl(w_s^{\alpha,\theta}\bigr)}\,
\biggr]
\\
\wt{\Xi} \,&:=\,
\sup_{\wt w} 
\max_{\alpha,\delta,\gamma}\, 
\Bigl(
\absb{\1\partial^2_{\alpha\delta}f(\wt w)}
\,+\,
(\11+\abs{w(\alpha)})\,
\absb{\partial^3_{\alpha\delta\gamma}f(\wt w)}\Bigr)
\,,
}
where $ w^{\alpha,\theta}_s \!:= w^{\alpha}_s + \theta\,(w_s-w^{\alpha}_s) $ for $ \theta \in [0,1] $, 
and $ \partial^k_{\alpha_1\cdots \alpha_k} =\frac{\partial}{\partial w(\alpha_1)}\cdots \frac{\partial}{\partial w(\alpha_k)} $.
\end{lemma}

\begin{Proof}
We will suppress the $ t $-dependence, i.e. we write $ w = w_t $, etc. Ito's formula yields
\bels{Ito for f(w)}{
\msp{-20}\dd f (w) 
\,=\,
\sum_\alpha 
\biggl(
-\frac{w(\alpha)}{2}\partial_\alpha f(w) + \frac{1}{2}\sum_{\delta} \Sigma(\alpha,\delta)\2\partial^2_{\alpha\delta} f(w)\biggr)
\1\dd t\,+\,\dd M
\,,
}
where $ \dd M = \dd M_t $ is a martingale term.
Taylor expansion around $ w = w^\alpha $ yields
\bea{
\partial_\alpha f(w) 
\,&=\, 
\partial_\alpha f(w^\alpha) 
\,+\msp{-10}
\sum_{\msp{10}\delta \in \cal{B}_2(\alpha)}\msp{-10}w(\delta)\2\partial^2_{\alpha\delta}f(w^\alpha)
\,+\msp{-16}
\sum_{\msp{16}\delta,\gamma \in \cal{B}_2(\alpha)}\msp{-16}
w(\delta)\1w(\gamma)\!\int_0^1\!(\11-\theta)\2\partial^3_{\alpha\delta\gamma}f(w^{\alpha,\theta})
\1\dd\theta
\\
\partial^2_{\alpha\delta} f(w) 
\,&=\,
\partial^2_{\alpha\delta} f(w^\alpha)
\,+\msp{-10}
\sum_{\msp{10}\gamma \in \cal{B}_2(\alpha)} \msp{-10}w(\gamma)\!
\int_0^1\!(\11-\theta)\,
\partial^3_{\alpha\delta\gamma} f(w^{\alpha,\theta})\1\dd\theta
\,.
}
By plugging these into \eqref{Ito for f(w)} and taking expectation, we obtain
\begin{subequations}
\label{derivative of EE f(w_t)}
\begin{align}
\frac{\dd}{\dd t}\E\,f(w)
\,&=\,
-\frac{1}{2}\sum_\alpha \E\,w(\alpha)\2\partial_\alpha f(w^\alpha)
\label{derivative of EE f(w_t): 1-st term} 
\\	
\label{derivative of EE f(w_t): 2-nd term} 
&\quad-\,\frac{1}{2}
\sum_{\alpha}\sum_{\delta \in \cal{B}_2(\alpha)} \E\Bigl[ \,\bigl(\1w(\alpha)\1w(\delta)-\Sigma(\alpha,\delta)\bigr)\2\partial^2_{\alpha\delta}f(w^\alpha)\,\Bigr] 
\\ 
\label{derivative of EE f(w_t): 3-rd term} 
&\quad
+\frac{1}{2}\sum_\alpha\sum_{\delta \notin \cal{B}_2(\alpha)} \Sigma(\alpha,\delta)\;\E\,\partial^{\12}_{\alpha\delta}f(w^\alpha)
\\
\label{derivative of EE f(w_t): 4-th term} 
&\quad-\;
\sum_\alpha \sum_{\delta,\gamma \in \cal{B}_2(\alpha)}
\int_0^1\! (\11-\theta)\,\E\Bigl[w(\alpha)\1w(\delta)\1w(\gamma)\2\partial^{\13}_{\alpha\delta\gamma}f(w^{\alpha,\theta})\Bigr]\,\dd \theta
\\
\label{derivative of EE f(w_t): 5-th term} 
&\quad
+
\frac{1}{2}\sum_{\alpha,\delta} \,\Sigma(\alpha,\delta) \msp{-8}\sum_{\gamma\in \cal{B}_2(\alpha)} 
\int_0^1\! (\11-\theta)\,
\E\Bigl[w(\gamma)\1\partial^{\13}_{\alpha\delta\gamma} f(w^{\alpha,\theta})\Bigr]\,\dd\theta
\,.
\end{align}
\end{subequations}
Now, we estimate the five terms on the right hand side of \eqref{derivative of EE f(w_t)} separately.

First, \eqref{derivative of EE f(w_t): 1-st term} is small since $ w(\alpha) $ is almost independent of $ w^\alpha $ by {\bf B3} and \eqref{def of w^alpha}: 
\[
\E\,w(\alpha)\2\partial_\alpha f(w^\alpha) 
\,=\, 
\E\,w(\alpha)\;\E\,\partial_\alpha f(w^\alpha) 
\,+\,\rm{Cov}(\1w(\alpha),\partial_\alpha f(w^\alpha)) 
=\,
\ord_{\eps,\nu}\bigl(\,\wt{\Xi}N^{-\nu}\,\bigr)
\,.
\]

In the term \eqref{derivative of EE f(w_t): 2-nd term},
if $ \delta \in \cal{B}_1(\alpha) $, then 
$ w(\alpha)\1w(\delta) $ is almost independent of $ w^\alpha $:
\[
\E\Bigl[ \bigl(\1w(\alpha)\1w(\delta)-\Sigma(\alpha,\delta)\1\bigr)\2\partial^2_{\alpha\delta}f(w^\alpha)\Bigr] 
\,=\,\rm{Cov}\bigl(\1w(\alpha)\1w(\delta)\2,\,\partial^2_{\alpha\delta}f(w^\alpha)\bigr)
\,=\,
\ord_{\eps,\nu}\bigl(\,\wt{\Xi}N^{-\nu}\2\bigr)
\,.
\]
If $ \delta \in \cal{B}_2(\alpha)\backslash \cal{B}_1(\alpha) $, then $ w(\alpha) $ is almost independent of $ (w(\delta),w^\alpha) $ and 
\bea{
\absB{\,\E\Bigl[ \bigl(\1w(\alpha)\1w(\delta)-\Sigma(\alpha,\delta)\bigr)\2\partial^2_{\alpha\delta}f(w^\alpha)\Bigr]\2} 
&\leq\,
\absb{\2\rm{Cov}\bigl(\1w(\alpha)\2,\,w(\delta)\partial^2_{\alpha\delta}f(w^\alpha)\bigr)} 
+\abs{\2\Sigma(\alpha,\delta)}\, \absb{\E\,\partial^2_{\alpha\delta}f(w^\alpha)} 
\\
&\leq\,
C(\eps,\nu)\, 
{\textstyle \sup_{\1w} \max_{\1\alpha,\1\delta,\gamma}}\,\bigl(\abs{\partial^2_{\alpha\delta}f(w)} +\abs{w(\alpha)}\abs{\partial^3_{\alpha\delta\gamma}f(w)}\bigr)\2N^{\1-\1\nu}
\,,
}
where we have used \eqref{size of mcl-B_2 and Sigma outside of it}. The last term containing derivatives is bounded by $ \wt{\Xi} $.

The term \eqref{derivative of EE f(w_t): 3-rd term} is negligible by $ \abs{\Sigma(\alpha,\delta)} \lesssim_{\eps,\nu} N^{-\nu} $ and $ \abs{\2\E\,\partial^2_{\alpha\delta}f(w^\alpha)\2} \leq  \wt{\Xi}  $.  
For \eqref{derivative of EE f(w_t): 4-th term} we use \eqref{size of mcl-B_2 and Sigma outside of it} and the definition of $\Xi$ to obtain
\[
\sum_\alpha \sum_{\delta,\gamma \in \cal{B}_2(\alpha)}
\int_0^1\! (\11-\theta)\,\absB{\,\E\Bigl[w(\alpha)\1w(\delta)\1w(\gamma)\2\partial^{\13}_{\alpha\delta\gamma}f(w^{\alpha,\theta})\Bigr]\,\dd \theta\,}
\;\leq\;
N^{1/2+C\eps}\,\Xi
\,.
\]
The last term \eqref{derivative of EE f(w_t): 5-th term} is estimated similarly
\[
\sum_{\gamma\in \cal{B}_2(\alpha)} 
\int_0^1\! (\11-\theta)\,
\E\Bigl[w(\gamma)\1\partial^{\13}_{\alpha\delta\gamma} f(w^{\alpha,\theta})\Bigr]\,\dd\theta
\;\leq\;
N^{-1/2+C\eps}\,\Xi
\,,
\]
and the double sum over $ \alpha,\delta $ produces a  factor of size $ CN $ due to the exponential decay of $\Sigma$. 
Combining the estimates for the five terms on the right hand side of \eqref{derivative of EE f(w_t)} we obtain \eqref{flow for matrix function}.
\end{Proof}

\begin{Proof}[Proof of   Corollaries~\ref{crl:Bulk universality} and \ref{crl:Bulk Gap Universality}]
We will only sketch the argument here as the procedure is standard. 
First we show that the matrix $ \bs{H}_t $ defined through \eqref{SDE for matrix flow} satisfies the bulk universality if $ t \ge t_N := N^{-1+\xi_1} $, for any $ \xi_1 > 0$. 
For simplicity, we  will focus on $t\le N^{-1/2}$ only. 
Indeed, from the fullness assumption {\bf B5} it follows that $ \bs{H}_t $ is of the form
\bels{H_t has G-component}{
\bs{H}_t = \wt{\bs{H}}_t + c(t)\1t^{1/2}\bs{U}
\,,
}
where $ c(t) \sim 1$  and $ \bs{U} $ is a GUE/GOE-matrix independent of  $ \wt{\bs{H}}_t $. 
For $ t \leq N^{-1/2} $ the matrix $ \wt{\bs{H}}_t $ has essentially the same correlation structure as $\bs{H}$, controlled by essentially the same model parameters. 
In particular the corresponding $ \wt{\cal{S}}_t $ operator is almost the same as $\cal{S}$. Let $ \wt{\bs{M}}_t $ solve the corresponding MDE with $ \cal{S} $ replaced by $ \wt{\cal{S}}_t$ and  let $ \wt{\rho}_t $ denote the function related to $ \wt{\bs{M}}_t $ similarly as $ \rho $ is related to $ \bs{M} $ (see Definition~\ref{def:Density of states}). 
Using the general stability for MDEs, Theorem~\ref{thr:Stability}, with
\[
\bs{\frak{G}}(\bs{0}) \,:=\,  \bs{M}
\,,\qquad
\bs{D} \,:=\, (\1\wt{\cal{S}}_t\msp{-1}-\cal{S}\1)[\2\wt{\bs{M}}_t]\,\wt{\bs{M}}_t\,,\qquad
\bs{\frak{G}}(\bs{D}) \,=\, \wt{\bs{M}}_t\,,
\]
(cf. \eqref{MDE stability2}) it is easy to check  that $ \wt{\bs{M}}_t $ is close to $\bs{M}$, in particular $ \wt{\rho}_t(\omega)  \ge \delta/2 $ when $ \rho(\omega) \ge \delta $. 
Moreover, the local law applies to  $ \wt{\bs{H}}_t $ as well, i.e. the resolvent $ \wt{\bs{G}}_t(\zeta) $ of $ \wt{\bs{H}}_t$ approaches  $ \wt{\bs{M}}_t(\zeta) $ for spectral parameters $\zeta$ with $\rho(\re \zeta)\ge \delta$. 
The bulk spectrum  of $ \wt{\bs{H}}_t $ is therefore the same as that of $ \bs{H} $ in the limit.
Combining these facts with the decomposition \eqref{H_t has G-component} we can apply  
 Theorem~2.2 from the recent work \cite{LSY} to conclude bulk universality for $ \bs{H}_t $, with $ t=t_N = N^{-1+\xi_1} $ in the sense of correlation functions as in Corollary~\ref{crl:Bulk universality}. 
 In order to prove the gap universality, Corollary~\ref{crl:Bulk Gap Universality}, we use Theorems~2.4 and 2.5 from \cite{LY}
 or Theorem~2.1 and Remark 2.2 from \cite{ES}. 

Second, we use  Lemma~\ref{lmm:Comparison flow} to show that $\bs{H}$ and $\bs{H}_t$ have the same 
local correlation functions in the bulk.
Suppose $ \rho(\omega) \ge \delta $ for some $ \omega \in \R $. 
We show that the difference 
\[
(\tau_1,\dots,\tau_k) \mapsto (\1\rho_{k\1;\1t_N}-\rho_k)\msp{-2}\bigl(\2\omega+{\textstyle \frac{\,\tau_{\11}\!}{N},\dots,\omega+\frac{\,\tau_{\1k}\!}{N}}\1\bigr)
\]
of the local $ k$-point correlation functions $ \rho_k $ and $ \rho_{k\1;\1t_N} $ of $ \bs{H} $ and $ \bs{H}_{t_N} $, respectively, converge weakly to zero as $ N \to \infty $.
This convergence follows from the standard arguments provided that 
\[
\abs{\2\E\,F(\bs{H}_{\1t})-\E\,F(\bs{H})\2} 
\,\to\, 0
\,,
\]
where $ F = F_N $ 
is a function of $ \bs{H} $ expressed as a smooth function $ \Phi $ of the following observables
\[
\qquad \frac{1}{N^p}\tr\,\prod_{j=1}^p(\1\bs{H}-\zeta_j^\pm\1)^{-1}
\,,\qquad
\zeta^\pm_j := \omega+\frac{\tau_{i_j}}{N}\pm\ii\1N^{-1-\xi_2}
,\quad j =1,\dots,p
\,,
\]
with $ p \leq k $ and $ \xi_2 \in (0,1) $ sufficiently small.
Here the derivatives of $ \Phi $ might grow only as a negligible power of $ N $ (for details see the proof of Theorem~6.4 in \cite{EYY}). 
In particular, basic resolvent formulas yield
\[
\text{RHS of \eqref{flow for matrix function}}
\,\leq\, C\1N^{\eps'}N^{C'\xi_2}\2\E\Bigl[(\11+\Lambda_t\1)^{C''}\Bigr]\,N^{1/2+\eps}\1t
\,,
\]
where $ \Lambda_t $ is defined like $ \Lambda $ in \eqref{definition of Lambda} but for the entries of $ \bs{G}_t(\zeta) := (\bs{H}_t-\zeta\1)^{-1} $ with $ \im\,\zeta \ge N^{-1+\xi_2} $.
In particular, we have used $ \abs{G_{xy}(t)} \leq \abs{m_{xy}(t)} + \Lambda_t \lesssim  1+\Lambda_t $ here.
The constant $ \Xi $ from \eqref{def of Gamma and wti-Gamma} is easily bounded by $ N^{\eps'+C\xi_2} $, where the arbitrary small constant $ \eps' > 0 $ originates from stochastic domination estimates for $ \Lambda_t $ and $ \abs{w_s(\alpha)}$'s.
The constant $ \wt{\Xi} $ from \eqref{def of Gamma and wti-Gamma}, on the other hand, is trivially bounded by $ N^{C} $ since the resolvents satisfy trivial bounds in the regime $ \abs{\1\im\,\zeta\1} \ge N^{-2} $, and the weight $ \abs{w(\alpha)} $ multiplying the third derivatives of $ f $ is canceled for large values of $ \abs{w(\alpha)} $ by the inverse in the definition $ \bs{G} = (\bs{A}+N^{-1/2}\bs{W}-\zeta\2\bs{1})^{-1} $.
Since the local law holds for $ \bs{H}_t $, uniformly in $ t \in [0,t_N]$, we see that $ \Lambda_t \leq N^{\eps'} (N\1\eta\1)^{-1/2} \leq N^{\eps'-\xi_2/2} $ with very high probability and hence
\bels{Comparison flow}{
\abs{\2\E\,F(\bs{H}_{t_N})-\E\,F(\bs{H})\2}  \leq C(\eps)\2N^{1/2+\eps}N^{-1+\xi_1}N^{C\1(\eps'+\xi_2)}
\,.
}
Choosing the exponents $ \eps,\eps',\xi_1,\xi_2 $ sufficiently small we see that the right hand side goes to zero as $N \to \infty $. This completes the proof of Corollary~\ref{crl:Bulk universality}.   Finally, the comparison estimate \eqref{Comparison flow} and the rigidity bound \eqref{bulk rigidity} allows
us to compare the gap distributions of $\bs{H}_{t_N}$ and $\bs{H}$, see Theorem~1.10 of \cite{KY2013}.  
This proves Corollary~\ref{crl:Bulk Gap Universality}.
\end{Proof}

\appendix

\section{Appendix}

\begin{Proof}[Proof of Lemma~\ref{lmm:Perturbed Combes-Thomas}]
Within this proof we adapt Convention~\ref{con:Comparison relation} such that $\varphi\lesssim \psi$ means  $\varphi\le C \psi$ for a constant $C$, depending only on ${\scr{P}}:=(\ul{\beta},P)$.  It suffices to prove \eqref{R inverse CT inequality} for $N \ge N_0$ for some threshold $N_0\lesssim 1$. Thus, we will often assume $N$ to be large enough in the following. 

We split $\bs{R}$ into a decaying  component $\bs{S}$ and an entrywise small component $\bs{T}$, i.e. we define 
\bels{splitting R into decaying and flat part}{
\bs{R}\,=\, \bs{S}+\bs{T}\,,\qquad s_{xy}\,:=\, r_{xy}\2\bbm{1}\pb{\2\abs{\1r_{xy}}\2\ge\2 {\textstyle \frac{2\1C_1}{N}}\1}\,.
}
The main part of the proof of Lemma~\ref{lmm:Perturbed Combes-Thomas} is to show that $\bs{S}$ has a bounded inverse,
\bels{bound on inverse of decaying part}{
\norm{\1\bs{S}^{-1}}\,\lesssim\, 1\,.
}
We postpone the proof of \eqref{bound on inverse of decaying part} and show first how it is used to establish \eqref{R inverse CT inequality}.

Since the entries of $\bs{S}$ are decaying as,  
$
\abs{s_{xy}}\,\le\,  \beta(\nu)\1(\11+d(x,y)\1)^{-\nu}
$,
for any $\nu \in \N$,  we can apply the standard Combes-Thomas estimate (Lemma~\ref{lmm:Perturbed Combes-Thomas} with $\alpha(0)=\beta(0)=0$) in order to get  the 
 decay of the entries of $\bs{S}^{-1}$ to arbitrarily high polynomial order $\nu \in \N$,
\bels{decaying entries of S inverse}{
\abs{(\bs{S}^{-1})_{xy}}\,\le\, \frac{C(\nu)}{(1+d(x,y))^\nu}\,.
}
In particular, we find that the $\norm{\2\cdot\2}_{1\vee \infty}$-norm (introduced in \eqref{ell-1, ell-infty and max-ell-1-infty matrix norms}) of $\bs{S}^{-1}$ is bounded
\bels{max 1 infty norm estimate on decaying part}{
\norm{\1\bs{S}^{-1}}_{1\vee \infty}\lesssim\, 1
\,.
}

We show now that $\norm{\bs{R}^{-1}\!-\bs{S}^{-1}}_{\rm{max}}\lesssim \frac{1}{N} $ which together with \eqref{decaying entries of S inverse} implies \eqref{R inverse CT inequality}. For a matrix $\bs{Q} \in \C^{N \times N}$  viewed as an operator mapping between $\C^{N}$ equipped with the standard Euclidean and the maximum norm we use the induced operator norms 
\[
\norm{\bs{Q}}_{2 \to \infty}\,:=\,\max_{\1x}{\textstyle \sqrt{\,\sum_{\1y}\2\abs{\1q_{xy}}^2}} 
\,,\qquad 
\norm{\bs{Q}}_{ \infty \to 2}\,:={\textstyle \sqrt{\,\sum_{\1x}\,\pb{\,\sum_{\1y}\2\abs{\1q_{xy}}\2}^{\!2}\,}}
\,.
\]
We write the difference between $\bs{R}^{-1}$ and $\bs{S}^{-1}$ as 
\bels{difference between S inverse and R inverse}{
-\bs{S}^{-1}\bs{T}\bs{R}^{-1}\,=\, \bs{R}^{-1}\!-\bs{S}^{-1}\,=\, -\bs{R}^{-1}\bs{T}\bs{S}^{-1}.
}
The first equality in \eqref{difference between S inverse and R inverse} implies 
\bels{infty norm bound on R}{
\norm{\bs{R}^{-1}}_\infty
\,\le\, 
\norm{\1\bs{S}^{-1}}_\infty
\bigl(\21+\norm{\bs{T}}_{2 \to \infty}\norm{\bs{R}^{-1}}_{\infty \to 2}\bigr)
\,\leq\,
 \norm{\1\bs{S}^{-1}}_\infty
 \bigl(\21+N\norm{\bs{T}}_{\rm{max}}\norm{\bs{R}^{-1}}\bigr)
\,\lesssim\, 1
\,,
}
where we used 
$
\norm{\bs{Q}}_{\infty \to 2}\le\sqrt{N}\norm{\bs{Q}}
$
and 
$ 
\norm{\bs{Q}}_{2 \to \infty}\le \sqrt{N}\norm{\bs{Q}}_{\rm{max}}
$
 for any $\bs{Q} \in \C^{N \times N}$, \eqref{max 1 infty norm estimate on decaying part} and $\norm{\bs{T}}_{\rm{max}}\lesssim N^{-1}$ from the definition of $\bs{T}$ in \eqref{splitting R into decaying and flat part}.
The second equality in \eqref{difference between S inverse and R inverse} on the other hand implies 
 $
\norm{\2\bs{R}^{-1}\!-\bs{S}^{-1}}_{\rm{max}}
\le
\norm{\1\bs{R}^{-1}}_\infty
\norm{\bs{T}}_{\rm{max}}
\norm{\1\bs{S}^{-1}}_1
\lesssim 
N^{-1}$, %
where  \eqref{max 1 infty norm estimate on decaying part} and  \eqref{infty norm bound on R} were  used in the second inequality. This finishes the proof of Lemma~\ref{lmm:Perturbed Combes-Thomas} up to verifying \eqref{bound on inverse of decaying part}.

We split $\bs{R}^{\!*}\bs{R}$ into a decaying and  an   entrywise small piece as we did with $\bs{R}$ itself in \eqref{splitting R into decaying and flat part},
\[
\bs{R}^{\!*}\bs{R}\,=\, \bs{L}+\bs{K}\,,\qquad \bs{L}\,:=\, \bs{S}^*\bs{S}\,,\qquad \bs{K}\,:=\, \bs{S}^*\bs{T}+\bs{T}^*\bs{S}+\bs{T}^*\bs{T}\,.
\]
From the related properties of $\bs{S}$ and $\bs{T}$ we can easily see that 
\bels{decay of entries of L and smallness of entries of K}{
\abs{\2l_{xy}}\,\le\, \frac{C(\nu)}{(1+d(x,y))^\nu}
\,,\qquad 
\norm{\bs{K}}_{\rm{max}}\,\lesssim\, \frac{1}{N}
\,,
}
where $ \bs{L} = (l_{xy})_{x,y} $. 
Using the a priori knowledge 
\bels{L+K inverse bound}{
{\bs{L}+\bs{K}}\,=\,\bs{R}^{\!*}\bs{R}
\;\gtrsim\; 
1
\,,
}
from the assumption $\norm{\bs{R}^{-1}}\lesssim 1$ of Lemma~\ref{lmm:Perturbed Combes-Thomas}, we will show that $\norm{\bs{L}^{-1}}\lesssim 1$, which is equivalent to \eqref{bound on inverse of decaying part}.   Note that both $\bs{L}$ and $\bs{K}$ are selfadjoint. 

Via spectral calculus we write $\bs{K}$ as a sum of a matrix  $\bs{K}_s$ with small spectral norm and a matrix  $\bs{K}_b$ with bounded rank
\bels{splitting K into finite rank and small part}{
\bs{K}\,=\, \bs{K}_s\msp{-2}+\bs{K}_b
\,,\qquad 
\bs{K}_s
\,:=\, 
\bs{K}\2 \bbm{1}_{(-\eps,\eps)}\msp{-1}(\bs{K})
\,.
}
with some $\eps>0$ to be determined later. 
Indeed,  from the Hilbert-Schmidt norm bound on the eigenvalues $\lambda_i(\bs{K})$ of $\bs{K}$,
\[
{\textstyle \sum_{\2i}}\2\lambda_i(\bs{K})^2
\,=\,
\tr \bs{K}^*\bs{K}
\,\le\, 
N^2\norm{\bs{K}}_{\rm{max}}^2
\,\lesssim\, 
1
\,,
\]
we see that $\rm{rank}\, \bs{K}_b \,\lesssim\, \frac{1}{{\eps^2}}$. On the other hand $\norm{\bs{K}_s}\le \eps$ by its definition in  \eqref{splitting K into finite rank and small part}. Since $\bs{L}+\bs{K}$ has a bounded inverse (cf. \eqref{L+K inverse bound}), so does  $\bs{L}+\bs{K}_b$ for small enough $\eps$, i.e.
\bels{L + finite rank is invertible}{
\norm{\1(\1\bs{L}+\bs{K}_b)^{-1}}\,\lesssim\, 1\,.
}

Now we fix $\eps \sim 1$ such that \eqref{L + finite rank is invertible} is satisfied.  In particular the eigenvalues of $\bs{L}+\bs{K}_b$ are separated away from zero. 
 Since $\rm{rank}\, \bs{K}_b \,\lesssim\,1$, we can apply  the interlacing property of rank one perturbations finitely many times  to 
 see that there are only finitely many eigenvalues of $\bs{L}$ in a neighborhood of zero, i.e.
\bels{finite rank of L below a threshold}{
\rm{rank}\1[\2\bs{L} \2\bbm{1}_{[0,c_1)}( \bs{L})]\,\lesssim\, 1\,,
}
for some constant $c_1 \sim 1$. In particular, there are constants $c_2\sim c_3 \sim 1$ such that $c_2+c_3 \le c_1$ 
 and $\bs{L}$ has a spectral gap at $[c_2,c_2+c_3]$,
\bels{spectral gap of L}{
\bs{L} \2\bbm{1}_{[c_2,c_2+c_3]}( \bs{L})\,=\, \bs{0}\,.
}
We split $\bs{L}$ into the finite rank part  $\bs{L}_s$   associated to the spectrum below the gap and the rest,
\[
\bs{L}\,=\, \bs{L}_s + \bs{L}_b\,,\qquad \bs{L}_s\,:=\, \bs{L} \2\bbm{1}_{[0,c_2)}( \bs{L})\,,\qquad \bs{L}_b\,:=\,\bs{L} \2\bbm{1}_{(c_2+c_3, \infty)}( \bs{L})\,.
\]

The rest of the  proof   is devoted to showing that $\bs{L}_s \gtrsim \bbm{1}_{[0,c_2)}( \bs{L})$, 
  which implies that $\bs{L}$ has a bounded inverse and thus shows \eqref{bound on inverse of decaying part}. More precisely, we  will  show that there are points $x_1, \dots, x_L$ with $L:=\rm{rank}\, \bs{L}_s\lesssim1 $ (cf. \eqref{finite rank of L below a threshold}) and a positive sequence $(C(\nu))_{\nu \in \N}$ such that 
\bels{decay of all eigenvectors close to zero}{
\norm{\bbm{1}_{[0,c_2)}( \bs{L})\1\bs{e}_x}\,\le\, \sum_{i=1}^L \frac{C(\nu)}{(1+d(x_i,x))^\nu}\,,
}
for any $\nu \in \N$ and $x\in \bb{X}$, where $(\bs{e}_x)_{x\in \bb{X}}$ denotes the canonical basis of $\C^N$.  
 Let $\bs{l}= (l_x)$ be any normalized eigenvector of $\bs{L}_s$ in the image of $\bbm{1}_{[0,c_2)}( \bs{L})$ with  associated eigenvalue $\lambda$. We need to show that $\lambda\gtrsim 1$.  
Since $\langle \bs{l},  \bbm{1}_{[0,c_2)}( \bs{L}) \bs{e}_x\rangle =  l_x$,
the  decay property \eqref{decay of all eigenvectors close to zero} of the spectral projection $\bbm{1}_{[0,c_2)}( \bs{L})$
 away from the finitely many centers $x_i$  implies that the components $l_x$ have arbitrarily high polynomial decay  away from  the points $x_1, \dots, x_L$.  
In particular, $\sum_{\1x} \abs{\1l_x}$ is bounded  and therefore we have (cf. \eqref{decay of entries of L and smallness of entries of K})
\[
\abs{\2\bs{l}^*\1\bs{K}\2\bs{l}\2}
\,\le\, 
\norm{\bs{K}}_{\rm{max}}\pb{\,\textstyle \sum_{\2x}\2\abs{\2l_x}\,}^2
\,\lesssim\, 
N^{-1}
\,.
\]
We infer that for the eigenvalue $\lambda$ we get a lower bound via
$
1 \lesssim \bs{l}^*(\bs{L}+\bs{K})\bs{l}= \lambda +\bs{l}^*\bs{K}\bs{l}
$, 
where we used \eqref{L+K inverse bound} for the inequality. Thus, $\lambda \gtrsim 1$ for large enough $N$. 

Now we prove \eqref{decay of all eigenvectors close to zero} by induction. We show that for any $l=0,\dots,L$ there is an $l$-dimensional subspace of the image of $\bbm{1}_{[0,c_2)}( \bs{L})$  such that   the associated orthogonal projection $\bs{P}_l$ satisfies
\bels{decay of subspace projection}{
\norm{\bs{P}_l\bs{e}_x}\,\le\, \sum_{i=1}^l \frac{C(\nu)}{(1+d(x_i,x))^\nu}\,.
}
The induction is over $l$. For $l=0$ there is nothing to show. Now suppose that \eqref{decay of subspace projection} has been established for some $l<L$. We will see now that it then holds for $l$ replaced by $l+1$ as well. 

We maximize the maximum norm of all vectors in the image of $\bbm{1}_{[0,c_2)}( \bs{L})-\bs{P}_l$ and pick the index $x_{l+1}$ where the maximum is attained,
\bels{definition of the maximum alpha}{
\xi
\,:=\, 
{\textstyle \max_{\2x}}\, \norm{(\bbm{1}_{[0,c_2)}\msp{-1}\msp{-1}( \bs{L})-\bs{P}_l)\2\bs{e}_x}
\,=\, 
\norm{(\bbm{1}_{[0,c_2)}\msp{-1}( \bs{L})-\bs{P}_l)\1\bs{e}_{x_{l+1}}}\,.
}
Here, $\xi>0$ since $l<L$.
Now we extend the projection $\bs{P}_l$ by the normalized vector $\bs{v}$ defined as
\bels{definition of v}{
\bs{P}_{l+1}\,:=\,\bs{P}_l+{\bs{v}\2\bs{v}^*} \,,\qquad
\bs{v}\,:=\, \frac{1}{\xi}\,
(\bbm{1}_{[0,c_2)}( \bs{L})-\bs{P}_l)\1\bs{e}_{x_{l+1}}\,.
}
The so defined 
vector $\bs{v}$ attains its maximum norm at the point $x_{l+1}$ and the value of this norm is $\xi$,  since for any $x$ we have 
\bels{estimate on components of v}{
\abs{v_x}
\,=\,
\abs{\,\bs{v}^*(\1\bbm{1}_{[0,c_2)}\msp{-1}( \bs{L})-\bs{P}_l)\1\bs{e}_{x}} 
\,\le\, 
\xi
\,=\,
\frac{1}{\xi}\2\bs{e}_{x_{l+1}}^*\!(\bbm{1}_{[0,c_2)}( \bs{L})-\bs{P}_l)\1\bs{e}_{x_{l+1}}
=\, 
v_{x_{l+1}}
\,.
}
 Here we used  \eqref{definition of the maximum alpha} and that $\bbm{1}_{[0,c_2)}( \bs{L})-\bs{P}_l\ge \bs{0}$ is an orthogonal projection. 

We will show that $\bs{P}_{l+1}$ satisfies \eqref{decay of subspace projection} with $l$ replaced by $l+1$.
We start by establishing that $\xi\gtrsim 1$. We write $\norm{\bs{v}}^2$ as a sum of contributions originating from the neighborhoods $B:=\bigcup_{i=1}^{l+1}B_{R}(x_i)$ of the points $x_i$  with some radius $R$ to be determined later and their complement. We  estimate the components of $\bs{v}$ by using \eqref{estimate on components of v} and the definition of $\bs{v}$ in \eqref{definition of v}, 
\bels{estimate for alpha start}{
1 \,=\, \norm{\bs{v}}^2 
\,=\, 
\sum_{y\1\in\1 B}\abs{v_y}^2 +\sum_{y}^B\2\abs{v_y}^2
\,\le\, 
\abs{B}\1 \xi^2 
+
\frac{1}{\xi^2}\sum_{y}^B\2\pb{\2\abs{\1\bs{e}_y^*\bbm{1}_{[0,c_2)}( \bs{L})\1\bs{e}_{x_{l+1}}}+\abs{\1\bs{e}_y^*\bs{P}_l\bs{e}_{x_{l+1}}}\2}^2
.
}
With the induction hypothesis \eqref{decay of subspace projection} the second summand in the sum on the right hand side of \eqref{estimate for alpha start}  is bounded by
\bels{decay estimate of second term}{
\abs{\2\bs{e}_y^*\1\bs{P}_l\bs{e}_{x_{l+1}}}
\,\le\, 
\norm{\1\bs{P}_l\bs{e}_y} 
\,\le\, 
{\textstyle C(\nu)\1\sum_{\2i\2=\21}^{\2l} (\11\!+\!R\1)^{-\1\nu}}
\,,
}
for $y\not \in B$. 
For the  other summand in \eqref{estimate for alpha start} 
 we use the decay estimate
\bels{decay estimate on P}{ 
\abs{\1\bs{e}_y^*\bbm{1}_{[0,c_2)}( \bs{L})\bs{e}_{x}}= 
\abs{( \bbm{1}_{[0,c_2)}(\bs{L}))_{yx}}
\,\le\, 
C(\nu)\2(\21\1+d(x,y))^{-\nu} 
\,,\qquad \nu \in \N\,.
}
The bound \eqref{decay estimate on P} follows from the  integral representation
\bels{integral representation of P}{
 \bbm{1}_{[0,c_2)}(\bs{L})\,=\, {\textstyle \oint_{\1\Gamma}} (\bs{L}-\zeta\1\bs{1})^{-1}\dd \zeta\,,
}
where the integral is over a closed contour $\Gamma$ encircling only the eigenvalues of $\bs{L}$ within $[0,c_2)$. 
Since $\bs{L}$ has a spectral gap above $c_2$ (cf. \eqref{spectral gap of L}) we may choose $\Gamma$ such that
\[
{\textstyle \max_{\1\zeta \1\in\1 \Gamma}}\norm{\1(\1\bs{L}-\zeta\1\bs{1})^{-1}}\,\lesssim\, 1\,.
\]
Since the entries of $\bs{L}$ are decaying by \eqref{decay of entries of L and smallness of entries of K}, we can apply 
the standard Combes-Thomas estimate
to see that the entries of $(\bs{L}-\zeta\1\bs{1})^{-1}$ decay as well. Then \eqref{decay estimate on P} follows from \eqref{integral representation of P}.

Using \eqref{decay estimate of second term} and \eqref{decay estimate on P} in \eqref{estimate for alpha start} yields
\bels{alpha bound}{
1\,\le\, \abs{B}\1\xi^2 +\frac{C(\nu)\abs{B}}{\xi^2R^\nu}\,\lesssim\, R^{P}\xi^2+\frac{C(\nu)}{\xi^2R^{\nu-P}}\,
}
where in the second inequality we estimated the size of $B$ with \eqref{polynomial ball growth}.  Now we choose $R:=\xi^{-1/{P}}$, $\nu:=\lceil 4P \rceil$.  Using $\xi \le 1$ (cf. \eqref{definition of the maximum alpha}), we obtain that the right hand side \eqref{alpha bound} is bounded
by a constant multiple of  $\xi$. Thus \eqref{alpha bound} proves $\xi \gtrsim 1$. 

We finish the induction by using the definition \eqref{definition of v} of $\bs{v}$ and  estimating
\[
\norm{\1\bs{P}_{\!l+1}\bs{e}_x}
\,\le\, 
\norm{\1\bs{P}_{l}\bs{e}_x} + \abs{\1v_x}
\,\le\,   
\norm{\1\bs{P}_{l}\bs{e}_x}+\frac{1}{\xi}\pb{\,\norm{\1\bs{P}_{l}\bs{e}_x} + \abs{\2 (\bbm{1}_{[0,c_2)}(\bs{L}))_{x_{l+1} x }}\,}
\,.
\]
Since $\xi \gtrsim 1$ and by the induction hypothesis \eqref{decay of subspace projection} as well as \eqref{decay estimate on P}  the bound \eqref{decay of subspace projection} with $l$ replaced by $l+1$ follows and Lemma~\ref{lmm:Perturbed Combes-Thomas} is proven. 
\end{Proof}

Now, we generalise the following result to the non-commutative setting:
\begin{lemma}[\cite{AEK1cpam}, Lemma~5.6]
\label{lmm:basic sp-gap for matrices}
A symmetric matrix $ \bs{S}=(s_{ij})_{i,j=1}^N $, with non-negative entries, has a spectral gap of size at least 
$ \norm{\bs{v}}/(N^{1/2}\norm{\bs{v}}_{\rm{max}})\, \min_{\1i,\1j} s_{ij} $, where $ \bs{v} \in \C^N $ satisfies $ \bs{S}\1\bs{v} = \norm{\bs{S}}\1\bs{v} $.
\end{lemma}
\begin{Proof}[Proof of Lemma \ref{lmm:Spectral gap}]
Since $\cal{T}$ leaves the cone $ \scr{C}_+ $ of positive semidefinite matrices invariant, the Perron-Frobenius  theorem  guarantees the existence of a normalized $\bs{T}\in \scr{C}_+$ with $\cal{T}[\bs{T}]=\bs{T}$. We  first  verify the bounds \eqref{bounds on eigenmatrix T} on this eigenmatrix.

From the upper and lower bounds \eqref{bounds on operator T} on $\cal{T}$ we infer
\bels{bounds on T in terms of average}{
\gamma\,\avg{\bs{T}}\2\bs{1} \,\leq\, \bs{T}\,\leq\, \Gamma\2\avg{\bs{T}}\2\bs{1}\,.
}
Multiplying by $\bs{T}$ on both sides of the second inequality and taking the normalized trace yields
$
1=\norm{\bs{T}}^2_{\rm{hs}}\leq \Gamma \1 \avg{\bs{T}}^2
$.
With the lower bound from \eqref{bounds on T in terms of average} on $\bs{T}$ we see that $\bs{T} \geq \frac{\gamma}{\sqrt{\Gamma}}\2\bs{1}$.
Furthermore, the normalization of $\bs{T}$ and the upper bound from \eqref{bounds on T in terms of average} imply 
\[
\bs{T} \,\leq\, \Gamma\1\avg{\bs{T}}\2\bs{1}\,\le\, \Gamma\1\norm{\bs{T}}_{\rm{hs}}\2\bs{1}\,=\, \Gamma\2\bs{1}.
\] 

Now we show the existence of a spectral gap and  that $1$ is a non-degenerate eigenvalue.
Showing \eqref{Spectral gap for T} is equivalent to proving that
\bels{equivalent to spectral gap}{
\qquad 
\scalar{\1\bs{R}}{(\1\rm{Id}\pm\cal{T}\2)[\bs{R}]\1}\,\geq\, %
\tsfrac{\gamma^{\16}}{2\1\Gamma^4}\,,
}
holds for all $ \bs{R}=\bs{R}^{\!\ast}\in \C^{N \times N} $ satisfying $ \norm{\bs{R}}_{\rm{hs}}=1$ and $ \scalar{\bs{T}}{\bs{R}}=0 $.
Here $ \bs{R} $ can be assumed to be self-adjoint since $\cal{T}$ preserves $\ol{\scr{C}}_+$, and thus $\cal{T}[\bs{R}]^*=\cal{T}[\1\bs{R}^{\!*}]\,$: 
\[
\scalar{\1\bs{R}}{(\1\rm{Id}\pm\cal{T}\2)[\bs{R}]\1}
\,=\, 
\scalar{\1\re \bs{R}}{(\1\rm{Id}\pm\cal{T}\2)[\re \bs{R}]\1} 
+ 
\scalar{\1\im \bs{R}}{(\1\rm{Id}\pm\cal{T}\2)[\im \bs{R}]\1}
\,.
\]

Let $\bs{R} $ be an arbitrary normalized self-adjoint matrix satisfying $ \scalar{\bs{T}}{\bs{R}}=0 $. 
We use the spectral representation 
$
\bs{R} =  \sum_{\1i} \varrho_i\,\bs{r}_i\2\bs{r}_i^*
$, 
with the orthonormal eigenbasis $(\bs{r}_i)_{i=1}^N$ of $\bs{R}$. Plugging  this spectral representation  into the right hand side of
\eqref{equivalent to spectral gap} reveals the identity 
\bels{scalar product in terms of S}{
\scalar{\1\bs{R}}{(\1\rm{Id}\1\pm\cal{T}\2)[\bs{R}]\1}
\,=\,
\bs{q}^*(\bs{1} \pm \bs{S})\bs{q}\,,
}
where we introduced the vector $\bs{q} \in \R^N$ 
of eigenvalues of $\bs{R}$ and the matrix $\bs{S} \in \R^{N \times N}$ with non-negative entries:
\[
q_i\,:=\,
N^{-1/2}\varrho_i
\,,\qquad s_{i j}\,:=\, {\bs{r}_i^*\cal{T}[\1\bs{r}_j\bs{r}_j^*]\2\bs{r}_i}\,.
\]
The vector $\bs{q}$ is normalized since $\norm{\bs{q}}=\norm{\bs{R}}_{\rm{hs}}=1$, and the matrix $\bs{S}$ is symmetric because of the self-adjointness of $\cal{T}$.
Furthermore, by \eqref{bounds on operator T} the entries of $\bs{S}$ satisfy lower and upper bounds,
\bels{bounds on matrix S}{
\gamma\1N^{-1} \leq s_{ij} \leq \Gamma\1N^{-1}.
}
In particular, by the Perron-Frobenius  theorem,  the matrix $\bs{S}$ has a unique normalized eigenvector $\bs{s}$ with 
positive entries  and with associated eigenvalue equal  to its spectral norm,  $ \bs{S}\2\bs{s} = \norm{\bs{S}}\1\bs{s} $.

We will now show that $\bs{S}$ has a spectral gap and $\bs{r}$ has a non-vanishing component in the direction orthogonal to $\bs{s}$. This will imply 
\bels{equivalent to spectral gap with S}{
\abs{\1\bs{q}^* \bs{S}\1\bs{q}\2}\,\le\, 1-\tsfrac{\gamma^{\16}}{2\1\Gamma^4}
\,,
}
which is equivalent to \eqref{equivalent to spectral gap} by \eqref{scalar product in terms of S} and therefore proves Lemma~\ref{lmm:Spectral gap}.

To verify \eqref{equivalent to spectral gap with S} we start with the observation that the norm of $\bs{S}$ is bounded by
\bels{bounds on norm of matrix S}{
c\,\le\, \bs{e}^*\bs{S}\2\bs{e}\,\le\, \norm{\bs{S}}\,=\, {\textstyle \sup_{\2\norm{\bs{w}}=1}\,}\bs{w}^*\cal{T}[\bs{w}\2\bs{w}^*]\bs{w}\,\le\, \norm{\cal{T}}_{\rm{sp}}\,=\,1\,,
}
where $\bs{e}=(\frac{1}{\sqrt{N}}, \dots,\frac{1}{\sqrt{N}})$,  and that the Perron-Frobenius eigenvector $\bs{s}=(s_i)_{i=1}^N $ satisfies 
\[
\max_i s_i
\,=\, 
\max_i \frac{ (\bs{S}\1\bs{s})_i}{\norm{\bs{S}}}
\,\le\, 
\frac{\Gamma}{\gamma\1N} \sum_i s_i\,\le\, \frac{\Gamma}{\2\gamma\sqrt{\msp{-1}N\1}}
\,,
\]
where we used  $ \bs{S}\2\bs{s} = \norm{\bs{S}}\1\bs{s} $, \eqref{bounds on norm of matrix S}, \eqref{bounds on matrix S} and $\norm{\bs{s}}=1$ in that order. 
 Applying, Lemma \ref{lmm:basic sp-gap for matrices} yields 
\bels{spectral gap in matrix S}{
\spec\p{\bs{S}}
\,\subseteq\, 
\sb{-\norm{\bs{S}}+\tsfrac{\1\gamma^3}{\2\Gamma^2\!},\norm{\bs{S}}-\tsfrac{\1\gamma^3}{\2\Gamma^2\!}}\cup \bigl\{\1\norm{\bs{S}}\1\bigr\}\,.
}

Finally we show that there is a non-vanishing component of $ \bs{q} $ in the direction orthogonal to $\bs{s}$. More precisely, we show that there is some sufficiently large vector $ \bs{w} \perp \bs{s} $ satisfying:
\bels{splitting of vector of eigenvalues of R}{
\bs{q}\,=\, (\11-\norm{\bs{w}}^2)^{1/2}\2\bs{s} +\bs{w}
\,.
}
Taking the scalar product with $\bs{t}:=(N^{-1/2}\bs{r}_i^* \bs{T}\bs{r}_i)_{i=1}^N$ and using $\bs{t}^*\bs{q}=\scalar{\bs{T}}{\bs{R}}=0 $ yields
\bels{w bound from below inequalities}{
(\11-\norm{\bs{w}}^2)\,\tsfrac{\gamma^{\12}}{\2\Gamma N\1}\pb{\2{\textstyle \sum_{\2i}}\1 s_i\,}^2
\,\le\,  
(\11-\norm{\bs{w}}^2)\2(\1\bs{t}^*\bs{s})^2
\,=\,
(\1\bs{t}^*\bs{w})^2 
\,\le\, 
\norm{\1\bs{t}\1}^2\norm{\bs{w}}^2
\,.
}
The first inequality in \eqref{w bound from below inequalities} follows from the lower bound on $\bs{T}$ in \eqref{bounds on eigenmatrix T}.  Since $\norm{\bs{t}}\le\norm{\bs{T}}_{\rm{hs}}=1$ and (cf. \eqref{bounds on norm of matrix S})
\[
\gamma 
\,\le\, 
\norm{\bs{S}}
\,=\, 
\bs{s}^*\bs{S}\2\bs{s}
\,\le\, 
\Gamma N^{-1}\pb{\,\textstyle \sum_{\2i} s_i\,}^2 ,
\]
we conclude that 
$\norm{\bs{w}}^2 \ge \frac{\gamma^{\13}}{\Gamma^2+\gamma^3} \ge  \frac{\gamma^{\13}}{2 \1\Gamma^{\12}} $,
where we used $\gamma \le \min\sett{1,\Gamma} $.
Combining this with $ \bs{w} \perp \bs{s} $, \eqref{spectral gap in matrix S} and \eqref{splitting of vector of eigenvalues of R}  yields
\[
\abs{\bs{q}^* \bs{S}\bs{q}}
\,\le\,
\norm{\bs{S}}\1(\11-\norm{\bs{w}}^2)+\pb{\norm{\bs{S}}-\tsfrac{\gamma^{\13}}{\Gamma^2\!}}\norm{\bs{w}}^2
\,\le\, 
1-\tsfrac{\gamma^{\13}}{\Gamma^2\!}\1
\norm{\bs{w}}^2
\,\le\, 
1-\tsfrac{\gamma^{\16}}{2\1\Gamma^4} = 1 -\theta
\,,
\]
where we also used $\norm{\bs{S}}\le 1$. Thus, \eqref{equivalent to spectral gap with S} is established and Lemma~\ref{lmm:Spectral gap} is proven.
\end{Proof}

\begin{lemma}[Linear large deviation]
\label{lmm:Linear large deviation}
Let $X=(X_x)_{x\in \bb{X}}$ and  $ b =  (b_x)_{x\in \bb{X}}$ be sequences of random variables that satisfy the following assumptions:
\begin{itemize}
\titem{i}
The entries of $X$ are centred, $\E\2X_x=0$. 
\titem{ii}
The entries of $ X $ have uniformly bounded moments, i.e.,  there is a sequence $ \ul{\beta}_1$ of positive constants such that  $ \E \2\abs{X_x}^\mu \,\le\, \beta_1(\mu)$, for all $ x $ and $ \mu \in \N $.
\titem{iii}
The correlations within $X$ decay, i.e.,  there is a sequence $ \ul{\beta}_2 $ of positive constants such that for every $ \eps > 0 $, every $A,B\subseteq\bb{X}$ satisfying $d(A,B):=\min\{\2d(x,y):x \in A,y\in B\1\} \ge N^\eps $, and all smooth  functions $\phi:\C^{\abs{A}} \to \C$, $\psi:\C^{\abs{B}} \to \C$, the quantities $X_A:=(X_x)_{x \in A} $ satisfy:
\bels{correlation decay in large deviation}{
\absb{\1\rm{Cov}(\phi(X_A),\psi(X_B))}\,\le\,\beta_2(\eps,\nu)\1 
\norm{\nabla \phi}_\infty\norm{\nabla \psi}_\infty N^{-\nu}
\,,\quad \nu \in \N
\,.
}
\titem{iv} 
The correlations between $X$ and $b$ are asymptotically small, i.e., there exists  a sequence $ \ul{\beta}_3 $ of positive constants, such that for  all  smooth functions $\phi,\psi:\C^N \to \C\,$ the following holds:
\bels{small correlation between X and b}{
\abs{\rm{Cov}(\phi(b),\psi(X))}
\,\le\,\beta_3(\nu)\1\norm{\nabla \phi}_\infty
\norm{\nabla \psi}_\infty N^{-\nu}
,\qquad \nu \in \N
\,.
}
\end{itemize}
Then the following large deviation estimate holds for every $\nu \in \N$:
\bels{Single sum large deviation}{
\absB{\sum_x b_x\2 X_x}\,\prec\,\pB{\2\sum_x\2\abs{\1b_x}^2}^{\!1/2}\msp{-8}+
\frac{1}{N^\nu}
\,.
}
\end{lemma}
\begin{Proof} Here we use Convention~\ref{con:Comparison relation} such that $\varphi\lesssim \psi$ means  $\varphi\le C \psi$ for a constant $C$, depending only on $\wt{\scr{P}}:=(\beta_1,\beta_2,\beta_3,P)$ (cf. \eqref{polynomial ball growth}). We divide the proof into three steps. 
\medskip\\
\emph{Step 1:} In this step we introduce a cutoff both for $X$ and $b$. We  show that
it suffices to prove the moment bound  
\bels{simplified large deviation lemma}{
\E\2\absb{\textstyle \sum_{\1 x}\msp{-3} b_x X_x}^{2\mu}\,\le\,C(\mu, \delta)N^{\delta}\,, \qquad \delta>0\,,
}
for two families  $X$ and $b$ of random variables that satisfy the upper bounds 
\bels{additional bounds}{
\max_x \abs{X_x}\,\le\, \sqrt{N}\,,\qquad {\textstyle \sum_{\2x}} \abs{\1b_x}^2\,\le\, 1\,,
}
in addition to 
the assumptions of Lemma~\ref{lmm:Linear large deviation}.

Indeed, for $X$ and $b$ as in Lemma~\ref{lmm:Linear large deviation} we define the new random variables
\bels{definition tilded variables}{
\wt{X}_x\,:=\, (1-\E)[\1X_x\, \theta\p{N^{-1}\abs{X_x}^2}\1]
\,,\qquad 
\wt{b}_x\,:=\, \frac{b_x}{\,(\,\sum_{\1y}\abs{\1b_y}^2\2)^{1/2}\!+N^{-\wt \nu}}\,,
}
where $\wt \nu \in \N$ and $\theta: [0,\infty) \to [0,1]$ is a smooth   cutoff  function such that $\theta|_{[0,1/2]}=1$ and $\theta|_{[1,\infty)}=0$.

It is easy to verify that $\wt X$ and $\wt b$ satisfy the assumption of Lemma~\ref{lmm:Linear large deviation}.
Now suppose that  \eqref{simplified large deviation lemma} holds with $X,b$ replaced by $\wt X, \wt b$. In particular, $\absb{\sum_{\1x} \wt{b}_x \wt{X}_x}\prec 1$. Then we see that 
\bels{conclusion form lemma with tilde}{
\textstyle
\absb{ \sum_{\1x} \wt{b}_x X_x}
\,\prec\,
\absb{\sum_{\1x} \wt{b}_x \wt{X}_x}+N^{-\nu}
\prec\, 1
\,,
}
for any $\nu \in \N$, 
where we used $\abs{X_x-\wt X_x}\prec N^{-\nu-1}$ 
and $\abs{\1\wt b_x}\le 1$. Plugging the definition \eqref{definition tilded variables} of $\wt{b}$ into \eqref{conclusion form lemma with tilde} yields 
\[
\textstyle 
\absb{\sum_{\2x} {b}_x X_x}\,\prec\,
\sqrt{\2\sum_{\2x} \abs{\1{b}_x}^2\2} +N^{-\wt \nu}
\,,
\]
and since $\wt \nu$ was arbitrary, Lemma~\ref{lmm:Linear large deviation} is proven, up to checking \eqref{simplified large deviation lemma} for random variables $X$ and $b$ that satisfy \eqref{additional bounds} in addition to the assumptions of the lemma. 
\medskip\\

\emph{Step 2:} In this step  we completely remove the weak dependence between $X$ and $b$, i.e. 
we show that it is enough to prove \eqref{simplified large deviation lemma} for a centered sequence $X$ independent of $b$ satisfying  \eqref{additional bounds}, the assumption (ii), %
and \eqref{correlation decay in large deviation}. Indeed, suppose that  $X$ and $b$ are not independent, but satisfy \eqref{small correlation between X and b} instead. Let $\wt b$ be a copy of $b$ that is independent of $X$ and $b$. We show that for any $\mu,\nu \in \N$, 
\bels{replace b by independent copy}{
\absB{\,\E\2\absb{\1{\textstyle \sum_{\1x}} \wt{b}_x X_x}^{2\mu}\!-\2\E\2\absb{\1{\textstyle \sum_{\1x}} b_x X_x}^{2\mu}}\,\le\, C(\mu, \nu)\1N^{- \nu}.
}
The bound \eqref{replace b by independent copy} implies \eqref{simplified large deviation lemma}, provided \eqref{simplified large deviation lemma} holds with $b$ replaced by $\wt b$. 

To prove \eqref{replace b by independent copy} we expand the powers on the left hand side, compare term by term and find
\[
\text{l.h.s. of \eqref{replace b by independent copy}}
\,\le\, 
N^{2 \mu}\max_{x_1, \dots , x_{2\mu}}\!\absb{\,\rm{Cov}(X_{x_1}\dots X_{x_{\mu}}\ol{X_{x_{\mu+1}}\dots X_{x_{2\mu}}},b_{x_1}\dots b_{x_\mu}\ol{b_{x_{\mu+1}}\dots b_{x_{2\mu}}})\2}\,,
\]
where the maximum is taken over all $x_1, \dots , x_{2\mu} \in \bb{X}$. Now we employ \eqref{small correlation between X and b} 
 as well as the bounds $\abs{b_x}\le 1$ and $\abs{X_x}\le \sqrt{N}$ to infer \eqref{replace b by independent copy}.
\medskip\\
\emph{Step 3:} By Step 1 and Step 2 we may assume for the proof of \eqref{simplified large deviation lemma} that $X$ is independent of $b$ and that these random vectors satisfy \eqref{additional bounds},  the hypothesis (ii)  and \eqref{correlation decay in large deviation}. In this final step we construct for every $\eps>0$ a partition of $\bb{X}$ into non-empty sets $I_1, \dots, I_K$ with the following properties:
\begin{itemize}
\item[(P1)]\label{P1:bound on size of partition} 
With a constant $C(\eps)$, depending only on $\eps$ and $P$ (cf. \eqref{polynomial ball growth}), the size of the partition is bounded by $ K \leq C(\eps)\1 N^{(P+1)\eps} $;
\item[(P2)]\label{P2:distance between indices in Ik} 
The indices within an element of the partition are far away from each other, i.e., if $ x,y \in I_k $, $ x\neq y $, then
$ d(x,y) \ge N^\eps $. %
\end{itemize}
 In other words, the elements within each $I_k$ are far from each other hence the corresponding components of $X$ and $b$ are 
practically independent. 

We postpone the construction of this partition to the end of the proof and explain first how it is used to get \eqref{simplified large deviation lemma}. We split the sum according to the partition and estimate
\[
\E\2\absB{\sum_x b_x X_x}^{2 \mu}\,=\, \E\2\absB{\sum_{k=1}^K\sum_{x \in I_k} b_x X_x}^{ 2\mu}\,\le\, K^{2\mu }\max_{k=1}^K\,\E\2\absB{\sum_{x \in I_k} b_x X_x}^{2 \mu}.
\]
By  the bound \hyperref[P1:bound on size of partition]{(P1)}  on the size of the partition and by choosing $\eps$ sufficiently small, it remains to show 
\bels{remains to show for Ik sum}{
\E\2\absB{\msp{-5}\sum_{\;x \in I_k}\! b_x\2 X_x}^{2 \mu}\,\le\, C(\mu)
\,.
}
For an independent sequence $(X_x)_{x \in I_k}$  satisfying the assumption (ii), the moment bound \eqref{remains to show for Ik sum} would be a simple consequence of  the Marcinkiewicz-Zygmund inequality. Therefore, \eqref{remains to show for Ik sum} follows from 
\bels{replace X by independent X family}{
\absbb{\,\E\2\absB{\msp{-5}\sum_{\;x \in I_k} \!b_x  X_x}^{2\mu}\!-\2\E\2\absB{\msp{-5}\sum_{\;x \in I_k}\! b_x \wt X_x}^{2\mu}\,}
\,\le\, 
\frac{C(\mu, \nu)}{N^{\1\nu}}
\,,
}
for all $\mu,  \nu \in \N$, 
where $\wt X=(\wt X_x)_{x}$ is an independent family of random variables, which is also independent of $X$ and $b$ and has the same marginal 
 distributions as $X$.

To show \eqref{replace X by independent X family} we expand the powers on the left hand side and use the independence of $b$ from $X$ and $\wt X$ as well as the upper bound $\abs{b_x}\le 1$, 
\bels{difference between X and wt X expectation}{
\msp{-15}
\text{l.h.s. of \eqref{replace X by independent X family}}
\leq\,
N^{2 \mu}\msp{-10}\max_{x_1, \dots , x_{2\mu}}\!
\absB{\1
\E\2
X_{x_1}\dots X_{x_{\mu}}\ol{X_{x_{\mu+1}}\dots X_{x_{2\mu}}\msp{-10}}\msp{10}
-\,
\E \2 \wt X_{x_1}\dots \wt X_{x_{\mu}}\ol{\wt X_{x_{\mu+1}}\dots \wt X_{x_{2\mu}}\msp{-10}}\msp{10}
}\,,
}
where the maximum is over all $\xi=(x_1, \dots,x_{2 \mu}) \in I_k^{{2 \mu}}$. For such a $\xi$ let $\xi_1, \dots, \xi_R\in I_k$   denote the indices appearing within $\xi$,  clearly $R\le 2\mu$. 
 Let furthermore the non-negative  integers $\mu_1,\dots,\mu_R$ and $\wt \mu_1,\dots,\wt\mu_R$ denote the corresponding numbers of appearances within $(x_1, \dots,x_{\mu})$ and $(x_{\mu+1}, \dots,x_{2\mu})$, respectively. 
 Then we can further estimate the term inside the maximum on the right hand side of \eqref{difference between X and wt X expectation} corresponding to $\xi$ by using the telescopic sum, 
\bels{telescopic sum of X and wt X}{
\msp{-10}\E\2\prod_{r=1}^RX_{\xi_r}^{\mu_r}\ol{X}_{\xi_r}^{\wt{\mu}_r}
-\,
\E\2\prod_{r=1}^R\wt X_{\xi_r}^{\mu_r}\ol{\wt X}_{\xi_r}^{\wt{\mu}_r}
\,=\, 
\sum_{r=1}^{R-1}
\rm{Cov}\pbb{ X_{\xi_r}^{\mu_r}\ol{X}_{\xi_r}^{\wt{\mu}_r} , \msp{-5}\prod_{\;s=r+1}^R
\msp{-10}
X_{\xi_s}^{\mu_s}\ol{X}_{\xi_s}^{\wt{\mu}_s} }
\prod_{s=1}^{r-1}\E\2X_{\xi_s}^{\mu_s}\ol{X}_{\xi_s}^{\wt{\mu}_s} 
\,.
}
To bound the covariance in \eqref{telescopic sum of X and wt X} we use $\abs{X_x}\le \sqrt{N}$ and \eqref{correlation decay in large deviation}
 with a sufficient large $\nu$   in combination with the estimate on the distance  \hyperref[P2:distance between indices in Ik]{(P2)}  between indices within one element $I_k$ of the partition. The claim \eqref{replace X by independent X family} follows.

We will now  inductively  construct the partition $I_1, \dots,I_K$ with the properties (i) and (ii) above. 
 Suppose  that the disjoint sets  $I_1, \dots,I_{k}$ have already been constructed. Then we pick an arbitrary $x_0 \in J_0:=\bb{X}\setminus (I_1\cup \dots \cup I_{k})$. Next we pick $x_1 \in J_1:=J_0 \setminus B_{N^{\eps}}(x_0)$, then $x_2 \in J_2:=J_1 \setminus B_{N^{\eps}}(x_1)$ and so on. The process stops at  some  step $L$ when $J_{L+1}$ is empty and we set $I_{k+1}:=\{x_0, \dots, x_{L}\}$. 
By construction,  property (ii) is satisfied for all elements $I_k$ of the partition. We verify the upper bound (i) on the number $K$ of such elements. For every $k$ we have
\bels{inclusion on step k}{
\bb{X}\setminus \bigl(\cup_{l=1}^{k-1}  I_l\2\bigr) \,\subseteq \, {\textstyle \bigcup_{x \1\in\1 I_k}}B_{N^{\eps}}(x)\,,
}
because otherwise another element of $I_k$ would have been chosen in the construction. The inclusion 
\eqref{inclusion on step k} implies 
\vspace{-0.6cm}
\bels{bound on size of Ik}{
N- \sum_{\2l\1=1}^{\1k-1}\abs{\1I_l\1}
\,\le\, 
\abs{\1I_k}\max_{x \in I_k}\abs{B_{N^{\eps}}(x)}\,\le\, \abs{I_k}N^{\eps P},
}
where we used \eqref{polynomial ball growth} for the last inequality. 
In particular, \eqref{bound on size of Ik} provides a lower bound on the size of $I_k$ which we use to  obtain that 
\[
N-\2\sum_{l=1}^{k}\2\abs{\1I_l\1}
\,\le\,
\bigl(\21- N^{-\eps P}\,\bigr)\pbb{N-\sum_{l=1}^{k-1}\abs{\1I_l\1}\,}
\,.
\]
Since $I_K$ contains at least one element, we infer by induction that
\[
1\,\le\, N-\sum_{l=1}^{K-1}\abs{I_l}
\,\le\,
(1-N^{-\eps P})^{K-1}N
\,\le\, 
N\ee^{-(K-1) N^{-\eps P}}. 
\]
We solve for $K$ and thus see that (i) holds true. This finishes the proof of Lemma~\ref{lmm:Linear large deviation}.
\end{Proof}

\begin{lemma}[Quadratic large deviation]
\label{lmm:Quadratic large deviation} Let $X=(X_x)_{x\in \bb{X}}$, $Y=(Y_x)_{x\in \bb{X}}$, $b=(\1b_{xy})_{x,y\in \bb{X}}$ be families of random variables that satisfy the following assumptions: 

\begin{itemize}
\titem{i} 
The families $X$ and $Y$ are centered, $\E\2X_x=\E\2Y_x=0$. 
\titem{ii} 
Both $X$ and $Y$ have uniformly bounded moments:  There is a sequence of constants $\ul{\beta}_1$ such that $ \E \2\abs{X_x}^\mu+\E \2\abs{Y_x}^\mu \,\le\, \beta_1(\mu) $, for all $ \mu \in \N $ and all $x\in \bb{X} $. 
\titem{iii}
The correlations within $ X$ and $ Y $ decay fast: There is a sequence $ \ul{\beta}_2 $ of  constants, s.t. for all $\eps > 0 $, $A,B\subseteq\bb{X}$, with $d(A,B)\ge N^{\eps}$, and smooth functions $\phi:\C^{\abs{A}} \to \C$, 
$\psi:\C^{\abs{B}} \to \C$,
\bels{correlation decay in quadratic large deviation}{
\max_{Z,Q \1\in\1\{X,Y\} }\abs{\1\rm{Cov}(\phi(Z_A),\psi(Q_B))}\,\le\, \beta_2(\eps,\nu)\frac{\norm{\nabla \phi}_\infty\norm{\nabla \psi}_\infty}{N^{\1\nu}\!}
\,,\qquad \nu \in \N\,,
}
where $Z_A:=(Z_x)_{x \in A}$ and $d(A,B)$ are defined as in Lemma~\ref{lmm:Linear large deviation}. %
\titem{iv} 
The correlations between $(X,Y)$ and $b$ are asymptotically small: There is a sequence of positive constants  $ \ul{\beta}_3$   such that for alll smooth functions $\phi:\C^{N^2} \to \C$, $\psi:\C^{2N} \to \C $, we have:
\bels{b uncorrelated to X ynd Y}{
\abs{\1\rm{Cov}(\1\phi(\1b\1),\psi(X,Y))}\,\le\,\beta_3(\nu)\2\norm{\nabla \phi}_\infty\norm{\nabla \psi}_\infty N^{-\nu},\qquad \nu \in \N
\,.
}
\end{itemize}
Then the following large deviation estimate holds for every $\nu \in \N$:
\bels{Double sum large deviation}{
\absB{\,\sum_{\2x,\1y}\2b_{xy}\1(X_xY_y-\E\2X_xY_y)}
\;\prec\; 
\pB{\,\sum_{\2x,\1y}\,\abs{\1b_{xy}}^2}^{\!1/2}\msp{-8}+\frac{1}{N^\nu}\,
.
}
\end{lemma}
\begin{Proof} 
We use Convention~\ref{con:Comparison relation} such that $\varphi\lesssim \psi$ means  $\varphi\le C \psi$ for a constant $C$, depending only on $\wt{\scr{P}}:=(\beta_1,\beta_2,\beta_3,P)$ (cf. \eqref{polynomial ball growth}).
The proof of Lemma~\ref{lmm:Quadratic large deviation} follows a similar strategy as the proof of Lemma~\ref{lmm:Linear large deviation}. Exactly as  in Step~1  
of the proof of Lemma~\ref{lmm:Linear large deviation}, we introduce new families of centered random variables
\[
\wt X_x \,:=\, (1-\E)[X_x\2\theta(N^{-1}\abs{X_x}^2)]
\,,
\qquad
\wt Y_x \,:=\, (1-\E)[Y_x\2\theta(N^{-1}\abs{Y_x}^2)]\,,
\]
and rescaled coefficients
\[
\wt{b}_{xy}\,:=\, 
\frac{b_{xy}}{(\2\sum_{u,v}\abs{\1b_{uv}}\1)^{1/2}+N^{-\wt \nu}}\,.
\]
In this way we reduce the proof of \eqref{Double sum large deviation} to showing the moment bound 
\bels{simplified quadratic large deviation}{
\E\2\absB{\sum_{x,y}b_{xy}\,(X_xY_y-\E\2X_xY_y)}^{2\mu}\,\le\,C(\mu,\delta)N^{\delta}\,,\qquad \delta >0\,,
}
for random variables $X,Y$ and $b$ that satisfy all assumptions of Lemma~\ref{lmm:Quadratic large deviation}, 
 and  the additional bounds 
\bels{additional bounds quadratic}{
\max_x\2 \abs{X_x}\,\le\, \sqrt{N}
\,,\qquad
\max_x\2 \abs{Y_x}\,\le\, \sqrt{N}
\,,\qquad 
\sum_{x,y} \,\abs{\1b_{xy}}^2\,\le\, 1\,.
}
 Following Step~2 of the proof of Lemma~\ref{lmm:Linear large deviation} and 
using \eqref{b uncorrelated to X ynd Y}  
we may also assume that $b$ is independent of $(X,Y)$. 

To show \eqref{simplified quadratic large deviation} we fix $\eps>0$ and choose the partition $I_1, \dots, I_K$ from Step~3 of the proof of Lemma~\ref{lmm:Linear large deviation} of the index set $\bb{X}$. In particular,  the properties \hyperref[P1:bound on size of partition]{(P1)} and \hyperref[P2:distance between indices in Ik]{(P2)} introduced in that proof  are satisfied.  We  split the sums over $x$ and $y$ in \eqref{simplified quadratic large deviation} according to this partition and estimate 
\[
\text{l.h.s. of \eqref{simplified quadratic large deviation}}
\,\le\, 
K^{4 \mu}\,\max_{k,l=1}^K\,\E\,\absB{\!\sum_{\;x \in I_k,y \1\in\1 I_l}\msp{-10}b_{xy}\2(X_xY_y-\E\2X_xY_y)\,}^{2\mu}
.
\]
By choosing $\eps$ sufficiently small and using  \hyperref[P1:bound on size of partition]{(P1)}  it suffices to show that for any fixed $k,l=1, \dots,K$ we have the moment bound 
\bels{restriction to Ik and Il}{
\E\,\absB{\msp{-5}\sum_{\;x \in I_k,y \in I_l}\msp{-10}b_{xy}\2(X_xY_y-\E\2X_xY_y)\,}^{2\mu}
\,\le\, 
C(\mu,\delta)\2N^{\delta}, \qquad \delta>0
\,.
}

For any $x \in I_k$ and $y \in I_l$ we introduce the relation
\[
x \,\Join\, y\qquad \text{whenever} \qquad   d(x,y) \2\le\2 {\textstyle \frac{\2N^{\eps}\!}{3}}\,.
\]
If $d(x,y)>\frac{N^{\eps}}{3}$, we correspondingly write $x \not \Join y$. Since the distances of indices within the set $I_k$ are bounded from below by 
 $N^{\eps}$ (c.f. the property \hyperref[P2:distance between indices in Ik]{(P2)}), 
we see that for every $x \in I_k$ there exists at most one $y \in I_l$  such that $x \Join y$ and the other way around. We set 
\begin{equation*}
\begin{split}
\iota(x)\,:=\,
\begin{cases}
1 & \text{ if there is } \; \wt y \in I_l \;\text{ s.t. }\; x \Join \wt y\,;
\\
0 & \text{ otherwise}\,,
\end{cases}
\,,\quad 
\iota(y)\,:=\,
\begin{cases}
1 & \text{ if there is } \; \wt x \in I_k \;\text{ s.t. }\; \wt x \Join y\,;
\\
0 & \text{ otherwise}\,,
\end{cases}
\end{split}
\end{equation*}
for any $x \in I_k$ and $y \in I_l$. Note that if $k=l$, then $\iota(x)=1$ for all $x \in I_k$. Furthermore,  let us define 
\[
\scr{S}\,:=\, 
\setB{\2(S,T): S\subseteq I_k\,,\;T \subseteq I_l\,,\;\text{such that } x \not \Join y  \text{ for all } x \in S \,, \; y \in T\2}
\,,
\]
the pairs of subsets with a distance of at least $\frac{N^{\eps}}{3}$. Inspired by Appendix B of \cite{EKYY3} we use the partition of unity
\bels{partition of unity}{
1\,=\,\frac{\,\sigma_{xy}\!}{\abs{\scr{S}}}\msp{-8}\sum_{\;(S,T) \in\1 \scr{S}}\msp{-10}\bbm{1}(x \in S) \bbm{1}(y \in T) \,, \qquad x \in I_k\,,\; y \in I_l\,, \; x \not \Join y\,,
}
where we introduced the numbers  $\sigma_{xy}$ to be $4,6,6$ and $9$ in the case when $(\iota(x),\iota(y))$ is $(0,0),(0,1),(1,0)$ and $(1,1)$, respectively.  
We split the sum in \eqref{restriction to Ik and Il} into a sum over pairs $(x,y)$ with $x \Join y$ and $x \not \Join y$. Afterwards we use \eqref{partition of unity} and find
\[
\sum_{x \in I_k,y \in I_l}\msp{-5}b_{xy}\,(X_xY_y-\E\2X_xY_y)\,=\, 
U + \frac{1}{\abs{\scr{S}}}\sum_{(S,T) \in \scr{S}}\msp{-5}V(S,T)\,,
\]
with the short hand notation
\[
U\,:=\msp{-5} \sum_{x \1\in\1 I_k,\2 y\1\in\1 I_l}%
\msp{-15}\bbm{1}( x \Join y)\,
b_{xy}\,(X_xY_y-\E\2X_xY_y)\,,\qquad 
V(S,T)\,:=\msp{-5} \sum_{x\1\in\1 S,\2 y \1\in\1 T}\msp{-7}\sigma_{xy}\,b_{xy}\,(X_xY_y-\E\2X_xY_y)\,.
\]
Thus, proving \eqref{restriction to Ik and Il} reduces to showing  for any pair of index sets $(S,T) \in \scr{S}$ that 
\[
\E\2\abs{\1U}^{2\mu}\!+ \E\2\abs{V(S,T)}^{2\mu}\,\le\, C(\mu).
\]
The moment bound on $U$ can be seen with exactly the same argument as \eqref{simplified large deviation lemma} in the proof of Lemma~\ref{lmm:Linear large deviation} since the family of centered random variables $X_xY_y-\E\2X_xY_y$  in this sum are almost independent. 
 The moment bound on $V(S,T)$  
 follows by comparing the moments of $V(S,T)$ with the moments of 
\[
\wt V(S,T)\,:=\msp{-5} \sum_{x\1\in\1 S,\2 y \1\in\1 T}\msp{-5}\sigma_{xy}\,b_{xy}\2\wt X_x \wt Y_y\,,
\]
where $\wt X_S =(\wt X_x)_{x \in S}$ and $\wt Y_T =(\wt Y_x)_{x \in T}$ are independent families of random variables, which are independent of $(b,X,Y)$ as well, with the same marginal distributions as $X_S=(X_x)_{x \in S}$ and $Y_T =(Y_x)_{x \in T}$, respectively. 
As  the result of this comparison,
$
\,\abs{\,\E\2 \abs{V(S,T)}^{2 \mu}- \E\2 \abs{\wt V(S,T)}^{2 \mu}\2} 
\le 
C(\mu,\nu)\2N^{-\nu},
$ 
because $X_S$ and $Y_T$ are essentially uncorrelated since $d(S,T) \ge \frac{N^{\eps}}{3}$ and because the families $X_S$ and $Y_T$ themselves are already essentially uncorrelated (cf. the property \hyperref[P2:distance between indices in Ik]{(P2)} from the proof of Lemma~\ref{lmm:Linear large deviation}).
Finally, the moments of $\wt V(S,T)$ satisfy the necessary bound by the Marcinkiewicz-Zygmund inequality as in the 
proof of Lemma~\ref{lmm:Linear large deviation}.  
The details are left to the reader. 
\end{Proof}


\end{document}